\documentclass[article, 12pt]{amsart}
\usepackage{amscd}
\usepackage{bbm}
\usepackage{dsfont}
\usepackage{enumerate}
\usepackage{latexsym}
\usepackage{srcltx}
\usepackage{amsfonts}
\usepackage{amssymb}
\usepackage{datetime}
\usepackage{setspace}
\usepackage{enumerate}
\usepackage{amsthm}
\usepackage{graphicx}
\usepackage{epstopdf}

\theoremstyle{definition}

\input epsf
\begin{document}

\title{ Conditional expectations and interpolation of linear operators\\
 on ordered ideals between $L^1(0.1)$ and $L^\infty(0,1)$.}

\author{A. A. Mekler}

\address{Wilsdruffer Str.7, 01067 Dresden, BRD/Germany}

\email{alexandre.mekler@mail.ru}


\date{\today}


\maketitle

\markboth{A. A. Mekler}{Conditional expectations and interpolation of linear operators }

\newpage

I am grateful to my friend professor Arkady Kitover for stimulating discussions, valuable advices and his huge work of translation the book into English.

\newpage

\centerline{\textbf{CONTENTS}}

\bigskip

\textbf{Introduction.}......................................................................................4

\bigskip

\textbf{Chapter 0. Definitions, notations, and some prerequisites.}........8

\bigskip

\textbf{Chapter 1. Principal symmetric and principal majorant ideals
generated by weakly regular and regular functions.}.............19

\bigskip

\textbf{Chapter 2. Two absolute constants for rearrangements of intervals.}.................................................................................................35

\bigskip

\textbf{Chapter 3. Double stochastic projections and interpolation property
 of ideals between $L^1$ and $L^\infty$}................................................39

 \bigskip

 \textbf{Chapter 4. On averaging principal symmetric ideals by countable partitions.}.............................................................................51

\bigskip

\textbf{Chapter 5. On symmetric ideals generated by averaging}
\textbf{ of principal symmetric ideals over countable partitions.}.....................60

\bigskip

\textbf{Chapter 6. Verifying and universal $\sigma$-subalgebras.}......................72

\bigskip

 \textbf {Chapter 7. Independent complement of an interval partition.}..........................................................................................................79

 \bigskip

 \textbf {Chapter 8. Verifying and universal complemented}

  \textbf{$\sigma$-subalgebras.}......................................................................................92

 \bigskip

 \textbf{Appendix}..............................................................................................98

 \bigskip

 \textbf{References}............................................................................................103

 \bigskip

 \textbf{Index}.....................................................................................................106

\newpage

\centerline{\textbf{Introduction}}

\bigskip

Symmetric \footnote{We will use the term ``symmetric'' rather then ``rearrangement invariant'' because in the western literature the last term means that the ideal space is a Banach lattice with monotone norm and Fatou properties (Properties (B) and (C) in~\cite{KA}).} Banach ideal spaces of measurable functions were introduced in Functional Analysis at the beginning of nineteen sixties by G.G. Lorentz (see~\cite{Lo}), who was inspired by the well known Hardy-Littlewood inequalities and by the results concerning Bernstein's polynomials. The study of these spaces was then continued by T. Shimogaki in~\cite{Sh}. During the next decade appeared three monographs devoted to the subject: by W.A.J. Luxemburg~\cite{Lu}, by K-M Chong and N.M.Rice~\cite{CR}, and by J. Lindenstrauss and L. Tzafriri~\cite{LT}. It is worth emphasizing that all these authors considered ideal symmetric \textbf{Banach} linear subspaces of $L^0$.

A new powerful impulse for studying these objects was given by the seminal paper of A.P. Calderon~\cite{Ca}. In that paper Calderon proved an important theorem about real interpolation of linear operators acting on vector spaces intermediate between $L^1$ and $L^\infty$. Recall that such a space $V,\ L^\infty \subseteq V \subseteq L^1$, is called a $(L^1,L^\infty)$ interpolation space if for \textbf{every admissible linear operator $T$} (i.e. $TL^\infty \subseteq L^\infty$ and $TL^1 \subseteq L^1$) we also have $TV \subseteq V$. Without going here into technical details we can say that Calderon proved that a vector space $V$, intermediate between $L^1$ and $L^\infty$, is an interpolation space if and only if it satisfies some \textit{majorant} property (in the sense of Hardy - Littlewood).

The theory of real linear interpolation was then successfully developed by many mathematicians, among them C. Bennet, D.W. Boyd, S.G. Krein, B.S. Mityagin, V.I. Ovchinnikov, J. Peetre, Y.I. Petunin, R. Sharply, A.A. Sedaev, E.M. Semenov, V.A. Schestakov. See also the bibliography in~\cite{KPS}.

In his paper of 1965~\cite{Ry} J.V. Ryff studied double stochastic operators and their orbits. Coupled with Calderon's theorem his approach shows that when verifying whether an intermediate vector space $V$ is an interpolation space it is sufficient to check that $TV \subseteq V$ not for every admissible operator but only for $T$ from the much smaller class of all \textbf{double stochastic operators}. Thus in order to $V$ be an $(L^1, L^\infty)$-interpolation space it is sufficient to prove that with each function $f$ $V$ must contain its orbit $\{Tf: T$ is a double stochastic operator $\}$. This is exactly the Calderon's majorant property which is equivalent to $(L^1, L^\infty)$-interpolation property.

Because the class of all double stochastic operators includes the composition operators generated by measure preserving transformations, an interpolation space must be symmetric. But the symmetric property \emph{is not sufficient}: there are Banach symmetric ideal spaces that are not interpolation spaces, ~\cite{Ru}.
 \footnote{On the other hand the rearrangement invariant intermediate Banach ideal spaces with monotone norm and Fatou properties always have the interpolation property (see e.g.~\cite{BS}).}

A considerable part of the results discussed in this paper is based on the idea of further restriction of the class of ``verifying'' operators. In particular, in Chapter 3 it is proved that an intermediate order ideal $V$ is an interpolation space if and only if $PV \subseteq V$ for every \textbf{double stochastic projection} $P$. \footnote{It is immediate to see that every interpolation space must be an order ideal in $L^1$ and therefore this result is an improvement of the Calderon-Ryff theorem.}

Let us emphasize that in this paper our measure space is $(I, \Lambda, \lambda)$ where $I$ is the interval $(0,1]$, $\lambda$ is the Lebesgue measure, and  $\Lambda$ is the $\sigma$-algebra of all Lebesgue measurable subsets of $I$. In this case every double stochastic projector is an operator of conditional mathematical expectation $E(\cdot|\mathcal{A})$ (we will also call it an \textbf{averaging} operator corresponding to a $\sigma$-subalgebra $\mathcal{A}$ of $\Lambda$). We will say that $\mathcal{A}$ \textit{averages} an ideal $V$ if $E(V|\mathcal{A})\subseteq V.$\\

The current monograph is mostly devoted to the study of averaging operators in the gap between vector order symmetric and interpolation ideals.\\

If a $\sigma$-subalgebra is purely discrete, i.e. generated by a finite or countable set of atoms, we will call it a finite, respectively countable, \textit{partition}. Because the objects we study are invariant under measure preserving transformations we can without loss of generality consider \textit{interval partitions}, i.e. such partitions whose atoms are pairwise disjoint intervals $(t_{n+1},t_n],\ t_n\downarrow 0$. Especially important for us among interval partitions will be \textit{monotonic} partitions. An interval partition is called monotonic if the lengths of the corresponding intervals do not increase when $I$ is passed in the direction from 1 to 0. We prove in Chapter 2 that every interval partition is equivalent to a subsequence of the (unique) monotonic partition that is obtained by ordering the intervals of the original partition in non-increasing (by length) order. And the smallest constant of equivalence is the golden ratio.

In Chapter 6 we introduce two disjoint classes of $\sigma$-subalgebras: verifying and universal $\sigma$-subalgebras. According to our definitions a universal $\sigma$-subalgebra averages any symmetric ideal, while a verifying one averages only the interpolation ideals.
Thus, if a vector space $V$ is a symmetric ideal, then to check whether it is an interpolation space it is enough to check that $E(V|\mathcal{A}) \subseteq V$ just for one verifying $\sigma$-subalgebra $\mathcal{A}$.

In Chapter 7 we prove a general result that \textit{if the independent complement to a countable partition averages a symmetric ideal then the partition itself averages this ideal}.

\textit{Moreover, if e.g. the symmetric ideal is contained in the space $Llog^+ L$ then the converse is also true}.

This result allows us to answer in Chapter 8 some questions raised by the author in~\cite{Me1}. In particular, by using this result we are able to establish a complete classification of verifying and universal $\sigma$-subalgebras in the class of all $\sigma$-subalgebras that have an independent complement in $\Lambda$ (for brevity we will call them \textit{complemented} $\sigma$-subalgebras). Notice that all at most countable partitions belong to this class.

 In Chapter 1 we introduce the notion of a \textit{principal symmetric ideal} $\mathcal{N}_f$ i.e the smallest (by inclusion) symmetric ideal that contains a function $f,\ f \in L^1(I)$. The notion of principal symmetric ideal can be actually considered as a generalization of the notion of a random variable. A function $f$ is called \textit{regular} if the corresponding principal ideal is an interpolation space. It is worth noticing that for a regular $f$ the ideal $\mathcal{N}_f$ coincides as a set with well known Marcinkiewicz space $\mathcal{M}_\psi$ where the function $\psi$ is  non-decreasing, concave down, continuous in 0, $\psi(0)=0$, and its derivative coincides with the non-increasing equimeasurable rearrangement of $|f|$ on $I$.

We consider the question of possibility to introduce a norm equivalent to a Banach lattice norm on a principal symmetric ideal. We also study the behaviour of regular functions under the action of operators of compression and dilation on $I$. In the same Chapter we prove that regular functions can be majorized by power functions and characterise the BMO space in terms of regular and weakly regular functions.

In Chapter 4 we provide a solution of the following problem: to state in terms of a function $f \in L^1(I)$ and a countable partition $\mathcal{F}$ necessary and sufficient conditions for $\mathcal{F}$ to average the principal symmetric ideal $\mathcal{N}_f$. If it is the case we call $f$ a $\mathcal{F}$-\textit{regular} function.

In Chapter 5 we study the properties of $\mathcal{F}$-regular functions, in particular, those of these properties that are analogous to properties of regular functions. We also study necessary and sufficient conditions for the image of a principal symmetric ideal under the action of an averaging operator to be itself a principal symmetric ideal.

The results contained in Chapters 7 and 8 are new and were not published yet. Other results were obtained by the author during the previous years, starting from 1974, and paper~\cite{Me1} can be considered as a prequel to the current monograph. Actually the goal of the monograph is to present in a unifying manner the main results of the author devoted to the topic named in the title.

Chapter 0 contains terminology and some information used in the sequel. At the beginning of each Chapter there are references and annotations pertinent to the contents of this Chapter. The general bibliography is located customary at the end of the monograph.\\

 Appendix contains the proof of a result of G. Ya. Lozanovsky and Yu. A. Abramovich that initiated the author's research represented in this monograph.\\

If a statement is given without a proof, reference, or clarification, that means that it is either trivial or sufficiently well known.\\

\newpage

\bigskip

\textbf{Chapter 0. Definitions, notations, and some prerequisites.}

\bigskip

We consider the measure space $(I, \Lambda, \lambda)$ where $I$ is the interval $(0,1]$, $\Lambda$ is the $\sigma$-algebra of all Lebesgue measurable subsets of $I$, and $\lambda$ is the probability Lebesgue measure.

A measurable map $\pi : I \to I$ is called a (metric) \textit{endomorphism} if $\lambda(\pi^{-1}(A))=\lambda(A)$, $A \in \Lambda$. If, additionally there is a measurable map $\pi^{-1}: I \to I$ such that $\pi(\pi^{-1}(x)=x$ a.e. on $I$ then $\pi$ is called an \textit{automorphism}. As usual the space of all (classes of a.e. coinciding) measurable functions is denoted by $L^0(I, \Lambda, \lambda)$= $L^0(I)$. Endowed with the partial order
$f \leq g \Leftrightarrow f(x) \leq g(x)$ a.e. the space $L^0(I)$ becomes a Dedekind complete vector lattice. The characteristic function (the indicator) of a measurable set $A$ is denoted by $1_A$, with $\mathbf{1}$ meaning $1_I$.

Two nonnegative functions $f, g \in L^0(I)$ are called \textit{equivalent} ($f \simeq g$) if there is a constant $C \geq 1$ such that $C^{-1}g \leq f \leq Cg$. Two functions $f, g \in L^0(I,\Lambda, \lambda)$ are called \textit{equimeasurable} (notation: $f \sim g$) if for any $\alpha > 0$ we have $\lambda(\{t: f(t) > \alpha\}) = \lambda(\{t: g(t) > \alpha\})$. Clearly, if $\pi$ is an endomorphism then the functions $f$ and $f\circ\pi$ are equimeasurable.
It is shown in~\cite{CR} that for any $f \in L^0(I, \Lambda, \lambda)$ there is the unique, continuous from the right, non-increasing function equimesurable with $|f|$. It is called $f^\star$ and named the \textit{non-decreasing rearrangement} of $f$. It is also proved in~\cite{CR} that there is an endomorphism $\pi$ such that $|f|=f^\star \circ \pi$.

As usual, we denote by $L^1(I,\Lambda, \lambda):= L^1(I)$ (respectively,$L^\infty(I,\Lambda, \lambda):= L^\infty(I)$) the vector subspaces of $L^0(I)$ that consist of functions such that $\int \limits_0^1 |f| d\lambda < \infty $ (respectively, $ess \; sup |f| < \infty$).

For $f \in L^1(I)$ the function $f^{\star \star}$ is defined as
$$f^{\star \star}(t) =\frac{1}{t}\int \limits_0^t f^\star d\lambda, \; t \in I.$$
Clearly, $f^\star \leq f^{\star \star}$.  The notation $g \prec f$ means that $g^{\star \star} \leq f^{\star \star}$; the notation $g \preceq f$ means that $g\prec f$ and $\int_0^1g d\lambda=\int_0^1f d\lambda.$

For any positive real number $s$ we introduce the compression-dilation operator ($c/d$-operator for brevity), $\rho_s: L^1(I) \rightarrow L^1(I)$ by the formula
$$(\rho_s f)(t)=
\left\{
  \begin{array}{ll}
  $f(st)$ , & \hbox{if $ 0<st \leq 1$;} \\
    0, & \hbox{if $st>1$,}
  \end{array}
\right.
f \in L^1(I), t \in (0,1]. \eqno{(0.1)}
$$

The statement of the next lemma can be easily checked. \bigskip

\textbf{Lemma 0.1}. For any $f,g \in L^0(I)$ we have
$$\begin{cases}
(\rho_2f^*)(t)\leq f^*(t)\leq(\rho_{\frac{1}{2}}f^*)(t)\leq \Big(f^*(t)+f^*(1-t)\Big)^*,\ t\in I;\\
0\leq f^*(t)\leq f^{**}(t)=(f{^{**}})^*(t)\geq t^{-1}\int_0^tfd\lambda,\ t\in I;\\
f^*\leq \rho_s f^*,\ \rho_s f^{**}\leq s^{-1}f^{**},\ s\in I;\\
(f+g)^*\leq \rho_sf^*+\rho_{1-s}g^*,\ 0<s<1;
\end{cases}\eqno(0.2)$$
$\Box$\\

\textbf{Lemma 0.2} For any $f \in L^0(I)$ and any $s \in I$
$$(\rho_sf)^* \leq \rho_sf^*.\eqno(0.3)$$
\begin{proof}For any $f,\ f\in L^0(I),$ and any $\alpha>0$ let
$A_f(\alpha):=\{t\in I:f(t)>\alpha\}$. Then for any $u\in I$ we have
$$A_{\rho_u f}(\alpha)= u^{-1}A_f(\alpha).$$
Indeed,  applying the following equivalencies
$$s\in u^{-1}A_f(\alpha)\Leftrightarrow us\in A_f(\alpha)\Leftrightarrow f(us)>\alpha\Leftrightarrow\rho_uf(s)>\alpha\Leftrightarrow s\in A_{\rho_uf}(\alpha)$$
to $(\rho_uf)^*$ we notice that from $(0.2)$  and properties of Lebesgue measure it follows that
$$\lambda\{\rho_uf>\alpha\}=\lambda\{ A_{\rho_uf}(\alpha)\}=\lambda\{ u^{-1}A_f(\alpha)\}\leq u^{-1}\lambda \{A_{f^*}(\alpha)\}=\lambda\{\rho_uf^*>\alpha\}.$$

\end{proof}

We will call a measurable function $f$ an \textit{elementary} function if for any $n \in \mathds{N}$ the set $f((0,1]) \cap (-n,n)$ is finite (or empty). In particular, the set of values of such a function is at most countable. It is immediate to see that for any measurable $f$ and for any $\varepsilon >0$ there is an elementary function $g$ such that $|f(x)-g(x)|<\varepsilon, x\in(0,1]$. If $f$ and $g$ are equimeasurable elementary functions then there is an automorphism $\pi$ such that $g=f\circ \pi$. In particular, if $f$ is an elementary function then there is such an automorphism $\pi$ that $f^\star = |f|\circ \pi$.

A subset $\mathcal{A}$ of $\Lambda$ is called a $\sigma$-subalgebra of $\Lambda$ (just $\sigma$-subalgebra, if it does not cause any ambiguity) if
\begin{enumerate}
  \item $A_n \in \mathcal{A}, n \in \mathds{N} \Rightarrow \bigcup \limits_{n=1}^\infty A_n \in \mathcal{A}$ and
$\bigcap \limits_{n=1}^\infty A_n \in \mathcal{A}$.
  \item $A \in \mathcal{A} \Rightarrow I \setminus A \in \mathcal{A}$.
  \item If $A \in \Lambda$ and $\lambda(A)=0$ then $A \in \mathcal{A}$.
\end{enumerate}

If $a \in \mathcal{A}$ then we call $a$ an \textit{atom} in $\mathcal{A}$ if for any $b\in \mathcal{A}$, $b \subset a$ either $\lambda(b)=0$ or $\lambda(a \setminus b)=0$. A $\sigma$-subalgebra $\mathcal{A}$ is called \textit{continuous} if it contains no atoms and is called \textit{discrete} if $\lambda((0,1] \setminus d)=0$ where $d$ is the set of all atoms in $\mathcal{A}$. We say that $\mathcal{A}$ is of \textit{mixed type} if it is neither discrete nor continuous.

For a subset $G$ of $\Lambda$ we denote $\sigma(G)$ the smallest $\sigma$-subalgebra in $\Lambda$ that contains $G$. We call $\sigma(G)$ the $\sigma$-subalgebra \textit{generated} by $G$.

Two $\sigma$-subalgebras $\mathcal{H}$ and $\mathcal{G}$ are called \textit{equimeasurable} if there is an automorphism $\pi$ such that $\pi\mathcal{H} = \mathcal{G}$. In this case we will write $\mathcal{H} \sim \mathcal{G}$.

If $\mathcal{H}$ and $\mathcal{G}$ are two $\sigma$-subalgebras and $\mathcal{H} \subseteq \mathcal{G}$ we will say that $\mathcal{H}$ is \textit{coarser} than $\mathcal{G}$  or that $\mathcal{G}$ is \textit{finer} than $\mathcal{H}$.

Let $\mathcal{A}$ be a $\sigma$-subalgebra. We say that $\mathcal{A}$ is a \textit{complemented} $\sigma$-subalgebra in $\Lambda$ if there is a $\sigma$-subalgebra $\mathcal{A}^\perp$ (called the \textit{idependent complement} of $\mathcal{A}$) such that
$$ \sigma(\mathcal{A} \cup \mathcal{A}^\perp) = \Lambda  $$
and
$$ \lambda(D \cap E) = \lambda(D) \times \lambda(E), \; D \in \mathcal{A}, E \in \mathcal{A}^\perp . $$

We refer the reader to~\cite{Ro} for necessary and sufficient conditions for a given $\sigma$-subalgebra of $\Lambda$ to have the independent complement in $\Lambda$. It follows from these conditions that any discrete $\sigma$-subalgebra of $\Lambda$ has the independent complement in it.

Let $\mathcal{F}$ be a discrete $\sigma$-subalgebra and $\{F_n: n \in \mathds{N}\}$ be the countable or finite set of all its atoms. Then we call $\mathcal{F}$ a countable (respectively, finite) \textit{partition}. We say that the partition $\mathcal{F}$ \textit{belongs to} the \textit{stochastic vector} $\vec{\phi}=[\phi_n]$ (notation:\ $\mathcal{F}\in \vec{\phi}$) , where $\phi_n = \lambda(F_n):\ n \in \mathds{N}$. Clearly, $\phi_n>0$ and $\sum \limits_{n=1}^\infty \phi_n =1$.

Let $\vec{\phi}$ and $\vec{\psi}$ be two stochastic vectors We say that $\vec{\phi}$ is finer (coarser) than $\vec{\psi}$ if there are discrete $\sigma$-algebras $\mathcal{A} \in \vec{\phi}$ and $\mathcal{B} \in \vec{\psi}$ such that $\mathcal{A}$ is finer (respectively, coarser) than $\mathcal{B}$.

Notice that every elementary function is measurable in respect to some at most countable partition.

If two partitions $\mathcal{A} \in \vec{\phi}$ and $\mathcal{B} \in \vec{\psi}$ are equimeasurable then there is a permutation $\tau$ of $\mathds{N}$ such that $\tau(\vec{\phi}) = \vec{\psi}$.

We will call a countable partition $\mathcal{B}$ an \textit{interval partition} if its atoms
 are intervals $(b_n, b_{n-1}]$ where $1=b_0 > b_1 > \cdots $ and $b_n \rightarrow 0$. Clearly, $\mathcal{B} \in \vec{\beta} =[\beta_n]$ where $\beta_n = b_n - b_{n-1},\ n \in \mathds{N}$. We will denote the interval partition $\mathcal{B}$ by $(b_n)$. It is immediate to see that $\mathcal{T} = (t_n)$ is finer than
$\mathcal{S} = (s_n)$ if and only if $\{s_n\} \subseteq \{t_n\}$.

Let $\mathcal{T}$ and $\mathcal{S}$ be two interval partitions, We call them equivalent and write $\mathcal{T} \simeq \mathcal{S}$ if there is a constant $C \geq 1$ such that
$C^{-1}t_n \leq s_n \leq Ct_n, n \in \mathds{N}$.

We say that an interval partition $\mathcal{S}$ is a \textit{multiple} of the interval partition $\mathcal{T}$ if there are $n_0 \in \mathds{N}$ and $p \in \mathds{R}_+$ such that $s_n = pt_n, n \geq n_0$.

An interval partition $\mathcal{B}$ is called \textit{monotonic} if $\beta_n \geq \beta_{n+1}, n \in \mathds{N}$. If $\vec{\beta}$ is a stochastic vector we will denote by $\vec{\beta}^\star$ the unique stochastic vector which coordinates are obtained by a permutation of coordinates of $\vec{\beta}$ and are non-increasing. Any countable partition $\mathcal{F}$ with stochastic vector $\vec{\beta}$ is equimeasurable with some monotonic interval partition $\mathcal{B}^\star$ with stochastic vector $\vec{\beta}^\star$.

We will call the points $2^{-n}, n \in \mathds{N}$ \textit{diadic} points. We denote the interval $(2^{-n}, 2^{-n+1}]$ by $D_n$. We call the interval partition $(2^{-m_n})$, where $m_n \in \mathds{N}$ and $0 = m_0 < m_1 < ...$, a \textit{diadic} partition. We denote by $\mathcal{D}$ the diadic partition $(2^{-n})$.

It is not difficult to see that for any $f = f^\star \in L^1$ we can find a $\mathcal{D}$-measurable $g,\ g = g^\star = \sum \limits_{n=1}^\infty a_n 1_{D_n}$ such that
$$ \rho_2 g(t) \leq f(t) \leq \rho_{1/2} g(t), t \in I.    \eqno{(0.4)} $$

To each interval partition $\mathcal{T} = (t_n)$ corresponds its \textit{(right) diadic projection} $\mathcal{T}_{(2)}$ generated by all distinct points $2^{-[log_2 t_n]+1}$, where $[r], [r] \leq r < [r]+1$ is the integer part of $r$. Thus to every countable partition $\mathcal{F}$ corresponds the unique diadic interval partition $\mathcal{F}^\star_{(2)}$.

A \textit{sample} from a countable partition $\mathcal{F}=\sigma(F_n,\ n\geq 1,)$ is the partition
$$\mathcal{F}_{(n_m)}:=\sigma(F_{n_m},\ m\geq 0),\ \textrm{where}\ F_{n_0}:= I\setminus\bigcup_{m\geq 1}F_{n_m}.$$

The following two propositions follow easily from the fact that $\lambda$ is a continuous measure, ~\cite{DS}. \\

\textbf{Proposition 0.3}. Let $A,B\in \Lambda,\ \lambda(A)=\lambda(B)>0,\ \lambda(A\cap B)=0,$ and let $\{A_n\}_{n\geq 1}$  be a partition of the set $A$ into pairwise disjoint subsets of positive measure. Then there is a partition $\{B_n\}_{n\geq 1}$ of the set $B$ into pairwise disjoint subsets of positive measure such that $\lambda(B_n)=\lambda(A_n),\ n\geq 1.$\\
$\Box$\\

\textbf{Proposition 0.4.}
 Let $A\in \Lambda,\ \lambda(A)>0$ and let $\{A_n\}_{n=1}^\infty$ be a partition of $A$ into pairwise disjoint subsets of positive measure. For every $f\in L^1(I)$ there is an $\tilde{f},\ \tilde{f}\sim f$, such that
$$\frac{1}{\lambda(A_n)}\int_{A_n}\tilde{f}d\lambda= \frac{1}{\lambda(A)}\int_Afd\lambda,\ n\geq 1.$$
$\Box$\\

For any $\sigma$-subalgebra $\mathcal{A}$ in $\Lambda$ the \textit{operator of conditional mathematical expectation}, or the \textit{averaging operator} $$E(\cdot|\mathcal{A}):L^1(I)\rightarrow L^1(I,\mathcal{A},\lambda)$$ is defined as Radon-Nikodym derivative $$E(f|\mathcal{A}):=\frac{d\lambda_f(A)}{d\lambda(A)},\ f\in L^1(I),$$
 where the measure $\lambda_f$ on $\mathcal{A}$ is defined by the formula  $\lambda_f(A):=\int_Afd\lambda, A\in \mathcal{A}.$\\
$\Box$\\

\textbf{Proposition 0.5.}

1). \textit{The projection property}.\\
If a $\sigma$-subalgebra $\mathcal{A}$ is finer than the $\sigma$-sublagebra $\mathcal{C}$, then
 $$E\Big(E(\cdot|\mathcal{A})|\mathcal{C})\Big)=E(\cdot|\mathcal{C}).$$

2). \textit{Change of variable}.\\

Let $\pi$ be an automorphism on $\Lambda$. For any $f\in L^1(I)$ we have
$$E(f\circ\pi|\mathcal{A}\circ\pi)=E(f|\mathcal{A})\circ\pi.$$
\\

3). If $\mathcal{F}=\sigma\Big(F_n,\ n\geq 1\Big)$ is an at most countable partition then
$$E(f|\mathcal{F})(t)=\sum_{n\geq 1}\Big(\frac{1}{\lambda(F_n)}\int_{F_n}fd\lambda \Big)\cdot 1_{F_n} (t),\ t\in I.\eqno(0.5)$$
$\Box$\\

Our next proposition follows from Propositions 0.4 and 0.5(2).\\

\textbf{Proposition 0.6.} 1). For any countable partition $\mathcal{F}$ and any function $f\in L^1(I)$ the equality $E(f|\mathcal{F})^*=E(\bar{f}|\mathcal{B})$ holds true where an interval partition $\mathcal{B}\sim \mathcal{F}$ and a function $\bar{f}=(\bar{f})^*\sim f.$\\

2). If countable partition $\mathcal{F}$ is finer than countable partition $\mathcal{G}$ then for any $f\in L^1(I)$ there is an $\tilde{f}\sim f$ such that $$E(f|\mathcal{G})=E(\tilde{f}|\mathcal{F}).\eqno(0.6)$$

3). If countable partitions $\mathcal{F}$ and $\mathcal{G}$ are equimeasurable then for any $f\in L^1(I)$ there is $g\in L^1(I),\ g\sim f$, such that $E(g|\mathcal{G})\sim E(f|\mathcal{F}).$\\
$\Box$\\

In the sequel we will need some notions from the interpolation theory of linear operators on $L^1(I).$\\

Let $T:L^1(I)\rightarrow L^1(I)$ be a bounded linear operator. We will call $T$ an \textit{admissible} operator if $T|_{L^\infty(I)}\subseteq L^\infty(I).$

A positive admissible operator $\mathcal{U}:L^1(I)\rightarrow L^1(I)$ is called \textit{double stochastic} if $\mathcal{U}(\textbf{1})=\mathcal{U}^*(\textbf{1})=\textbf{1},$ where $\mathcal{U}^* : L^\infty(I) \rightarrow L^\infty(I)$ is the operator adjoint to $\mathcal{U}$. \\

Examples of double stochastic operators are provided by composition operators generated by automorphisms on $I$ and averaging operators corresponding to $\sigma$-subalgebras of $\Lambda$. \\

We will denote the set of all double stochastic operators by $\mathfrak{D}$. Let $x \in L^1(I)$. The set $\{Tx\}_{T \in \mathfrak{D}}$ is called (see~\cite{Ry}) the \textit{Ryff's orbit} (or simply the orbit) of $x$ and will be denoted by $\Omega_x$.

A vector space  $X$ such that $L^1(I) \supseteq X \supseteq L^\infty(I)$ $\Big($ respectively, a Banach space $(X,||\cdot ||_X)$, continuously embedded between $L^1(I)$ and $L^\infty(I)$ i.e. $L^1(I) \supseteq X \supseteq L^\infty(I)$ and $\|x\|_{L^1(I)} \leq \|x\|_X \leq \|x\|_{L^\infty(I)},\ x \in X$ $\Big),$ is called an \textit{interpolaton space} $\Big($ respectively, a \textit{strongly interpolation space}$\Big)$ between $L^1(I)$ and $L^\infty(I)$, if the restriction on $X$  of any admissible operator $T$ maps $X$ into itself $\Big($ respectively, if  $$||T||_{X\rightarrow X}\leq \max[||T||_{L^1(I)\rightarrow L^1(I)}, ||T||_{L^\infty(I)\rightarrow L^\infty(I)}]\Big).$$
 The Calderon's interpolation theorem~\cite{Ca}, states that to verify that $X$ is an interpolation space $\Big($ respectively that $(X,||\cdot ||_X)$ is a strongly interpolation space $\Big)$ it is enough to prove that
$$ x \in X, \; y \prec x\ \Rightarrow y \in X, $$
respectively, that
$$ x \in X, \; y \prec x \Rightarrow \|y\|_X \leq \|x\|_X. $$

 Even before the appearance of Calderon's paper~\cite{Ca}  Ryff proved~\cite{Ry} that if $y \preceq x$ than $y=Tx$ where $T \in \mathfrak{D}$. Thus, the \textbf{Ryff - Calderon Theorem} states that to verify that a space is an interpolation (respectively, strongly interpolation) space it is enough to check that $TX \subseteq X$ only for every $T$ in $\mathfrak{D}$,  instead of verifying it for every admissible operator.\\

A vector space $X$, $L^\infty(I) \subseteq X \subseteq L^1(I)$ $\Big($respectively, a Banach space $X$ such that $L^1(I) \supseteq X \supseteq L^\infty(I)$ and $\|x\|_{L^1(I)} \leq \|x\|_X \leq \|x\|_{L^\infty(I)},\ x \in X\Big)$ is called an \textit{order ideal} (respectively, \textit{Banach ideal}) in $L^1(I)$ if
$$ x \in X, |y| \leq |x| \Rightarrow y \in X $$
and, in the case of a Banach ideal, additionally $\|y\| \leq \|x\|$, i.e. the norm $\|\cdot\|$ is \textit{order monotonic}. In the sequel by an ideal we will always mean an order ideal.

It is immediate to see that any interpolation space (a strongly interpolation space) is an ideal (respectively, a Banach ideal). Indeed, if $|y| \leq |x|$ then $y = Mx$ where $M$ is the operator of multiplication by some function $m \in L^\infty(I)$.

We will call interpolation ideals (respectively, Banach interpolation ideals) \textit{majorant} ideals (respectively, \textit{strongly majorant} Banach spaces). By Calderon - Ryff's theorem every majorant ideal together with a function f contains its orbit $\Omega_f$ and, in Banach case, $y \in \Omega_f \Rightarrow \|y\| \leq \|f\|$.

An ideal that is not majorant will be called \textit{non-majorant}.

An ideal $X$, $X \subseteq L^1(I)$, is called \textit{symmetric} if
$$ x \in X, y^\star \leq x^\star \Rightarrow y \in X. $$
It follows from properties of double stochastic operators that every majorant ideal is symmetric.
A symmetric ideal $X$ endowed with a monotonic norm $\|\cdot \|_X$ is called a \textit{symmetric space} if it is norm complete and $\|x\|_X = \|x^\star\|_X, \; x \in X$.

A strongly interpolation Banach ideal $X$, $L^\infty(I) \subseteq X \subseteq L^1(I)$, is called a \textit{strongly majorant} symmetric space.\\

From Lemmas 0.1 and 0.2 follow the inequalities
 $$(\rho_2f^*)(t)\leq f^*(t)\leq(\rho_{\frac{1}{2}}f^*)(t)\leq \Big(f^*(t)+f^*(1-t)\Big)^*,\ t\in I,\ f\in L^1(I),\ \eqno(0.7)$$

It follows from (0.7) that symmetric ideals are invariant under the action of c/d operators.

Let $\{X_\gamma\}_{\gamma\in\Gamma}$ be a family of ideals in $L^1(I)$. The sum of this family is defined by the formula $$\sum_{\gamma\in\Gamma}X_\gamma:=\{\sum_{i=1}^nx_{\gamma_i}:x_{\gamma_i}\in X_{\gamma_i},\ i=1,...,n;\ (\gamma_1,...,\gamma_n)\in \Gamma^n,\ n\geq 1\}.\eqno(0.8)$$
\\

\textbf{Lemma 0.7.} The sum $\sum_{\gamma\in\Gamma}X_\gamma$ of ideals is again an ideal in $L^1(I)$. If $X_\gamma$ are symmetric (respectively, majorant) ideals than their sum is also a symmetric (respectively, majorant) ideal.\\

$\Box$\\

\textbf{Lemma 0.8.} For any ideal $X,\ X\subseteq L^1(I),$ and any $\sigma$-subalgebra $\mathcal{A}$ the image $E\Big(X|\mathcal{A}\Big)$ is an ideal in the space $L^1(I,\mathcal{A},\lambda)$.\\
$\Box$\\

\textbf{Lemma 0.9.} For any set $\mathfrak{X}=\{X_\gamma\}_{\gamma\in \Gamma}$ of ideals $X_\gamma\subseteq L^1(I)$ and any $\sigma$-subalgebra $\mathcal{A}$ the following formula is valid
$$E\Big(\sum_{\gamma\in \Gamma}X_\gamma|\mathcal{A}\Big)=\sum_{\gamma\in \Gamma}E\Big(X_\gamma|\mathcal{A}\Big),\eqno(0.9)$$
where in the right part the sum of ideals is taken in the space $L^1(I,\mathcal{A},\lambda)$.\\
$\Box$\\

\textbf{DEFINITION.} 1. If for a vector ideal $X$ and for operator $E(\cdot|\mathcal{A})$ of averaging by a $\sigma$-subalgebra $\mathcal{A}$ we have $E(X|\mathcal{A})\subseteq X$ then we will say that $\mathcal{A}$ \emph{averages} $X$ or that $X$ is \emph{averaged} by $\mathcal{A}$.\\

2. We will say that a symmetric space $(X,\|\cdot \|_X)$ is \textit{strongly averaged} by a $\sigma$-subalgebra $\mathcal{A}$ $\Big{(}$ or that  $\mathcal{A}$ \textit{strongly averages} symmetric space $(X,||\cdot ||_X)$ $\Big{)}$, if $\mathcal{A}$ averages the symmetric ideal $X$ and also $||E(\cdot|\mathcal{A})||_{X\rightarrow X}\leq 1$.\\

\textbf{Remark 0.10.} Let $(\Omega,\Sigma,\mu)$ be a continuous probability space, i.e. a measure space with a complete non-atomic measure $\mu$, $\mu(\Omega)=1.$ The definitions of a $\sigma$-subalgebra in $\Sigma$ and a complemented $\sigma$-subalgebra in $\Sigma$ remain exactly the same as in the case of $(I, \Lambda, \lambda)$. As above we can introduce the notion of equimeasurable functions $f\in L^1(I)$ and $\tilde{f}\in L^1(\Omega,\Sigma,\mu)$. The sets $A\subseteq L^1(I)$ and $\tilde{A}\subseteq L^1(\Omega,\Sigma,\mu)$ are called $\stackrel{e}{\sim}$\emph{equivalent} if for any $f\in A$ there is $\tilde{f}\in \tilde{A},$ equimeasurable with $f$, and vice versa, for any $\tilde{f}\in \tilde{A}$ there is an equimeasurable $f\in A$ (we will use the notation: $A\stackrel{e}{\sim}\tilde{A}$). Let $\tilde{f}\in L^1(\Omega,\Sigma,\mu)$. Then the non-increasing rearrangement of $\tilde{f}$ is  the function $f^*$, where $f\in L^1(I)$ and $f$ is equimeasurable with $\tilde{f}$. For the space $L^1(\Omega,\Sigma,\mu)$ the definitions of the orbit of a function, symmetric and majorant ideals, and also ideals averaged by a $\sigma$-subalgebra $\tilde{\mathcal{B}}\subseteq\Sigma$ are verbatim repetitions of the corresponding definitions in the case of $(I, \Lambda, \lambda)$. There is a one-to-one correspondence between symmetric ideals in $L^1(I)$ and $\stackrel{e}{\sim}$equivalent to them symmetric ideals in $L^1(\Omega,\Sigma,\mu),$ (see~\cite{Me3} and~\cite{Me1}). Two $\sigma$-subalgebras  $\mathcal{A}\subseteq\Lambda$ and $\tilde{\mathcal{A}}\subseteq\Sigma$ are called $\stackrel{e}{\sim}$\emph{equivalent} if the sets of all indicator functions of all elements of the $\sigma$-subalgebras  $\mathcal{A}$ and $\tilde{\mathcal{A}},$ are $\stackrel{e}{\sim}$equivalent.\\
$\Box$\\

\textbf{Remark 0.11.}\\

1). $\stackrel{e}{\sim}$equivalent symmetric ideals  $X$ and $\tilde{X}$ are either both majorant or both non-majorant.  A $\sigma$-subalgebra  $\mathcal{A}\subseteq\Lambda$ and an $\stackrel{e}{\sim}$equivalent to it  $\sigma$-subalgebra $\tilde{\mathcal{A}}\subseteq\Sigma$ either both average or both do not average  $\stackrel{e}{\sim}$equivalent symmetric ideals $X$ and $\tilde{X}$.\\

2). In virtue of Remark 0.5 if a $\sigma$-subalgebra $\mathcal{A}$ is equimeasurable with another $\sigma$-subalgebra $\mathcal{B}$ and $X$ is a symmetric ideal, then $\mathcal{A}$ averages $X$ if and only if $\mathcal{B}$ does. In particular, countable partitions with the same stochastic vector either all average $X\subseteq L^1(I)$ or all do not. Because symmetric ideals are invariant under the action of contraction/dilation operators all interval partitions equivalent to an interval partition $\mathcal{B}$ average a symmetric ideal $X$ if and only if $\mathcal{B}$ averages it. \\

3) If an interval partition $\mathcal{S}$ is a multiple of the interval partition $\mathcal{T}$, then $\mathcal{S}$ averages a symmetric ideal $X$ if and only if $\mathcal{T}$ averages it.\\

4). The binary projection $\mathcal{B}_{(2)}$  of an interval partition $\mathcal{B}$ averages a symmetric ideal $X\subseteq L^1(I)$ if and only if $\mathcal{B}$ averages it (see~\cite{Me4}).\\

5). It follows from Proposition 0.6 and part 2 of the current remark that if a countable partition $\mathcal{H}$ does not average some ideal then any countable partition $\mathcal{F}$ that is finer then $\mathcal{H}$, as well as any countable partition that is equimeasurable with $\mathcal{F}$, does not average it.  On the other hand, if a countable partition $\mathcal{H}$ averages an ideal $X$, then every sample $\mathcal{H'}$ from $\mathcal{H}$ also averages $X$.\\
$\Box$\\

Because the whole space $L^1(I)$ is a majorant (and therefore, a symmetric) ideal for every set $Z\subseteq L^1(I)$ there exists the smallest (by inclusion) symmetric ideal (majorant ideal) that contains $Z$. We denote this ideal by $\mathcal{N}_Z$ (respectively, $\mathcal{M}_Z$) and say that $Z$ \emph{generates} the symmetric ideal $\mathcal{N}_Z$, (respectively, the majorant ideal $\mathcal{M}_Z$). If $Z$ is a singleton, $Z=\{f\},\ f\in L^1(I),$ we will write $\mathcal{N}_f$ instead of $\mathcal{N}_{\{f\}}$ (respectively, $\mathcal{M}_f$ instead of $\mathcal{M}_{\{f\}}$) and we call this ideal \emph{the principal symmetric} (respectively, \emph{the principal majorant}) ideal \emph{generated by} $f.$ If $\mathcal{A}$ is a non-atomic $\sigma$-subalgebra of the $\sigma$-algebra $\Lambda$, then to avoid any ambiguity we will use the symbol $(\mathcal{A})\mathcal{N}_X\ (\textrm{respectively,}\ (\mathcal{A})\mathcal{M}_X)$ to denote the symmetric (respectively, majorant) ideal in $L^1(I,\mathcal{A},\lambda)$, generated by the set $X\subseteq L^1(I,\mathcal{A},\lambda)$.\\
$\Box$\\

\textbf{Lemma 0.12.} Let $(\Omega,\Sigma,\mu)$ be a non-atomic probability space. For any set $\mathfrak{X}=\{X_\gamma\}_{\gamma\in \Gamma}$ of ideals  $X_\gamma\subseteq L^1(\Omega,\Sigma,\mu)$ and for any $f\in L^1(\Omega,\Sigma,\mu)$ we have
$$\begin{cases}\mathcal{N}_{\sum_{\gamma\in \Gamma}X_\gamma}=\sum_{\gamma\in \Gamma}\mathcal{N}_{X_\gamma};\\
 \mathcal{N}_f=\sum_{\bar{f}\sim f}\mathcal{N}_{\bar{f}}.\end{cases}\eqno(0.10_{\mathcal{N}})$$
$$\begin{cases}\mathcal{M}_{\sum_{\gamma\in \Gamma}X_\gamma}=\sum_{\gamma\in \Gamma}\mathcal{M}_{X_\gamma};\\
 \mathcal{M}_f=\sum_{\bar{f}\sim f}\mathcal{M}_{\bar{f}}.\end{cases}\eqno(0.10_{\mathcal{M}})$$
$\Box$\\

The Lemmas 0.1 and 0.2 allow us to describe the principal symmetric ideal $\mathcal{N}_f$ (respectively, the principal majorant ideal $\mathcal{M}_f$) as follows:\\

\textbf{Lemma 0.13.} For any $f\in L^1(I)$ we have
$$\begin{cases}\mathcal{N}_f=
\mathcal{N}_f=\{z\in L^1(I):\ \textrm{there is a}\ q=q(z)>1,\ \textrm{such that}\ z^*\leq q\rho_{q^{-1}}f^*\}.\\

\mathcal{M}_f=\{z\in L^1(I):\ \textrm{there is a}\ q=q(z)>1,\ \textrm{such that}\ z^{**}\leq q\rho_{q^{-1}}f^{**}\} = \\

=\{z\in L^1(I):\ \textrm{there is a}\ q=q(z)>1,\ \textrm{such that}\ z^{**}\leq qf^{**}\}.\\
\end{cases}\eqno(0.11)$$

The first equality in (0.11) allows us to define $\mathcal{N}_f$ even when $f\in L^0(I)$ but $f \not \in L^1(I)$.\\

\textbf{Lemma 0.14.} For any countable partition $\mathcal{F}$ and any $f\in L^1(I)$ we have
$$\mathcal{N}_{E(\mathcal{N}_f|\mathcal{F})}=\sum_{\hat{f}\sim f,\ \mathcal{\hat{F}}\sim \mathcal{F}}\mathcal{N}_{E(\hat{f}|\mathcal{\hat{F}})},
\ \mathcal{M}_{E(\mathcal{N}_f|\mathcal{F})}=\sum_{\hat{f}\sim f,\ \mathcal{\hat{F}}\sim \mathcal{F}}\mathcal{M}_{E(\hat{f}|\mathcal{\hat{F}})}. \eqno(0.12)$$
$\Box$\\

\textbf{Remark 0.15.} It follows from Lemmas 0.12 and 0.13 that every symmetric ideal $X$ (every majorant ideal $Y$), that is $\stackrel{e}{\sim}$\emph{equivalent} to a principal symmetric ideal $\mathcal{N}_f$ (respectively, to a principal majorant ideal $\mathcal{M}_g$) is itself a principal symmetric ideal $\mathcal{N}_{\bar{f}}$ (respectively, a principal majorant ideal $\mathcal{M}_{\bar{g}}),$ where $\bar{f}\sim f,\ \bar{g}\sim g.$.\\
$\Box$\\

The last equality in $(0.11)$ shows that the principal majorant ideal $\mathcal{M}_f$ generated by $f$ coincides as a subset in $L^1(I)$ with the \emph{Marcinkiewicz space} $\textsl{M}_f$, see~\cite{Lo}. Recall that if $f\in L^1(I)$ then $\textsl{M}_f$ is defined as
$${\textsl{M}_f:=\{x\in L^1(I):\ \|x\|_{\textsl{M}_f}:=\sup_{t\in I}\frac{\int_0^tx^*d\lambda}{\psi(t)}=\sup_{t\in I}\frac{x^{**}(t)}{f^{**}(t)}<\infty\},}\eqno(0.13)$$
where $\psi(t):=\int_0^tfd\lambda$ is an \emph{\textsl{M}-function}, i.e. concave and continuous\\
function on $[0,1]$ such that $\psi(0)=0,\psi(1)=1.$

It follows from the definition $(0.13)$ that every Marcinkiewicz space is a strongly majorant symmetric space.\\

\textbf{Lemma 0.16.} (\cite{Lo}). For any $x,y\in\textsl{M}_f$ we have
$$\|x^*-y^*\|_{\textsl{M}_f}\leq \|x-y\|_{\textsl{M}_f}.$$
$\Box$\\

The closure of $\mathcal{N}_f$ in $\textsl{M}_f$ in the norm $\|\cdot\|_{\textsl{M}_f}$ is denoted by $\textsl{M}_f^1$.\\

\textbf{Lemma 0.17.} (\cite{BM}). The space $\textsl{M}_f^1$ is a symmetric space with the norm $\|\cdot\|_{\textsl{M}_f}$ induced on it from $\textsl{M}_f$.\\
$\Box$\\

\newpage

\bigskip

\centerline{\textbf{ Chapter 1. Principal symmetric and principal majorant ideals}}
\centerline{\textbf{generated by weakly regular and regular functions}}

\bigskip

\centerline{\textbf{The material of this chapter is based on papers ~\cite{Me5}~-\cite{Me9}.}}
\bigskip

\emph{In this chapter we consider the properties of regularity and weak regularity of functions generating principal symmetric and principal majorant ideals.}

\bigskip

Let $s>0$. By $\hat{s}$ we denote $\max[s,s^{-1}]$ (therefore,  $\hat{s}^{-1}=\min[s,s^{-1}]$).\\

\textbf{Lemma 1.1.} For any $f,\ f\in L^1(I),$ the principal symmetric ideal generated by $f$ can be defined by the formula
 $$\mathcal{N}_f=\{x\in L^1:n_f(x)<\infty\},\ \textrm{where}\   n_f(x):=\inf_{s>0}\{x^*(t)\leq \hat{s}\cdot f^*(\hat{s}^{-1}\cdot t),\ t\in I\}.$$

The (nonlinear) functional $n_f$ on $L^1(I)$ introduced in the statement of Lemma 1.1 is called \textbf{the modular}.

Let us list some obvious properties of the modular $n_f$. For any $x,y\in L^1(I)$ and any real $r,\ r\neq 0,$ we have
$$1).\ n_f(x)\geq 0;\ n_f(x)=0\Leftrightarrow x=0;\ 2).\ n_f(r\cdot x)\leq \hat{|r|}n_f(x);$$
$$3).\ n_f(x+y)\leq 2\cdot [n_f(x)+ n_f(y)];$$
$$4).\ n_f(y)\leq n_f(x),\ \textrm{if}\ |y|\leq |x|;\ 5).\ n_f(x)=n_f(y),\ \textrm{if}\ x\sim y.$$
Thus, the modular $n_f$ is monotonic and symmetric on $\mathcal{N}_f$.

It is easy to see that the modular $n_f$ defines a locally convex topology on $L^1(I)$ and that this topology has a countable base of absorbing neighborhoods of zero
$U_k=\{x:\ x^*(t)\leq \max[ 2^k,2^{-k}]\cdot f^*(\min[ 2^k,2^{-k}]\cdot t),\ t\in I\},\ k=0,\pm 1,\pm 2, ...\ .$\\

\textbf{Lemma 1.2} The modular $n_f$ has the Levi - Fatou property, i.e. \\

$6).$ \emph{If$\{x_k\},\ \{x_k\}\subseteq\mathcal{N}_f,$ is a nondecreasing sequence of functions such that $C:=\sup_{k\geq 0}n_f(x_k)<\infty,$ then} $\sup_{k\geq 0}x_k\in \mathcal{N}_f.$\\

$7).$ the modular $n_f$ is complete on $\mathcal{N}_f$.\\

\emph{Proof.} Let $\sup_{k\geq 0}x_k:=x,\ \sup_{k\geq 0}x^*_k:=\bar{x}.$ It is immediate that $\bar{x}$ is a measurable non-increasing function, satisfying the inequality $\bar{x}(t)\leq C\cdot f^*(C^{-1}\cdot t),\ t\in I.$ Thus, $\bar{x}\in\mathcal{N}_f\ \textrm{и}\ n_f(\bar{x})\leq C.$ But $\bar{x}\sim x$, (see~\cite{CR}), and therefore $x\in\mathcal{N}_f\ \textrm{and}\ n_f(x)\leq C.$ For the proof of the implication $(6) \Rightarrow (7)$ see\cite{Am}. \\
$\Box$\\

\textbf{Definition 1.1.} 1. A function $f,\ f\in L^1(I),$ is called \emph{weakly regular} if there is a constant $C=C(f)>0$ such that $$f^*(\frac{t}{2})\leq C\cdot f^*(t),\ t\in I.$$
2. A function $f,\ f\in L^1(I),$ is called \emph{regular} if there is a constant $K=K(f)>0$ such that $f^{**}(t)\leq K\cdot f^*(t),\ t\in I.$\\
$\Box$\\

It is easy to see that weak regularity of $f$ is equivalent to the following: for any $c$, $0<c<1$, there is a $C$, $C>0$, such that $\rho_cf^*(t)\leq C\cdot f^*(t),\ t\in I.$ Therefore, for a weakly regular $f$ we have $$\mathcal{N}_f=\{z\in L^1(I,\Lambda,\lambda):\ z^*\leq q\cdot f^*\ \textrm{for a suitable}\ q=q(z),\ q>0\},\eqno (1.1)$$
and $$n_f(x)=\inf\{s>0:x^*(t)\leq s\cdot f^*(t)\}.\eqno(1.2)$$

Let $f=f^*\in L^1(I),\ \psi(t)=\int_0^tfd\lambda,\ t\in I,\ \textsl{M}_\psi$ be the Marcinkiewicz space with its natural norm $\|\cdot\|_{\textsl{M}_\psi},$ see $(0.12)$. Consider the following function of $s,\ s\in I,$ \\

$\|\rho_s\|_{\textsl{M}_\psi\rightarrow \textsl{M}_\psi}=\sup\{\|\rho_sx\|_{\textsl{M}_\psi\rightarrow \textsl{M}_\psi} :\ \|x\|_{\textsl{M}_\psi,}=1\}$.\\

\textbf{Lemma 1.3}  For $s\in I$ we have $\|\rho_s\|_{\textsl{M}_\psi\rightarrow \textsl{M}_\psi}=\|\rho_sf\|_{\textsl{M}_\psi}$.\\

\emph{Proof}. For any $x\in \textsl{M}_\psi,\  \|x\|_{\textsl{M}_\psi}=1,$ in virtue of Lemma 0.2 we have $\int_0^t(\rho_sx)^*d\lambda=\int_0^t\rho_sx^*d\lambda\leq \int_0^t\rho_sfd\lambda$. Now, the statement of the lemma follows from the fact that  $\|f\|_{\textsl{M}_\psi}=1$ .\\
$\Box$\\

\textbf{Theorem 1.4} For any $f\in L^1(I)$ the following conditions are equivalent.\\

$\imath)$. There is a norm $||\cdot||$ on $\mathcal{N}_f$ such that the topology defined by this norm is equivalent to the topology defined by the modular $n_f$;\\

$\imath\imath).$ The principal symmetric ideal $\mathcal{N}_f$ coincides with the principal majorant ideal $\mathcal{M}_f$;\\

$\imath\imath\imath)$.\  Any $\sigma$-subalgebra of $\Lambda$ averages the principal symmetric ideal $\mathcal{N}_f$;\\

$\imath v)$.  $\textrm{the function}\ f \ \textrm{is regular},$ i.e., $f^{**}\in \mathcal{N}_f.$\\

$v).\ \textrm{There are constants}\ k>0\ \textrm{and}\ p\in (0,1),\ \textrm{such that}\ \|\rho_s \|_{M_\psi\rightarrow M_\psi}\leq k\cdot s^{-p},\ 0<s\leq 1.$\\

$v\imath).\ \liminf_{t\rightarrow 0}\frac{\psi(2t)}{\psi(t)}>1.$\\

\emph{Proof.} $\imath)\Rightarrow\imath\imath).$ Let $||\cdot||\simeq n_f$ on $\mathcal{N}_f$ and let $G$ be a convex bounded neighborhood of zero in the topology induced by the modular $n_f$, such that $G$ is the unit ball for the norm $||\cdot||$. Then $\sup\{n_f(g):\ g\in G\}:=C<\infty.$ Let $H=\{h\in \mathcal{N}_f:\ |h|\leq |g|\ \textrm{for a suitable}\   g\in G \}$ and denote by $F$ the convex hull of $H$. Then $\sup\{n_f(x):\ x\in F\}\leq 2\cdot C, \textrm{because} \sup\{n_f(x):\ x\in H\}\leq C.$ Thus $F$ is a convex and bounded neighborhood of zero in the topology generated by the modular $n_f$. A routine reasoning shows that $F$ is an order ideal set. Therefore the Minkowski's functional of $F$,
\cite{Rud}, generates on $\mathcal{N}_f$ a monotonic norm $||\cdot ||'$ equivalent to $n_f$. $\mathcal{N}_f$ endowed with this norm becomes a Banach space, and because the norm $||\cdot ||'$ has the Levi-Fatou property, $\mathcal{N}_f$ is a majorant ideal (see e.g.~\cite{CR}). Thus, by definition  $\mathcal{N}_f=\mathcal{M}_f$ and the implication is proved.\\

The inverse implication $\imath\imath)\Rightarrow\imath)$ follows from the completeness and monotonicity of both $||\cdot||_{\textsl{M}_\psi}$ and $n_f$ on $\mathcal{M}_f$.\\

The equivalence $\imath\imath)\Leftrightarrow\imath v)$ in one dirction is trivial and in the opposite direction will be proved below when we prove the implication $\imath\imath\imath)\Rightarrow\imath v).$ \\

$\imath\imath)\Rightarrow \imath\imath\imath).$ If $\mathcal{N}_f=\mathcal{M}_f$ then $\mathcal{N}_f$, being a majorant ideal, is averaged by any $\sigma$-subalgebra of $\Lambda$, because the average of a function is contained in its orbit. \\

We will now prove the implication $iii)\Rightarrow iv).$\\

Assume that for some $f=f^*\in L^1(I)$ the inclusion $f^{**}\in\mathcal{N}_f$ is false. We will construct such an interval partition $\mathcal{W}$ that $E(f|\mathcal{W})\notin \mathcal{N}_f$ and thus come to a contradiction.\\

By definition of the principal symmetric ideal $\mathcal{N}_f$ if $f^{**}\notin\mathcal{N}_f$, then we can find a sequence $\{t_n\}\subseteq I,\ t_{n+1}<t_n,\ n\geq 1,\ t_1=1$ such that for any $m\geq 1$ we have
$$\sup_n\frac{f^{**}(t_n)}{f(\frac{t_n}{m})}=\infty.$$
It follows that for a suitable sequence of natural numbers $\{n_m\},\ n_1=1,$ we have
$$1<c_m:=\frac{f^{**}(s_m)}{f(\frac{s_m}{m})}\uparrow\infty,$$
where $s_m=t_{n_m},\ m\geq 1,\ s_1=t_{n_1}=1.$\\

Let us fix a $\varepsilon,\ 0<\varepsilon<c_1.$ On the interval $(0,s_1)$ we consider the continuous nondecreasing function
$$c_1-\varepsilon<J_1(u):=\frac{(s_1-u)^{-1}\int_u^{s_1}fd\lambda}{f(\frac{s_1}{1})}.$$
There is a point $s_{m_2},\ s_{m_2}<\frac{s_1}{2}$ such that
$$\frac{(s_1-s_{m_2})^{-1}\int_{s_{m_2}}^{s_1}fd\lambda}{f(\frac{s_1}{1})}>c_1-\varepsilon.$$
Assume that the sequence of points $\{s_{m_j}\},\ s_{m_j}<\frac{s_{m_{j-1}}}{2},\ j=2...,k,$ is constructed in such a way that $$\frac{(s_{m_j}-s_{m_{j-1}})^{-1}\int_{s_{m_{j-1}}}^{s_{m_j}}fd\lambda}{f(\frac{s_{m_j}}{m_j})}>c_{m_j}-\varepsilon,$$
we will define on the interval $(0,s_{m_k})$ a continuous nondecreasing function
$$J_{k+1}(u):=\frac{(s_{m_k}-u)^{-1}\int_u^{s_{m_k}}fd\lambda}{f(\frac{s_{m_k}}{m_k})}$$
and chose a point $s_{m_{k+1}}<\frac{s_{m_k}}{2}$ such that
$$\frac{(s_{m_k}-s_{m_{k+1}})^{-1}\int_{s_{m_{k+1}}}^{s_{m_k}}fd\lambda}{f(\frac{s_{m_k}}{m_k})}>c_{m_k}-\varepsilon.$$
Because for the inductively constructed sequence $\{s_{m_k}\}_k$ we have $s_{m_k}\downarrow 0,$ by setting $w_k:=s_{m_k},\ k\geq 1$ we obtain that the interval partition $\mathcal{W}=(w_k)$ does not average $\mathcal{N}_f$.\\

The equivalence $\imath v)\Leftrightarrow v)$ was proved in~\cite{Lo1} and ~\cite{Sh}.\\

$\imath v)\Rightarrow v\imath)$. Assume that there is a connstant $C\geq 1$ such that $f^{**}\leq Cf$. Then for the function $\psi(t):=\int_0^tfd\lambda,\ t\in I,$ we have
$$\frac{\psi(t)}{\psi(t/2)}=\frac{\int_0^tfd\lambda}{\int_0^{t/2}fd\lambda}=1+\frac{\int_{t/2}^t fd\lambda}{\int_0^{t/2}fd\lambda}\geq 1+\frac{f(t)}{f^{**}(t/2)}\geq 1+\frac{f(t)}{2f^{**}(t)}\geq 1+\frac{f(t)}{2Cf(t)}=1+\frac{1}{2C},$$
and the implication is proved. \\

Let us prove the inverse implication: $v\imath)\Rightarrow\imath v).$ Because $f^*(t)\geq \ 1/t \int_t^{2t}f^*d\lambda,\ t\in (0,2^{-1}],$ the implication follows from the relations $$\liminf_{t\rightarrow 0}\frac{_{\int_t^{2t}f^*d\lambda}}{\psi(t)}>0\Leftrightarrow \limsup_{t\rightarrow 0}\frac{\psi(t)}{\int_t^{2t}f^*d\lambda}<\infty\Leftrightarrow\limsup_{t\rightarrow 0}\frac{f^{**}(t)}{1/t\int_t^{2t}f^*d\lambda}<\infty\Rightarrow$$
$$\Rightarrow\limsup_{t\rightarrow 0}\frac{f^{**}(t)}{f^*(t)}<\infty,$$
which is exactly the condiiton $\imath v)$.\\

The equivalence $v\imath)\Leftrightarrow\imath)$ folows from  $\imath v)\Leftrightarrow v\imath),\ \imath\imath)\Leftrightarrow\imath v)$ and $\imath\imath)\Leftrightarrow\imath)$.\\
$\Box $\\

\textbf{Corollary 1.5} A function $f=f^*\in L^1(I)$ is regular if and only if the function $f^{**}$ is regular.\\

\emph{Proof.} The implication ($f$ is regular) $\Rightarrow$ ($f^{**}$ is regular) is obvious. Let us prove the inverse implication. Denote $\psi^{**}(t):=\int_0^tf^{**}d\lambda,\ t\in I$, and notice that in virtue of Theorem 1.4 $\imath v)\ f^{**}\in \textsl{M}_\psi$, therefore there is a $c>0$ such that $c^{-1}\|x\|_{\textsl{M}_\psi}\leq\|x \|_{\textsl{M}_{\psi^{**}}}\leq c\|x\|_{\textsl{M}_\psi}$ for any $x\in \textsl{M}_\psi$. Next, it follows from the monotonicity of the norm $\|\cdot\|_{\textsl{M }_\psi}$, from the equality $f=f^*\leq f^{**}$, and from Lemma 1.3 and Theorem 1.4 that we can chose $k>0$ and $p\in (0,1)$ such that
$$\|\rho_s \|_{M_\psi\rightarrow M_\psi}=\|\rho_sf\|_{M_\psi}\leq \|\rho_sf^{**}\|_{M_\psi}\leq c\|\rho_sf^{**} \|_{\textsl{M}_{\psi^{**}}}=c\|\rho_s\|_{\textsl{M}_{\psi^{**}}\rightarrow\textsl{M}_{\psi^{**}}}\leq c\cdot ks^{-p},\ s\in I.$$
It remains to apply again Theorem 1.4 (v) to see that $f$ is regular.\\
$\Box$\\

In~\cite{Me6} and~\cite{Me8} the author considered the connections between the properties of regularity and weak regularity of functions from $L^1(I)$.\\

\bigskip

\textbf{Theorem 1.6.} \\

1). For any $f\in L^1(I)$ and any $\varepsilon>0$ there is a weakly regular function $g\in L^1(I)$, such that $\mathcal{M}_f=\mathcal{M}_g$ and also
$$g^*(t/2)\leq (2+\varepsilon)\cdot g^*(t),\ t\in I.$$
2). If $f=f^*\in L^1(I)$ then there is a $g=g^*$, equivalent to $f$ and such that
$$g^*(t/2)\leq 2\cdot g^*(t),\ t\in I,$$
if and only if $f\simeq h^{**}$ for some $h=h^*\in L^1(I)$.\\

3). If $f\in L^1(I)$ then there are $g\in L^1(I)$ and $\varepsilon,\ 0<\varepsilon<2$ such that $\mathcal{M}_f=\mathcal{M}_g$ and also
$$g^*(t/2)\leq (2-\varepsilon)\cdot g^*(t),\ t\in I,$$
if and only if $f^*\simeq f^{**}$, i.e. if and only if $f$ is a regular function.\\

\emph{Proof.} Let us prove 1). Without loss of generality we can assume that $f^*$ is $\mathcal{D}$-measurable: $f^*(t)=\sum_{n\geq 1}a_n\cdot \textbf{1}_{D_n}.$\\

Let us represent $f^*$ as $f^*=\sum_{n\geq 0}b_n\cdot f_n$, where $f_n=\textbf{1}_{(0,2^{-n}]},\ n\geq 0.$ Let us fix $\varepsilon>0$ and let $c=2+\varepsilon.$ Define the set $\Phi(c):=\{g=g^*\in L^1(I):\ g^*(t/2)\leq c\cdot g^*(t),\ t\in I\}$ and notice its obvious properties:\\

a). $g_1,g_2\in\Phi(c),\ \alpha_1,\alpha_2>0\Rightarrow \alpha_1\cdot g_1+\alpha_2\cdot g_2\in\Phi(c);$\\

b). If the sequence $\{g_n\}\subset\Phi(c)\ \textrm{and }\ g_n\uparrow x\in L^1(I),\ \textrm{then}\ x\in\Phi(c).$\\

Let $g_0=1$ and\\
$$g_n(t):=\begin{cases}1, &\mbox{if}\ 0<t\leq 2^{-n};\\
c^{k-n},& \mbox{if}\ 2^{-k-1}<t\leq 2^{-k}: k=0,1,...,n-1, \end{cases}t\in (0,1],$$

$n=1,2,...$ and define $g:=\sum_{n\geq 0}b_n\cdot g_n.$ It is obvious from our construction that $g_n\in\Phi(c), $ and therefore, according to  a) and b) $g\in\Phi(c).$ Also, in virtue of obvious inequalities $f_n\leq g_n,\ n\geq 1$ we have $f^*\leq g$, and therefore $\mathcal{M}_f\subseteq \mathcal{M}_g.$ If we prove that
$$\int_0^tgd\lambda\leq\frac{c-1}{c-2}\int_0^tfd\lambda,\ t\in I,\eqno(1.3)$$
then the inverse inclusion and thus the statement 1 of the theorem would be proved.\\

For any $ t\in I$ we have
$$\int_0^tgd\lambda=\sum_{n\geq 0}b_n\int_0^tg_nd\lambda;\ \int_0^tfd\lambda=\sum_{n\geq 0}b_n\int_0^tf_nd\lambda,$$
and therefore (1.3) would follow  from the inequalities
$$\int_0^tg_nd\lambda\leq \frac{c-1}{c-2}\int_0^tf_nd\lambda,\ n\geq 0,$$
which we will prove next. \\

By definition
$$\int_0^tf_nd\lambda=\begin{cases} t, &\mbox{if}\ 0<t\leq 2^{-n};\\
2^{-n},& \mbox{if}\ 2^{-n}<t\leq 1. \end{cases}. $$

$$\int_0^tg_nd\lambda=\begin{cases} t &\mbox{if}\ 0<t\leq 2^{-n};\\
2^{-n}+c^{-1}(t-2^{-n})& \mbox{if}\ 2^{-n}<t\leq 2^{-n+1};\\
2^{-n}+c^{-1}2^{-n}+c^{-2}(t-2^{-n+1})& \mbox{if}\ 2^{-n+1}<t\leq 2^{-n+2};\\
........................................\\
2^{-n}+c^{-1}2^{-n}+c^{-2}(t-2^{-n+1})+...+c^{-n}(t-2^{-1})& \mbox{if}\ 2^{-1}<t\leq 1.
\end{cases}. $$

Thus $\int_0^tg_nd\lambda=\int_0^tf_nd\lambda$ if $0<t\leq 2^{-n},$ and
$$\int_0^tg_nd\lambda=2^{-n}[1+c^{-1}(1+2/c+2^2/c^2+...)]=2^{-n}\frac{c-1}{c-2}=\frac{c-1}{c-2}\int_0^tf_nd\lambda,$$
where $2^{-n}<t\leq 1.$  Statement 1 is thus proved.\\

2). By Lemma 0.1 for any $h\in L^1(I)$
we have the inequality $h^{**}(t/2)\leq 2 h^{**}(t),\ t\in I$ and from it statement  2) follows in one direction. Conversely, let $f^*\simeq g^*$ where $g^*\in L^1(I)$ is such that $g^*(t/2)\leq 2\cdot g^*(t),\ t\in I.$  Without loss of generality we can assume that $g^*$ is a $\mathcal{D}$-measurable function: $g^*(t)=\sum_{n\geq 1}a_n\cdot \textbf{1}_{D_n}:\ D_n=(2^{-n},2^{-n+1}],\ 0<a_0\leq a_1\leq...\leq a_n\uparrow \infty,\ \sum_{n\geq 1}2^{-n}a_n<\infty.$ The inequalities $a_{m+n}\leq 2^ma_n,$ which are true for $g^*$ for any $m,n\geq 1$ guarantee that $t_0g^*(t_0)\leq 4t_1g^*(t_1)$ for all $0<t_0\leq t_1\leq 1.$\\

Let $\varphi(0)=0,\ \varphi(t)=tg^*(t)$ if $0<t\leq 1$ and  $\varphi(t)=tg^*(1)$ if $t>1.$ The above inequalities for $g^*$ guarantee that the function $\varphi$ is quasi-concave on $[0,\infty)$ and on this interval it is equivalent to its smallest concave majorant
$$\tilde{\varphi}:\ 2^{-1}\tilde{\varphi}\leq\varphi\leq\tilde{\varphi},\eqno(1.4)$$
(see~\cite{KPS}. Because $\lim_{t\rightarrow 0}\tilde{\varphi}(t)=\lim_{n\rightarrow \infty}2^{-n}a_n=0,$ we also have $\lim_{t\rightarrow 0}\varphi(t)=0.$ Then the concave function $\tilde{\varphi}$ on $[0,1]$ can be represented as $\tilde{\varphi}(t)=\int_0^thd\lambda,$ where $h=h^*\in L^1(I,\Lambda,\lambda).$ Now, it follows from $(1.4)$ that $f^*(t)\simeq g^*(t)\simeq\frac{\varphi(t)}{t}\simeq \frac{\tilde{\varphi}(t)}{t}=h^{**},\ t\in I.$ Statement 2 is proved.\\

3). Let $\psi(t)=tf^{**}(t),\ t\in I$. It follows from Theorem 1.4 that the relation $f^*\simeq f^{**}$ is equivalent to the condition $\liminf_{t\rightarrow 0}\frac{\psi(t)}{\psi(t/2)}>1,$ or, equivalently, $\liminf_{t\rightarrow 0}\frac{f^{**}(t) t}{f^{**}(t/2) t/2 }>1,\ t\in I.$ By considering, if necessary,  instead of $f$ an equivalent to it function, and by choosing an appropriate $0<\delta<1$ we can write the last inequality as $\frac{f^{**}(t)t}{f^{**}(t/2)t/2 }>1+\delta,\ t\in I,$ or $f^{**}(t/2)\leq \frac{2}{1+\delta}f^{**}(t),\ t\in I.$ By renaming the constants we get that $f^{**}(t/2)\leq 2(1-\varepsilon)f^{**}(t),\ t\in I,\ 0<\varepsilon<1.$ The theorem is completely proved.\\
$\Box$\\

\textbf{Corollary 1.7} Let $f$ be an arbitrary function from $L^1(I)$. For any $\gamma,\ 0<\gamma<1,$ there is $g\in L^1(I)$, such that $\mathcal{M}_g=\mathcal{M}_f$ and the function $h:=g^\gamma$ is regular.\\

\emph{Proof.} By Theorem 1.6.1) we can find an $\varepsilon$, such that for some $\delta,\ 0<\delta<1,$ $(2+\varepsilon)^\gamma< 2-\delta.$ Assume that $g=g^*\in L^1(I)$ is a weakly regular function such that $\mathcal{M}_g=\mathcal{M}_f$ и $g(t/2)\leq (2+\varepsilon)g(t),\ t\in I.$ Let $h=g^\gamma,$ then
$$h(t/2)=[g(t/2)]^\gamma\leq (2+\varepsilon)^\gamma[g(t)]^\gamma< (2-\delta)h(t),\ t\in I.$$
It remains to apply criterion 3 of Theorem 1.6.\\

\textbf{Theorem 1.8} The following conditions are equivalent

(1) $f\in L^1(I)$ is regular.

(2) $\mathcal{M}_g=\mathcal{M}_f \Rightarrow g$ is regular.

(3)$\mathcal{M}_g=\mathcal{M}_f \Rightarrow g$ is weakly regular.\\

\emph{Proof.} (1) $ \Rightarrow (2) $. Let $f$ be regular. By Theorem 1.4 the equality $\mathcal{M}_g=\mathcal{M}_f$ implies that the inequalities $\liminf_{t\rightarrow 0}\frac{\psi_g(2t)}{\psi_g(t)}>1$ and $\liminf_{t\rightarrow 0}\frac{\psi_f(2t)}{\psi_f(t)}>1$ are equivalent. Thus, $g$ is regular.\\

The implication $(2) \Rightarrow (1)$ is trivial.

To prove the implication  $(3) \Rightarrow (1)$ we will consider an additional construction.

Let us fix $h=h^*\in L^1,\ \int_0^1hd\lambda=1$ and let $\psi_h(t)=\int_0^thd\lambda,\ t\in I$. For any $n\geq 1$ let\\

$q_h(n):=$ \emph{the number of points of the sequence }$\{\psi_h(2^{-k})\}_{k\geq 1}$, \emph{in the binary interval}  $D_n=(2^{-n},2^{-n+1}].$\\

Because the function $\psi_h$ is concave we have $q_h(n)\geq 1,\ n\geq 1.$ It is immediate that the equalities $q_g(n)=q_h(n),\ n\geq 1$  for the functions $g=g^*\in L^1(I)$ and $h=h^*\in L^1(I)$ imply that $\mathcal{M}_g=\mathcal{M}_h.$\\

\textbf{We will now interrupt the proof of Theorem 1.8 and prove two lemmas we need to finish the proof of implication} $(3) \Rightarrow(1)$.

\bigskip

\textbf{Lemma 1.9} For any sequence $\{q_n\}_{n\geq 1}$ of positive integers there is a $g=g^*\in L^1(I,\Lambda,\lambda)$, such that $q_n=q_g(n),\ n\geq 1.$\\

\emph{Proof. } Assume that in each interval $D_n$ we have constructed points $\{u_i^n:\ 2^{-n}<u^n_{q_n}<...<u^n_1\leq 2^{-n+1};\ n\geq 1\}$ in such a way that if we number them in decreasing order (from $1$ to $0$)  we obtain a strongly decreasing to zero sequence $\{v_k\}_{k\geq 0}$ and also
$$v_0=1,\ v_k-v_{k+1}\geq \frac{v_{k-1}-v_k}{2},\ k\geq 1.\eqno(1.5) $$
It follows from $(1.5) $ that $\mathcal{D}$-measurable function $g$ defined as
$$g(t):=\sum_{k\geq 1}2^k(v_{k-1}-v_k)1_{D_k}(t),\ t\in I,\eqno(1.6)$$
is a non-negative non-increasing on $I$. Additionally we have the equalities
$$\int_0^{2^{-k}}gd\lambda =\sum_{i\geq k+1}2^i2^{-i}(v_{i-1}-v_i)=v_k,\ k\geq 1.\eqno(1.7)$$
If we define $\psi_g(0)=0$, $\psi_g(t)=\int_0^tgd\lambda,\ t\in I$  and take as the sequence $\{q_n\}$ the sequence $\{q_{\psi_g}(n)\}$ then we obtain that $\mathcal{M}_g=\mathcal{M}_f.$ \\

To construct the points $\{u_i^n\}$ by the given sequence $\{q_n\}=\{q_f(n)\}$ let us introduce the following terminology. We will call a binary half-segment $D_n$ \emph{a single-point-like}, if $q_n=1$ and \emph{multi-points-like} otherwise. A maximal by inclusion element of the set of all sequences of the form $\{(D_n,D_{n+1},...,D_{n+l-1}):\ q_n=q_{n+1}=q_{n+l-1}=1,\ l\geq 1\}$ will be called \emph{a single-point-like block of length} $l$. A half-segment $D_{n-1}$ will be called \emph{initial} for this block. Consider now the construction of points $\{u_i^n\}$ in the following four cases.\\

$\mathcal{U}_1$. If there is a single-point-like block $(D_1,D_2,...,D_l)$ of length $l,\ 1\leq l\leq\infty$ let $u_i^n=2^{-i+1},\ 1\leq i\leq l.$\\

$\mathcal{U}_2$.  If a half segment $D_n$ is multi-point-like and not an initial one, then $u_i^n=2^{-n}(1+2^{-i+1}),\ 1\leq i\leq q_n.$\\

$\mathcal{U}_3$. If a multi-point-like half segment $D_n$ is initial for a single-point-like block$\{D_{n+1},...,D_{n+l-1},...\}$ of infinite length, then $u_i^n=2^{-n}(1+2^{-i+1}),\ 1\leq i\leq q_n$ and $u_i^{n+i}=u_{q_n}^n\cdot 2^{-i},\ i\geq 1.$\\

$\mathcal{U}_4$. Let now the half segment$D_n$ be initial for a single-point-like block of  length $l$. Define now a single point on every half segment of this block as $u_1^{n+i}=2^{-n-i}(1+2^{-l+i-1}),\ i=1,...,l$. We will construct now points on $D_n$ depending on the relation between $l$ and $q_n$.

$$u_i^n=\begin{cases} \textrm{if}\ q_n< l+3,\ \textrm{then}\ 2^{-n}(1+2^{-i+1})\ \textrm{for}\ i=1,2,...,q_n-1\ \textrm{and}\       2^{-n}(1+2^{-l-1})\ \textrm{for}\ i=q_n;$$\\
\\     \textrm{if }\ q_n\geq l+3,\ \textrm{then}\ 2^{-n}(1+2^{-i+1}),\ \textrm{for}\ i=1,2,...,q_n.\end{cases}. $$

Clearly the set $\mathcal{U}=\bigcup \{u_i^n\}$ is completely defined by the rules $\mathcal{U}_1 -         \mathcal{U}_4.$ The sequence $\{v_n\}$ obtained as indicated above from the points of $\mathcal{U}$ is strongly decreasing to 0 (and $v_0=1$). It is immediate to check that  inequlities $(1.5)$ are satisfied.\\
$\Box$\\

\textbf{Lemma 1.10} The following condiitons are equivalent.
(1) A $\mathcal{D}$-measurable function $f=f^*$ is regular.

(2) The sequence $\{q_{\psi_f}(n)\}$ is bounded.\\

\emph{Proof.} $(2) \Rightarrow (1)$. Assume that there is an integer $d$ such that $\{q_{\psi_f}(n)\}\leq d,\ n\geq 1.$ For an arbitrary $t\in (0,2^{-1}]$ we chose $k\geq 1$ such that $2^{-k-1}<t\leq 2^{-k}$. Let  $v_{k-1}\in D_n$. By definition $v_{k-1}-v_k\geq 2^{-n-q_{\psi_f}(n)}$. On the other hand it follows from $(1.6)$ that $f(2t)=f(2^{-k+1})=2^k(v_{k-1}-v_k)$ and therefore by $(1.7)$
$$\psi(t)\leq \psi(2^{-k})=v_k<v_{k-1}\leq 2^{-n+1}\leq 2^{q_{\psi_f}(n)+1}(v_{k-1}-v_k)=2^{q_{\psi_f}(n)+1}2^{-k}f(2t)\leq 2^{d+2}tf(2t);$$
$$\frac{\psi(2t)}{\psi(t)}=1+\frac{\int_t^{2t}fd\lambda}{\psi(t)}\geq 1+\frac{tf(2t)}{\psi(t)}\geq 2^{-d-2}.$$
Thus $\liminf\frac{\psi(2t)}{\psi(t)}>1$ and therefore $f$ is regular.\\

$(1) \Rightarrow (2)$. Assume that $\lim_{k\rightarrow\infty} q_{\psi_f}(n_k)=\infty.$ Our construction guarantees that
$$u^{n_k}_{q_{\psi_f}(n_k)-1}-u^{n_k}_{q_{\psi_f}(n_k)}\leq 2^{-n_k}2^{-q_{\psi_f}(n_k)+2};\ u^{n_k}_{q_{\psi_f}(n_ k)}>2^{-n_k};\ u^{n_k}_{q_{\psi_f}(n_k)}-u^{n_k+1}_1\geq 2^{-n_k+3}.\eqno(1.8)$$
Let $u^{n_k}_{q_{\psi_f}(n_k)}=v_{m_k};\ u^{n_k}_{q_{\psi_f}(n_k)-1}=v_{m_k-1}.$ It follows from $(1.8)$ that

 $$\frac{\psi(2\cdot 2^{-m_k})}{\psi(2^{-m_k})}=\frac{v_{m_k-1}}{v_{m_k}}=1+\frac{v_{m_k-1}-v_{m_k}}{v_{m_k}}\leq 1+2^{-q_{\psi_f}(n_k)+2}.$$
Thus $\liminf\frac{\psi(2t)}{\psi(t)}=1$ and by Theorem 1.4 $f$ cannot be regular.\\
$\Box$\\

\textbf{Remark} 1.11. In view of $(0.4)$ the condition in the statement of Lemma 1.10 that $f$  is $\mathcal{D}$-measurable can be omitted.\\

We can now finish the proof of Theorem 1.8. Let us remind the reader that assuming $\lim\frac{\psi(2t)}{\psi(t)}=1$ we need to construct a function $h=h^*\in L^1(I,\Lambda,\lambda)$ such that $h$ is not weakly regular and $\mathcal{M}_h=\mathcal{M}_f,$ where $\psi(t)$  is the concave function from the proof of Lemma 1.10 and $f(t)=\frac{d\psi}{dt}.$ We will show that we can just put $h=f$.\\

By using the notations from the proof of Lemma 1.10 and applying the inequalities $(1.8)$ we obtain
$$\frac{f(2^{-m_k})2^{-m_k-1}}{f(2^{-m_k}+1)2^{-m_k}}=\frac{v_{m_k}-v_{m_k+1}}{v_{m_k-1}-v_{m_k}}\geq \frac{2^{-n_k-3}}{2^{-n_k}2^{-q_{\psi_f}(n_k)+2}}.$$
Thus, $\frac{f(2^{-m_k})}{f(2^{-m_k+1})}\geq 2^{q_{\psi_f}(n_k)-5}$ and by Lemma 1.10 $\limsup_{k\rightarrow\infty} q_{\psi_f}(n_k)=\infty$. Therefore $f$ cannot be weakly regular.
Theorem 1.8 is now completely proved.\\
$\Box$\\

\textbf{Remark 1.12.} 1). The regularity of a function $f$ is a topolical invariant of the Marcinkiewicz space $\textsl{M}_f$. I.e.,if for the functions $f$ and $g$ we have $\mathcal{M}_f=\mathcal{M}_g$ and $f$ is regular, then $g$ is also regular.\\
2). Because the equality $\mathcal{N}_f=\mathcal{N}_g$ implies that $\mathcal{M}_f=\mathcal{M}_g$ it follows from Theorem 1.6.3) that if symmetric ideals generated by $f$ and $g$ coincide and one of these functions is regular then the second is also regular.\\
$\Box$\\

Our next result concerning regularity is stated in terms of averaging operators generated by interval partitions.\\

\textbf{Theorem 1.13} (compare with Theorem 5.12 below). Let $f=f^*\in L^1,\ \mathcal{B}=(b_n)$ - be some interval partition. The following conditions are equivalent.\\

1). $f$ is regular and $\mathcal{N}_f=\mathcal{N}_{E(f|\mathcal{B})}$;\\

2). The function $E(f|\mathcal{B})$ is regular.\\

\emph{Proof.}  The implication $1)\Rightarrow 2)$ follows immediately from Remark 1.12.2).\\

Let us prove that $2)\Rightarrow 1).$ Let $g:=E(f|\mathcal{B})\subseteq \mathcal{M}_f$. Assume first that $f^{**}\in\mathcal{N}_g$. Then $f^{**}\in\mathcal{M}_g\subseteq \mathcal{M}_f,$ and the regularity of $f$ follows from Theorem 1.4. \\

Assume now that $f^{**}\notin\mathcal{N}_g$. Let us represent $g$ in the form $g=\sum_{n\geq 1}\alpha_n1_{(b_n,b_{n-1}]},$ where $\alpha_n=\frac{\int_{b_n}^{b_{n-1}}fd\lambda}{b_{n-1}-b_n},\ n\geq 1.$ It follows from our assumption that there is a sequence $\{n_k\}_{k\geq 0},\ n_0=0$ such that $b_{n_k}\downarrow 0$ and
$$\frac{f^{**}(b_{n_k})}{g(b_{n_k})}\geq 2^k,\ \textrm{or}\ \frac{f^{**}(s_k)}{g(s_k)}\geq 2^k,\ k\geq 0.\eqno (_1.9)$$
where, to simplify our notations, we put $b_{n_k}:=s_k,\ k\geq 0.$ Without loss of generality we can assume that $s_1<2^{-1}.$ Applying inequalities $(_1.9)$ we can inductively construct the sequence $\{k_m\}$ of natural numbers such that for the interval partition $\mathcal{U},\ \mathcal{U}=(u_m),$ where $u_0=1,\ u_m=s_{k_m},\ m\geq 1,$ we have $E(f|\mathcal{U})\notin\mathcal{N}_g.$ Because by definition $\mathcal{U}$ is coarser than $\mathcal{B}$ we get $E(f|\mathcal{U})=E(g|\mathcal{U})\in\mathcal{M}_g$,
and,because $g$ is regular $\mathcal{M}_g=\mathcal{N}_g$ - a contradiction.\\

Let $k_1=1$ and assume that for $m\geq 1$ the natural number $k_m$ has been already constructed. On the interval $[0,u_m)$ consider the function $\Phi,\ \Phi(u):=\frac{(u_m-u)^{-1}\int_u^{u_m}fd\lambda}{g(u_m)}.$ Obviously $\Phi$ is a non-increasing function continuous on $[0,u_m)$ and it follows from $(_1.9)$ that $\Phi(0)\geq 2^{k_m}.$ Let us chose the number $k_{m+1}$ so large that $k_{m+1}>k_m,\ \Phi(s_{k_{m+1}})>2^{k_m}-1.$\\

Now, when the sequence $\{k_m\}_{m\geq 1}$ has been constructed let us put $v_m=\min[2u_m,\tau_m]$ where $\tau_m$ is the closest to $u_m$ from the right point of the partition $\mathcal{B},\ m\geq 1.$\\
Let $h=h^*=E(g|\mathcal{U})$. For the monotonic functions $g$ и $h$ from the construction follow the inequalities
$$\frac{h(v_m/2)}{g(v_m)}\geq  2^{k_m}-1,\ m\geq 1.$$
Therefore, $h(\frac{t}{2})\notin\mathcal{N}_g,$ and consequently $h\notin\mathcal{N}_g,$ in contradiction with Theorem 1.4.\\

Thus $f^{**}\in\mathcal{N}_{E(f|\mathcal{B})},$ and from it follows, as was shown above, the regularity of $f$. Finally, because  $f^{**}\simeq f,$ we obtain $\mathcal{N}_f=\mathcal{N}_{E(f|\mathcal{B})}$.\\
$\Box$\\

\textbf{Remark 1.14.} 1. It is worth noticing that for any unbounded function $f=f^*\in L^1(I)$ there is an interval partition $\mathcal{B}$ such that the function $E(f|\mathcal{B})$ is not only not regular but even not weakly regular.\\
2. It follows from Remark 1.12.2 that the statement of Theorem 1.13 remains true not only for principal symmetric ideals but for principal majorant ideals. Namely, we have the following theorem.\\

\textbf{Theorem 1.13.}$^\prime$ Let $f=f^*\in L^1$ and let $\mathcal{B}=(b_n)$ be some interval partition. the following conditions are equivalent.\\

1). $f$ is regular and $\mathcal{M}_f=\mathcal{M}_{E(f|\mathcal{B})}$;\\

2). The function $E(f|\mathcal{B})$ is regular.\\

A function  $f,\ f\in L^1(I),$ is called a \emph{ function of bounded mean oscillation on} $I$ $\Big(\textrm{and we write}: f\in BMO(I)\Big)$ if

$$||f||_{BMO(I)}:=\sup_{[a,b]\subseteq I}\frac{1}{b-a}\int_a^b \Big{|}f(s)-\frac{1}{b-a}\int_a^bf(u)du \Big{|}ds< \infty.$$

\bigskip

\textbf{Theorem 1.15}. Let $g=g^*\in L^1(I)$. The following conditions are equivalent.\\

$(\imath)\ g\in BMO(I);$\\

$(\imath\imath)\ g^{**}(t) - g(t)\leq C,\ t\in I,$ where $0<C$ is a constant;\\

$(\imath\imath\imath)\ g(\frac{t}{2}) - g(t)\leq K,\ t\in I,$ where $0<K$ is a constant;\\

$(\imath v)\ g(t)=\log_2 f(t),\ t\in I,$ where $f=f^*\in L^1(I)$ is a weakly regular function.\\

\emph{Proof.} The equivalence of $(\imath)$ and $(\imath\imath)$ was proved in \cite[Theorem 7.10, p.382]{BS}. To prove the implication $(\imath\imath)\Rightarrow(\imath\imath\imath)$ we will need the following lemma that is an analogue of the well known in mathematical statistics \textit{property of median} (See e.g.~\cite[pp.178-179]{Cra}). We omit the proof of the lemma, that can be easily obtained by comparison of areas under the corresponding graphs.\\

\textbf{Лемма 1.16}. Let $g=g^*\in L^1(I).$ For any $t\in I$ we have
$$\min_a\int_0^t|g(s)-a|ds=\int_0^t |g(s)-g(\frac{t}{2})|ds.$$\\

\bigskip

We return to the proof of Theorem 1.15. Assume  $(\imath\imath).$ Then, in virtue of $\imath)$ we have $\|g(t)\|_{BMO}<\infty$ and therefore from the definition of $BMO(I)$ and Lemma 1.16 we obtain that for any $t\in I$
$$|g^{**}(t)-g(\frac{t}{2})|=|\frac{1}{t}\int_0^t[g(s)-g(\frac{t}{2})]ds|\leq\frac{1}{t}\int_0^t|g(s)-g(\frac{t}{2})|ds\leq$$
$$\frac{1}{t}\int_0^t|g(s)-\frac{1}{t}\int_0^tg(u)du|ds\leq \|g(t)\|_{BMO}.$$

From it and the triangle inequality we obtain the implication $(\imath\imath)\Rightarrow(\imath\imath\imath)$. Next, if we consider $f=\exp_2g$ then it follows from already proved implications that $(\imath\imath\imath)\Rightarrow(\imath v).$ It remains to prove that $(\imath v) \Rightarrow (\imath\imath).$ In view of $(0.4)$ we can assume that $f:=\sum_{n\geq 1}\alpha_n\cdot 1_{D_n}$, where $D_n=[2^{-n},2^{-n+1})$ и $1\leq \alpha_1\leq ...\leq \alpha_n\leq \alpha_{n+1}\leq Q\cdot \alpha_n,\ n\geq 1,\ Q>1.$\\

Let $\beta_n=\log_2\ \alpha_n,\ n\geq 1;\ g=\sum_{n\geq 1}\beta_n\cdot 1_{D_n}.$ For any $n\geq 1$ we have $\beta_{n+1}\leq \tilde{Q}+ \beta_n$, where $\tilde{Q}=\log_2\ Q>0$. Therefore,
$$g^{**}(2^{-n})=\sum_{k\geq n+1}\beta_k\cdot 2^{-k}\leq\sum_{k\geq 1}
\beta_n\cdot 2^{-k}+\tilde{Q}\cdot\sum_{k\geq 1}k\cdot 2^{-k}=\beta_n+C=g(2^{-n})+C,\ n\geq 1$$
где $C=\tilde{Q}\cdot\sum_{k\geq 1}k\cdot 2^{-k}<\infty.$

Thus, for any $t\in D_n$ we have $g^{**}(t)\leq g^{**}(2^{-n})\leq g(2^{-n})+C=g(t)+C$.\\
$\Box$\\

\textbf{Corollary 1.17}. For any $g=g^*\in BMO(I)$ we can find positive constants $C_1$ and $C_2$ such that
$$g(t)\leq C_1- C_2\cdot \ \ln\ t,\ t\in I.$$
\emph{Proof.} Let $g(t)=\log_2\ f(t)$, where $f=f^*\in L^1,\ f(\frac{t}{2})\leq Q\cdot f(t),\ t\in I.$ Without loss of generality we can assume that $f(1)=1,$ or $g(1)=0.$ Let $\gamma,\ 0<\gamma<1,$ be such that $Q^\gamma<2$ and let $h=f^\gamma.$ Because $h(\frac{t}{2})\leq Q^\gamma\cdot h(t),\ t\in I,$ from Theorem 1.6 we obtain that $h$ is a regular function. Therefore we can find a constant $K>1$ such that
 $$\frac{h(t)}{\int_0^thd\lambda}\geq \frac{1}{K\cdot t},\ t\in I.$$
By integrating this inequality by $t$ from $s$ to $1$, where $s$ is an arbitrary number from $I$ we obtain
$$\ln\ s^{-\frac{1}{K}}\leq \ln\ \frac{\int_0^1hdt}{\int_0^shdt}\ ,$$
and therefore $s\cdot h^{**}(s)=\int_0^shd\lambda\leq (\int_0^1hd\lambda)\cdot s^{\frac{1}{K}}.$ Thus, $f^\gamma(s)=h(s)\leq h^{**}(s)\leq (\int_0^1hd\lambda)\cdot s^{(\frac{1}{K}-1)}$. It remains to take logarithms of both parts of the last inequality and change the names of constants, if nrcessary.\\
$\Box$\\

\textbf{Theorem 1.18}. For any regular function $f=f^*\in L^1(I)$ we can find the indicator $\delta=\delta(f),\ 0<\delta<1$, and a constant $C>0$, such that $f^{**}(t)\leq Ct^{-\delta},\ t\in I.$\\

\emph{Proof.} Let $f^{**}\leq Kf^*,\ K\geq 1.$ Without loss of generality we can assume that $f$ is not a constant function, i.e. $K>1$. Then,
$$f^*(t)[\int_0^tfd\lambda]^{-1}\geq [Kt]^{-1},\ t\in I.$$
By integrating this inequality over the interval $[r,1]$, where $r$ is an arbitrary number from $I$, we obtain
$$\int_r^1\frac{f^*(t)}{\int_0^tf^*d\lambda}d\lambda \geq K^{-1}\int_r^1\frac{dt}{t}=\ln [(\frac{1}{r})^{\frac{1}{K}}].$$
By representing the left part in the form
$$\int_r^1\frac{d\int_0^tf^* d\lambda}{\int_0^tf^* d\lambda}=\ln\frac{\int_0^1f^* d\lambda}{\int_0^rf^* d\lambda},$$
we see that
$$\frac{\int_0^1f^* d\lambda}{\int_0^rf^* d\lambda}\geq (\frac{1}{r})^{\frac{1}{K}}.$$
Therefore, $f^{**}(r)\leq Cr^{-\delta}$, где $C=\int_0^1f^* d\lambda,\ \delta=1-\frac{1}{K}.$\\
$\Box$\\

\newpage

\bigskip \centerline{\textbf{Chapter 2. Two absolute constants for rearrangements of intervals.}}
\bigskip
\centerline{\textbf{The results of this chapter are based on paper~\cite{AM}.}}
\bigskip

\emph{In this Chapter we compute two absolute constants for rearranging a countable family of subintervals of $I$ in such an order that their lengths are non-increasing when moving from 1 to 0.}

\bigskip

Let $\vec{a}$ be an arbitrary stochastic vector and let $\mathcal{B}=(b_n)$ be an interval partition with stochastic vector $\vec{a}$. Let $\mathcal{B}^*=(b^*_n)$ be the monotonic rearrangement of $\mathcal{B}$ with stochastic vector $\vec{a}^*$.\\

\textbf{Theorem 2.1.}. The golden ratio $\alpha:=\frac{5^{1/2}+1}{2}$ is the smallest among all constants $\delta$ such that.
$$\begin{cases}\textrm{For\ any\ integer}\ n\geq 1\ \textrm{there}\ \textrm{is} \ \textrm{an} \ \textrm{integer}\ m\geq 0,\ \textrm{such}\ \textrm{that}\\
\delta^{-1}\cdot b^*_m<b_n<\delta\cdot b^*_m.\\
\end{cases}\eqno (2.1)$$
\bigskip
\emph{Proof.} First notice that
$$\alpha-\alpha^{-1}=1. \eqno (2.2)$$
Assume now that the golden ratio $\alpha$ does not satisfy  $(2.1)$ and let $\vec{a}=(a_n)$ be a stochastic vector for which (2.1) becomes false if $\delta=\alpha$. Then, because (2.1) is equivalent to $|\log_{\alpha} b_n - \log_{\alpha} b^*_m|<1$, we can find an integer $N>0$ such that for any integer $m\geq 0$ we have
 $$|\log_{\alpha}b_N - \log_{\alpha} b^*_m|\geq 1. \eqno (2.3)$$
Let us fix an integer $m\geq 0$ such that
$$b^*_{m+1}<b_N\leq b^*_m. \eqno (2.4)$$
In virtue of ($2.3$), $\alpha b_{m+1}^*\leq b_N\leq \alpha^{-1}b_m^*$ and combining it with ($2.2$) we obtain $a^*_{m+1}=b^*_m-b^*_{m+1}\geq \alpha b_N-\alpha^{-1}b_N=b_N.$ From it, from the monotonicity of the stochastic vector $\vec{a}^\star$, and taking into consideration the positions of points of the interval partition $\mathcal{B}$ in $I$ we get
$$a_1^*\geq a_2^*\geq...\geq a^*_{m+1}\geq  b_N> a_{N+j},\ j\geq 1. \eqno (2.5)$$
Let $\gamma$ be a bijection of the set of natural numbers such that
$$a_n^*=a_{\gamma(n)},\ n\geq 1. \eqno (2.6)$$
It follows from ($2.5$) and ($2.6$) that for any integers $i=1,...,m+1$ and any integer $j\geq 1$  $a_{\gamma(i)}=a_i^*\neq a_{N+j}.$ Or, equivalently,
$$\gamma(i)\neq N+j,\ i=1,...,m+1;\ j\geq 1. \eqno (2.7)$$
Finally, from ($2.7$) we can conclude  $\gamma(i)\leq N,\ i=1,...,m+1.$ Now we get a contradiction with ($2.4$): $1-b_{m+1}^*=1-b^*_1+b^*_1-b^*_2...+b_m^*-b^*_{m+1}=a_1^*+...+a^*_{m+1}=a_{\gamma(1)}+...+a_{\gamma(m+1)}\leq a_1+...+a_N=1-b_1+b_1-b_2+...+b_{N-1}-b_N=1-b_N.$\\

We have proved that $\alpha$ satisfies condition ($2.1$). It remains to prove that for any $\delta:\ 1\leq\delta<\alpha,$ we can find a stochastic vector $\vec{a}$ not satisfying (2.1).\\

Fix $\varepsilon,\ 0<\varepsilon<\min [\frac{\alpha-1}{\delta}+\alpha-2,10^{-1}],$ and consider the stochastic vector with coordinates $a_1=2-\alpha, a_2=\alpha-1-\varepsilon,a_n=\varepsilon\cdot 2^{-n+2},\ n\geq 3.$ Simple computations show that $a_2>a_1>a_3>a_4>...$; therefore $a_1^*=\alpha-1-\varepsilon,a_2^*=2-\alpha,a_n^*=a_n,\ n\geq 3.$ Applying ($2.2$) we see that: $\delta^{-1}b_0^*=\delta^{-1}>\alpha^{-1}=\alpha-1=b_1=\delta\cdot (\frac{\alpha-1}{\delta}+\alpha-2+2-\alpha)>\delta\cdot(2-\alpha+\varepsilon)=\delta\cdot b_1^*.$\\

Thus $\delta^{-1}b_0^*>b_1>\delta b_1^*>\delta b_2^*>...>\delta b_n^*>...,$ and for any integer $m\geq 0$ it cannot be true that $\delta^{-1}b_m^*<b_1<\delta b_m^*$.\\
$\Box$\\

Let $\vec{a},\ \vec{a}^*,\ \mathcal{B}=(b_n),\ \mathcal{B}^*=(b^*_n)$ denote the same objects as above. We associate with a function $f=f^*\in L^1(I)$  three M-functions: $\psi(t)=\int_0^tfd\lambda;\ \psi_{\mathcal{B}}(t)=\int_0^tE(f|\mathcal{B})d\lambda,\ \psi_{\mathcal{B}^*}(t)=\int_0^tE(f|\mathcal{B}^*)d\lambda,\ t\in I.$\\

Notice that if $\delta\geq \alpha$ then from Theorem 2.1 easily follows the inequality
$$\psi_{\mathcal{B}}(t)\leq \delta\cdot\psi_{\mathcal{B}^*}(t),\ t\in I.\eqno (2.8)$$
Therefore it is quite natural to assume that the smallest possible constant in ($2.8$) is $\alpha$.\\
But, interestingly, this assumption is not correct and instead the following statement takes place. \\

\textbf{Theorem 2.2}. The smallest possible constant in ($2.8$) is $\beta:=4/3<\alpha$.\\

\emph{Proof.} We will prove first that $\beta$ satisfies ($2.8$). Consider the inequality
$$\psi_{\mathcal{B}}(b_n)\leq 4/3\cdot\psi_{\mathcal{B}^*}(b_n),\ n\geq 0. \eqno (2.9)$$
We will prove now that if $(2.9)$ is correct than we have ($2.8$) with $\beta$ instead of $\delta$.\\

 Let $\varphi(t):=4/3\cdot\psi_{\mathcal{B}^*}(b_n),\ b_n<t\leq b_{n-1},\ n\geq 1.$ It follows from $\varphi(b_n)= 4/3\cdot\psi_{\mathcal{B}^*}(b_n),\ n\geq 0$ that the graph of the concave function $4/3\cdot\psi_{\mathcal{B}^*}(t)$, that connects the points $(b_n,\varphi(b_n))$, must be not lower than the graph of the piecewise linear nondecreasing function $\varphi(t)$, that connects the same points. I.e.  $\varphi(t)\leq 4/3\cdot\psi_{\mathcal{B}^*}(t),\ t\in I.$ Therefore if ($2.9$) is not satisfied then for any $t\in (b_{n+1},b_n],\ n\geq 0$ we have
$$\psi_{\mathcal{B}}(t)=\frac{t-b_{n+1}}{b_n-b_{n+1}}\psi(b_n)+\frac{b_n-t}{b_n-b_{n+1}}\psi(b_{n+1})\leq 4/3\frac{t-b_{n+1}}{b_n-b_{n+1}}\psi_{\mathcal{B}^*}(b_n)+$$
$$+4/3\frac{b_n-t}{b_n-b_{n+1}}\psi_{\mathcal{B}^*}(b_{n+1})=\frac{t-b_{n+1}}{b_n-b_{n+1}}\varphi(b_n)+\frac{b_n-t}{b_n-b_{n+1}}\varphi(b_{n+1})=\varphi(t)\leq 4/3\cdot\psi_{\mathcal{B}^*}(t).$$
Thus, we need to prove ($2.9$). Let us fix an integer $N\geq 0$. If for some $m\geq 0$ we have  $b_N=b^*_m$ then
$$\psi_{\mathcal{B}}(b_N)=\psi(b_N)=\psi(b^*_m)=\psi_{\mathcal{B}^*}(b^*_{m})\leq 4/3\psi_{\mathcal{B}^*}(b^*_{m})=4/3\psi_{\mathcal{B}^*}(b_N).$$
assume now that the inequality $b_N\neq b^*_m$ is false for any integer $m\geq 0$. Fix an $m\geq 0$ such that $b^*_{m+1}<b_N<b^*_m,$ and introduce the notations: $d=b^*_{m+1},\ b=b^*_m,\ w=b_N.$\\

To finish the proof of Theorem 2.2 we will need the following three lemmas.\\

\textbf{Lemma 2.3.} If $u>0$ and $v>0$ are such that $u^{-1}+v^{-1}<4$ then $u+v>1$.\\

\emph{Proof.} Assume that $u+v\leq 1.$ Then $(u+v)(u^{-1}+v^{-1})=1+uv^{-1}+vu^{-1}+1<4\Leftrightarrow uv^{-1}+vu^{-1}<2$. If $a:=uv^{-1}$ then  $vu^{-1}=a^{-1}$, and we would obtain that $a+a^{-1}<2\Leftrightarrow a^2-2a+1<0\Leftrightarrow (a-1)^2<0,$ a contradiction.\\
$\Box$\\

\textbf{Lemma 2.4.} The following inequality holds.
$$w[(w-d)^{-1}+(b-w)^{-1}]\geq 4. \eqno (2.10)$$
\emph{Proof.} Assume to the contrary that: $[\frac{w-d}{w}]^{-1}+[\frac{b-w}{w}]^{-1}<4$. Then it follows from Lemma 2.3 that
$\frac{b-d}{w}=\frac{b-w}{w}+\frac{w-d}{w}>1.$ Therefore,
$$a_1^*\geq a_2^*\geq ...\geq a_{m+1}^*=b-d>w=b_N>a_{N+j},\ j\geq 1,$$
i.e. we have inequalities ($2.5$). Now, as in the proof of Theorem 2.1, we obtain a contradiction with the inequality $b_N>b^*_{m+1}.$ \\
$\Box$\\

\textbf{Lemma 2.5.} From inequality $(2.10)$ follows the inequality
$$w[\frac{w-d}{b-d}w+\frac{b-w}{b-d}d]^{-1}\leq\frac{4}{3} \eqno (2.11)$$
\emph{Proof.} Let us write $(2.10)$ in an equivalent form:
$$w[(w-d)^{-1}+(b-w)^{-1}]\geq 4\Leftrightarrow \frac{w}{w-d}+\frac{w}{b-w}\geq 4\Leftrightarrow w(w-d+b-w)\geq 4(w-d)(b-w)\Leftrightarrow$$
$$\Leftrightarrow w(b-d)\geq 4(w-d)(b-w)\Leftrightarrow 3(w-d)(b-w)\leq \frac{3}{4}w(b-d).$$
Let us also use the following equivalent form of $(2.11)$:
$$w[\frac{w-d}{b-d}w+\frac{b-w}{b-d}d]^{-1}\leq\frac{4}{3}\Leftrightarrow \frac{3}{4}w(b-d)\leq (w-d)w+(b-w)d.$$
Thus $(2.10)$ implies $(2.11)$ if
$$3(w-d)(b-w)\leq (w-d)w+(b-w)d\Leftrightarrow \frac{w}{b-w}+\frac{d}{w-d}\geq 3\Leftrightarrow \frac{w}{b-w}+\frac{d}{w-d}+\frac{w-d}{w-d}\geq 4\Leftrightarrow$$
$$\Leftrightarrow \frac{w}{b-w}+\frac{w}{w-d}\geq 4.$$
But the last statement exactly means that $(2.10)$ is satisfied.\\
$\Box$\\
Let us continue the proof of Theorem 2.2.

Combining the conditions $\psi(b)\geq \psi(w),\ \frac{\psi(d)}{d}\geq\frac{\psi(w)}{w}$ for the function $\psi(t),\ t\in (0,1]$ with inequality $(2.11)$ we get
$$\frac{\psi_{\mathcal{B}}(w)}{\psi_{\mathcal{B}^*}(w)}=\frac{\psi(w)}{\frac{w-d}{b-d}\psi(b)-\frac{b-w}{b-d}\psi(d)}\leq \frac{\psi(w)}{\frac{w-d}{b-d}\psi(w)-\frac{b-w}{b-d}\frac{d}{w}\psi(w)}=$$$$=w[\frac{w-d}{b-d}w+\frac{b-w}{b-d}d]^{-1}
\leq\frac{4}{3}.$$
Thus we have proved that ($_2.8$) is true with $\beta$ instead of $\delta$. Let us prove now that for any $\gamma,\ 1\leq \gamma<\beta$ we can find an interval partition $\mathcal{B}=(b_n)$, an $M$-function $\psi$, and a $t\in (0,1]$ such that ($_2.8$)is false if $\delta=\gamma$я.\\

Let us define the function $\varphi(t),\ \varphi(t)=t\cdot 1_{[0,2/3]}(t)+2/3\cdot 1_{[2/3,1]}(t),\ t\in (0,1].$ Fix $\varepsilon,\ 0<\varepsilon<\min[1/3,(8-6\gamma)(6-3\gamma)]$ and let $b_0=1,\ b_1=2/3,\ b_n=\varepsilon\cdot 2^{-n},\ n\geq 2.$ Then $\varphi_{\mathcal{B}}=\varphi,$ and a simple computation shows that $b_0^*=1,\ b_1^*=\varepsilon/2+1/3,\ b_n^*=b_n,\ n\geq 2;\ \frac{\varphi_{\mathcal{B}}(2/3)}{\varphi_{\mathcal{B}^*}(2/3)}
=\frac{(8-6\varepsilon)}{(6-3\varepsilon)}>\gamma.$\\

By putting $\psi=3/2\cdot\varphi$ we obtain the required interval partition $\mathcal{B}=(b_n)$, the $M$-function $\psi$, and the point $t=2/3.$\\
$\Box$\\

\bigskip

\newpage

\bigskip

\bigskip

\centerline{ \textbf{Chapter 3. Double stochastic projections and interpolation property}}

\centerline{\textbf{ of ideals between $L^1$ and $L^\infty$}.}
\bigskip
\centerline{\textbf{The results of this chapter are based on papers ~\cite{Me11},~\cite{Me12},\ and~\cite{Me13}.}}
\bigskip

\emph{We prove that if every countable partition averages a vector ideal $I$ then $I$ is a majorant ideal.}\\

\bigskip

\textbf{Proposition 3.1}. If every at most countable partition averages vector ideal $X,\ X\subseteq L^1(I),$ then $X$ is a symmetric ideal.\\

\emph{Proof.} 1). We need to prove the implication
$$x\in X,\ y\ \textrm{is equimeasurable with}\ x\Rightarrow y\in X.\eqno(3.1)$$
Without loss of generality we can assume that $x,y\geq \textbf{0}.$ Moreover, we can assume that both functions $x$ and $y$ are elementary functions. Indeed, for any $x$ and $y$ we can find elementary functions $x_1$ and $y_1$ such that
$$\textbf{0}\leq x_1\leq x\leq x_1+\textbf{1},\ \textbf{0}\leq y_1\leq y\leq y_1+\textbf{1},\eqno(3.2)$$
Thus, if implication $(3.1)$ is true for elementary functions then by $(3.2)$ it will be true for any functions $x$ and $y$.\\

2). Let $x \in X$ be an elementary function of the form $x=\sum_{i=1}^nr_i\cdot 1_{B_i}$ and let $y \in L^1(I)$ be equimeasurable with $x$. Then  $y=\sum_{i=1}^nr_i\cdot 1_{C_i}$, where $C_i,B_i\in \Lambda$ and $\lambda(C_i)=\lambda(B_i);\ i=1,...,n$. Notice that the assumptions of the Proposition guarantee that $\mathbf{1} \in X$.  Therefore $x,y\leq \bar{r}\cdot\textbf{1}\in X$, where $\bar{r}:=\max_{i=1}^nr_i$, and thus $y\in X.$\\

3). Let $x=\sum_{i=1}^\infty r_i\cdot 1_{B_i}$, where $\sigma(B_i,\ i\geq 1)$ is an infinite countable partition of $I,$ and let $y$ be a function from $L^1(I)$ equimeasurable with $x$, i.e. $y=\sum_{i=1}^\infty r_i\cdot 1_{C_i},$ where $0<\lambda(C_i)=\lambda(B_i),\ 1\leq i<\infty.$\\

Assume first that there is an index $n_0$ such that $\lambda(B_i\cap C_j)=0$ if $i,j\geq n_0.$ It follows from part 2 of the proof that it is enough to prove that the function $\bar{y}$ defined as
$$\bar{y}(t):=\begin{cases} \sum_{i\geq n_0}r_i\cdot 1_{C_i}(t), &\mbox{if } t\in\bigcup_{i\geq n_0} C_i;\\
0,& \mbox{otherwise}. \end{cases}. $$
is in $X$, if we know that the function $\bar{x}$ defined as
$$\bar{x}(t):=\begin{cases} \sum_{i\geq n_0}r_i\cdot 1_{B_i}(t), &\mbox{if } t\in\bigcup_{i\geq n_0} B_i;\\
0,& \mbox{otherwise}. \end{cases}. $$
is in $X$.
Notice that $\bar{x}(t)=0$ if $t\in C_i,\ i\geq n_0$.\\

Let us consider the partition $\mathcal{F}:=\sigma \Big(F_0,F_i;\ i\geq n_0\Big),$ where $F_i:=B_i\cup C_i,\ i\geq n_0,\ F_0:=I\setminus\bigcup_{i\geq n_0}F_i$. By assumption $E(\bar{x}|\mathcal{F})\in X$ and for any $t\in F_i,\ i\geq n_0,$ we have
$$E(\bar{x}|\mathcal{F})(t)=\frac{1}{\lambda(F_i)}\int_{F_i}\bar{x}d\lambda= \frac{1}{2\lambda(B_i)}\int_{B_i\cup C_i}\bar{x}d\lambda=\frac{1}{2\lambda(B_i)}\int_{B_i}\bar{x}d\lambda=\frac{r_i}{2}.$$
If, on the other hand, we assume that $t\in F_0$ then $E(\bar{x}|\mathcal{F})(t)=0$. Thus, $\textbf{0}\leq\bar{y}\leq 2E(\bar{x}|\mathcal{F})$ and therefore $\bar{y}\in X.$\\

4). Let us consider the general case. Let again $F_i:=B_i\cup C_i,\ i\geq 1$. Because $\lambda(\bigcup_{i\geq j}F_i)\leq \sum_{i\geq j}[\lambda(B_i)+\lambda(C_i)]=2\sum_{i\geq j}\lambda(B_i)\downarrow_{j\rightarrow\infty}0,$ there is an index $n_0$ such that $\lambda(\bigcup_{i\geq n_0}F_i)<\frac{1}{2}.$\\
By Proposition 0.3 we can find pairwise disjoint subsets $B_i',\ i\geq n_0$, of the set $I\setminus\bigcup_{i\geq n_0}F_i$ such that $\lambda(B'_i)=\lambda(B_i),\ i\geq n_0.$ Then the functions $y'$ and $x'$ are equimeasurable, where
$$y'(t):=\begin{cases} r_i, &\mbox{if }\ t\in B'_i,\ i\geq n_0;\\
0,& \mbox{if}\ t\notin \bigcup_{i\geq n_0}B'_i. \end{cases}; $$
$$x'(t):=\begin{cases} r_i, &\mbox{if }\ t\in B_i,\ i\geq n_0;\\
0,& \mbox{if}\ t\notin \bigcup_{i\geq n_0}B_i. \end{cases}. $$

Because $x\in X$, it follows from part 3 of the proof that $y'\in X$. On the other hand, $y'$ is equimesurable with the function
$$\tilde{y}(t):=\begin{cases} r_i, &\mbox{if}\ t\in C_i,\ i\geq n_0;\\
0,& \mbox{if}\  t\notin \bigcup_{i\geq n_0}C_i. \end{cases}. $$
Applying again part 3 of the proof we obtain that $\tilde{y}\in X.$ From it and from part 2 follows that $y\in X.$ \\
$\Box$\\

\textbf{Theorem 3.2.} Let $X$ be a vector ideal in $L^1(I)$ (respectively, a Banach ideal) such that any at most countable partition averages $X$. Then \\

I. Together with every function ideal $X$ contains its orbit and therefore by the Calderon - Ryff theorem it is a majorant ideal.\\

II. Respectively, if $||E(\cdot|\mathcal{F})||_{X\rightarrow X}\leq 1$ for any at most countable partition $\mathcal{F}$ then $(X,||\cdot ||_X)$ is a strongly majorant ideal.\\

\textit{Proof}. Proof of part I. By Proposition 3.1. $X$ is a symmetric ideal. we need to prove that it is a majorant ideal.\\

First let us prove the implication $y^*\prec x^*\in X\Rightarrow y\in X $ when $y^*\in L^1(I)$ is an elementary function.\\

At this point we need to interrupt the proof of Theorem 3.2 in order to state and proof several lemmas.

\bigskip

\textbf{Lemma 3.3.} Let $[a,b]$ be an arbitrary closed finite interval and let $f$ and $g$ be positive functions on $[a,b]$ satisfying the following conditions.\\

1. $f,g \in L^1(a,b);$\\
2. For some natural $n$
$$\left\{\begin{array}{rcl}
g=\sum_{j=1}^n\alpha_j\cdot 1_{\triangle_j},\ \textrm{where}\ \alpha_1\geq \alpha_2\geq ...\geq\alpha_n\geq 0;\  \triangle_j:=     [t_{j-1},t_j],\ j=1,...,n;\\ a=t_0<t_1<...<t_n=b;\ t_j\leq\frac{t_{j-1}+t_{j+1}}{2},\ j=1,2,...,n.\\
\end{array}
\right. (3.3)$$
3. $$
\left\{\begin{array}{rcl}
\int_{\triangle_j}fd\lambda\geq\int_{\triangle_j}gd\lambda,\ j=1,...,n-1;\\
\int_a^b fd\lambda>\int_a^b gd\lambda;\ \int_{t_1}^b fd\lambda<\int_{t_1}^b gd\lambda.\\
\end{array}
\right.
$$
Then we can find a function $\bar{f}$ on $[a,b]$ that is equimeasurable with  $f$ and such that $$\int_{\triangle_j}\bar{f}d\lambda\geq\int_{\triangle_j}gd\lambda,\ j=1,...,n.$$

\emph{Proof}. The proof will consist of two parts.\\

I. Let us prove that under the assumptions of Lemma 3.3 we can find a function  $f'\sim f$ such that\\

$\int_{\triangle_j}f'd\lambda\geq\int_{\triangle_j}gd\lambda,\ j=1,...,n-1,\ \int_{t_1}^b f'd\lambda\geq\int_{t_1}^b gd\lambda.$\\

First let us notice that from the assumptions of Lemma 3.3 follows the inequality $\int_{\triangle_n}fd\lambda<\int_{\triangle_n}gd\lambda.$
Without loss of generality we can assume that $f$ is non-increasing on $\triangle_n$. Indeed, we can consider monotonic rearrangement of $f$ on $\triangle_n$ without changing it on $[a,t_{n-1}];$ the resulting function is equimeasurable with $f$.\\

It follows from inequalities $(3.3)$ that $t_1-a\leq b-t_{n-1}\leq b-t_1.$ We define the following two functions of $r$:\\

 $$
\left\{\begin{array}{rcl}
u(r):=\int_{b-r}^{b}fd\lambda+\int_{a+r}^{t_1}fd\lambda -\int_{a}^{t_1}gd\lambda;\\
v(r):=\int_a^{a+r} fd\lambda+\int_{t_1}^{b-r}fd\lambda-\int^{b}_{t_1}gd\lambda,\\
\end{array}
\right.
r\in [0,t_1-a].$$
Let us notice they following obvious properties of these functions.\\

a). They are both continuous on $[0,t_1-a];$\\

b). $u(0)=\int_{\triangle_1}fd\lambda-\int_{\triangle_1}gd\lambda\geq 0;$\\
$v(0)=\int_{t_1}^bfd\lambda-\int_{t_1}^bgd\lambda<0;$\\

c). On $[0,t_1-a]$ we have the identity
$$u(r)+v(r)=\int_a^bfd\lambda-\int_a^bgd\lambda\equiv\ const>0.$$
Assume that we can find a point $r_0\in [0,t_1-a]$ such that $u(r_0)\geq 0,\ v(r_0)\geq 0.$ Then the rearrangement of $f$ on the intervals of equal length  $[a,a+r_0]$ and $[b-r_0,b]$ provides the needed function $f'.$\\
Let us prove that we can always find a point $r_0$ with such properties. If not, then it would follow from properties a) - c) of functions $u$ and $v$ that $u(r)>0$ for any $r\in [0,t_1-a].$ In particular, $u(t_1-a)>0,$ i.e.
$$\int_{b-t_1+a}^bfd\lambda-\int_a^{t_1}g\lambda>0.$$
Because $g$  is non-increasing on $[a,b]$ then
$$\int_{b-t_1+a}^bfd\lambda-\int_{b-t_1+a}^bgd\lambda>0.\eqno(3.4)$$
Because $[b-t_1+a,b]\subset\triangle_n,$ from the previous inequality and from the inequality $\int_{\triangle_n}fd\lambda<\int_{\triangle_n}gd\lambda$ follows that
$$\int_{t_{n-1}}^{b-t_1+a}fd\lambda<\int^{b-t_1+a}_{t_{n-1}}gd\lambda.\eqno(3.5)$$
The assumption that $f(r)\leq g(r)$ for any $r\in[b-t_1+a,b]$ would contradict the inequality $(3.4)$.
Therefore there is  a point $r'\in[b-t_1+a,b]$ such that $f(r')> g(r')=\alpha_n.$ But then, because $f$ is non-increasing on $\triangle_n$, the inequality $f(r)> g(r)$ remains true on the whole interval $[t_{n-1}, r']$ that contains the interval $[t_{n-1},b-t_1+a]$, in contradiction with($3.5$). That concludes the first part of the proof.\\

II. We will finish the proof of Lemma 3.3 applying induction by $n$. For $n=2$ the statement of the lemma follows from part I. Assume that the statement of the lemma is true for any finite interval for indices  2,3,...,$n-1$. Let us prove it for index $n$. By part I we can construct on the interval $[a,b]$ a function $\tilde{f}$ equimeasurable with $f$ and such that $$\int_{t_1}^b\tilde{f}d\lambda\geq \int_{t_1}^bg d\lambda,\ \int_{\triangle_j}\tilde{f}d\lambda\geq \int_{\triangle_j}gd\lambda,\ j=1,...,n-1.$$
Let $k$ be the smallest natural number in the interval $[2,n-1]$ such that $$\int_{t_k}^b\tilde{f}d\lambda<\int_{t_k}^bgd\lambda.$$
(If no such $k$ exists then $\tilde{f}$ is the function we need). On the interval $[t_{k-1},b]$ the functions $\tilde{f}$ and $g$ satisfy all the conditions of the current lemma. By inductive assumption there is a function $\hat{f}$ defined on the interval $[t_{k-1},b]$, equimeasurable with the restriction of $\tilde{f}$ on this interval, and such that
 $$\int_{\triangle_j}\hat{f}d\lambda\geq \int_{\triangle_j}gd\lambda,\ j=k,k+1,...,n.$$
But then the function $\bar{f}$ defined on $[a,b]$ as follows
$$\bar{f}(w):=\begin{cases}\tilde{f}(w),\ &\mbox{if }\ w\in [a,t_{k-1}];\\
\hat{f}(w),& \mbox{if}\ w\in [t_{k-1}, b] \end{cases} $$
has the properties we need.\\
$\Box$

\textbf{Lemma 3.4.} Let $f$ and $g$ be defined on $[a,b]$, satisfy conditions 1 and 2 of Lemma 3.3, and such that
$$4.\ \int_a^{t_j}fd\lambda>\int_a^{t_j}gd\lambda,\ j=1,...,n.$$
Then we can find a function $\bar{f}$ defined on $[a,b]$, equimeasurable with $f$, and such that  $$\int_{\triangle_j}\bar{f}d\lambda\geq \int_{\triangle_j}gd\lambda,\ j=1,...,n.$$
\emph{Proof} We prove the lemma using induction by $n$. Let $n=2$. If $\int_{\triangle_2}fd\lambda\geq \int_{\triangle_2}gd\lambda$ then $\bar{f}=f$. Otherwise, the existence of $\bar{f}$ follows directly from Lemma 3.3.\\

Assume that the statement of the lemma is true for any finite closed interval for $2,3,..., n-1.$ Let us prove that it remains true for $n$.\\

First we will prove that without loss of generality we can assume that the inequality $$\int_{\triangle}fd\lambda<\int_{\triangle}gd\lambda$$
is true only on the last interval $\triangle_n.$ (If the inequality is false for all the intervals $\triangle_j,\ j=1,...,n$, then we can take $\bar{f}=f.$) Indeed, let $j_0$ be the largest natural number less than $n$ and such that
$\int_{\triangle_{j_0}}fd\lambda<\int_{\triangle_{j_0}}gd\lambda$. Consider the functions $f$ and $g$ on the interval $[a,t_{j_0}].$ On this interval $f$ and $g$ satisfy all the assumptions of Lemma 3.3 and by the induction assumption we can find a function $f_{j_0}$ defined on $[a,t_{j_0}]$, equimeasurable with the restriction of $f$ on this interval, and such that $$\int_{\triangle_i}f_{j_0}\geq \int_{\triangle_i}gd\lambda,\ i=1,...,j_0.$$
Let
$$\bar{f}_{j_0}(s):=\begin{cases}f_{j_0}(s), &\mbox{if }\ s\in[a,t_{j_0}];\\
f(s),\ &\mbox{if}\ s\in[t_{j_0}, b].\end{cases}$$
Then the function $\bar{f}_{j_0}$ is equimeasurable with  $f$ on $[a,b]$, it satisfies all the conditions of the current lemma, and the inequality $$\int_{\triangle}\bar{f}_{j_0}d\lambda<\int_{\triangle}gd\lambda$$
can be true only on the interval $\triangle_n.$\\

Thus, we assume that for the functions $f$ and $g$ we have the inequalities
$$\int_{\triangle_j}fd\lambda\geq\int_{\triangle_j}gd\lambda,\ j=1,...,n-1.\eqno(3.6)$$
Therefore we can assume that
 $$\int_{t_1}^b fd\lambda<\int_{t_1}^bgd\lambda\eqno(3.7).$$
Indeed, if $\int_{t_1}^b fd\lambda\geq\int_{t_1}^bgd\lambda$ then for the functions $f$ and $g$ considered on the interval $[t_1,b]$ all the conditions of the current lemma are satisfied (taking into consideration inequalities $(3.6)$). By the induction assumption there is a function $f_1$ defined on $[t_1,b]$, equimeasurable on this interval with $f$, and such that
$$\int_{\triangle_i}f_1d\lambda\geq\int_{\triangle_i}gd\lambda,\ j=2,...,n.$$
Then the function
$$\bar{f}(s):=\begin{cases}f(s), &\mbox{если }\ s\in[a,t_1];\\
f_1(s),&\mbox{если}\ s\in[t_1, b] \end{cases} $$
has the claimed properties.\\

Finally, we can assume that for the functions $f$ and $g$ are satisfied the inequalities ($3.6$) и $(3.7)$. But then we can apply Lemma 3.3 to prove the existence of the function $\bar{f}$.\\
$\Box$\\

\textbf{Lemma 3.5.} Let  $[a,b]$ be an arbitrary finite interval and let $f$ and $g$ be nonnegative integrable functions on $[a,b]$ such that
$$g=\sum_{j=1}^n\alpha_j\cdot 1_{\triangle_j},\ \textrm{где}\ \alpha_1\geq \alpha_2\geq...\geq \alpha_n>0.$$
$$\triangle_j=[t_{j-1},t_j],\ j=1,2,...,n; a=t_0<t_1<...<t_n=b,$$
and for every $t\in [a,b]$
$$\int_a^tfd\lambda>\int_a^tgd\lambda.$$
Then there is a function $\bar{f}$ defined on $[a,b]$, equimeasurable with $f$, and such that
$$\int_{\triangle_j}\bar{f}d\lambda\geq \int_{\triangle_j}gd\lambda,\ j=1,...,n.$$

\bigskip

\emph{Proof.} If the partition $\{\triangle_j\}_{j=1}^n$ of $[a,b]$ into intervals where $g$ is constant satisfies condition $(3.3)$ of Lemma 3.3 then the statement of the lemma follows directly from Lemma 3.4.\\

If the partition $\{\triangle_j\}_{j=1}^n$  is arbitrary $[a,b]$ we can always construct a finer partition $\{\triangle'_j:=[t'_{j-1},t'_j]\}_{j=1}^m$ such that we have the inequalities
$$t'_j\leq \frac{t'_{j-1}+t'_{j+1}}{2},\ j=1,...,m$$
the function $g$ is constant on the elementa of the new partition, and all the conditions of Lemma 3.4 are fulfilled. Therefore we can find a function $\bar{f}$ defined on $[a,b]$, equimeasurable with $f$, and such that
$$\int_{\triangle'_j}\bar{f}d\lambda\geq \int_{\triangle'_j}gd\lambda,\ j=1,...,m.$$
Finally, we notice that, because the new partition is finer than the old one, we have the inequalities
$$\int_{\triangle_j}\bar{f}d\lambda\geq \int_{\triangle_j}gd\lambda,\ j=1,...,n.$$

$\Box$\\

\textbf{Lemma 3.6.} Let $0\leq x, y\in L^1(I)$ and assume that for any $t\in I$
$$\int_0^txd\lambda>\int_0^tyd\lambda.$$
Then there is  a sequence $\{r_n\}^\infty_{n=1}\subseteq I,\ 1=r_1>r_2>...,\ \lim_{n\rightarrow\infty} r_n=0$ such that for any $n=1,2,...$ and for any $t\in(r_n,1]$ we have
$$\int_{r_n}^txd\lambda>\int_{r_n}^tyd\lambda.$$

\bigskip
\emph{Proof}. Assume to the contrary that a sequence with required properties does not exist. Let $1=r'_1>r'_2>...,\ \lim_{n\rightarrow\infty} r'_n=0$ be a sequence of points from $I$. It follows from our assumption that there is an index $N$ such that if $n\geq N$ then all the sets
$$\mathcal{D}_n=\{t\in (r'_n, 1]: \int_{r'_n}^txd\lambda\leq\int_{r'_n}^tyd\lambda$$
are not empty.\\

 Let $\delta_n:=\sup\mathcal{D}_n;$ clearly $\delta_n\in\mathcal{D}_n,\ n\geq N.$ Let us prove that
$$\delta:=\liminf\delta_n>0.$$
Indeed, otherwise we can find a decreasing to $0$ subsequence $\delta_{n_k}.$ By assumption it does not satisfy the property indicated in the statement of the lemma. Therefore we can find a natural $i$ and a point $\gamma\in(\delta_{n_i}, 1]$ such that
 $$\int_{\delta_{n_i}}^\gamma xd\lambda\leq\int_{\delta_{n_i}}^\gamma yd\lambda.$$
Fro there we obtain
$$\int_{r'_{n_i}}^\gamma x d\lambda=\int_{r'_{n_i}}^{\delta_{n_i}}x d\lambda+\int_{\delta_{n_i}}^\gamma x d\lambda\leq \int_{r'_{n_i}}^{\delta_{n_i}}y d\lambda+\int_{\delta_{n_i}}^\gamma y d\lambda=\int_{r'_{n_i}}^\gamma y d\lambda,$$
in contradiction with the definition of $\delta_{n_i}$.\\

Thus, $\delta>0;$ let $\delta_{m_k}\rightarrow\delta.$ By taking the limit in the inequality
$$\int_{r'_{n_k}}^{\delta_{m_k}}x d\lambda\leq\int_{r'_{n_k}}^{\delta_{m_k}}y d\lambda$$
we obtain, in contradiction with the conditions of the lemma, that
$$\int_0^\delta x d\lambda\leq\int_0^\delta y d\lambda.$$
$\Box$\\

\textbf{Lemma 3.7.} Let $x=x^*,y=y^*\in L^1(I)$, assume additionally that $y$ is a function with at most countable range and $\mathcal{F}$=measurable, where $\mathcal{F}$ is an interval partition. Assume that for any $t\in I$ we have
$$\int_0^txd\lambda\geq\int_0^tyd\lambda.$$
Then there is a function $\bar{x}$ defined on $I$, equimeasurable with $x$ and such that
$$y\leq E(\bar{x}|\mathcal{F}).$$

\bigskip
\emph{Proof.} By considering $x+ \varepsilon$, where $\varepsilon$ is an arbitrary small positive comstant, we can assume that the inequality in the statement of the lemma is strict. By Lemma 3.6 we can find a strongly decreasing to $0$ sequence $\{r_n\}_{n=0}^\infty,\ r_0=1$ of points from $(0,1]$  such that  for the restrictions $x_n$ and $y_n$ of the functions $x$ and $y$, respectively on the intervals $\triangle_n=[r_n,r_{n-1}]$ the inequality
$$\int_{r_n}^tx_nd\lambda>\int_{r_n}^ty_n d\lambda$$
is valid for any $t\in [r_n,r_{n-1}],\ n\geq 1.$ Denote by  $\mathcal{F}_n$ the partition of $\triangle_n$ into the intervals where the function $y_n$ is constant, and by $x_n$ the restriction of function $x$ on $\triangle_n,\ n=1,2,... .$ Now we can apply Lemma 3.5 to conclude that for any
 $ n=1,2,...$ there is a function $\bar{x}_n$ defined on $\triangle_n$, equimeasurable with $x_n$, and such that
$$y_n\leq E(\bar{x}_n|\mathcal{F}_n).\eqno(3.8)$$
Next we define on $[0,1]$ the function $\bar{x}$ by taking $\bar{x}(0)=\bar{x}(1)=0,\ \bar{x}(s)=\bar{x}_n(s),\ \textrm{if}\ s\in \triangle_n,\ n=1,2,... .$ Let
$$\mathcal{\bar{F}}=\{F: F\in \mathcal{F}_n\ \textrm{for}\ \textrm{some}\ n=1,2,...\} .$$
Clearly the partition $\mathcal{\bar{F}}$ of the interval $[0,1]$ is finer than the partition $\mathcal{F}$. It follows from ($3.8$) that $y\leq E(\bar{x}|\mathcal{\bar{F}})$ and therefore $y\leq E(\bar{x}|\mathcal{F}).$ It remains to notice that $\bar{x}$ and $x$ are equimeasurable on $I$.\\
$\Box$\\

\bigskip

We return to the proof of Theorem 3.2. Without loss of generality and like in the proof of Lemma 3.7 we can assume that $x=x^*,\ y=y^*$  and for any  $t\in I$ we have
$$\int_0^txd\lambda<\int_0^tyd\lambda.$$
Let us fix an $\varepsilon>0$.
We can find  a $y_{\varepsilon} = y^*_{\varepsilon}$ with at most countable range and such that
$$\textbf{0}\leq y_{\varepsilon} \leq y\leq y_{\varepsilon} +\varepsilon\cdot \textbf{1}.\eqno(3.9)$$
Clearly the functions $y_{\varepsilon}$ and $x$ satisfy all the conditions of Lemma 3.7 and therefore we can find a function $x_{\varepsilon}$ equimeasurable with $x$ and such that
$$y_{\varepsilon}\leq E(x_{\varepsilon}|\mathcal{F}_{\varepsilon}),\eqno(3.10)$$
where $\mathcal{F}_{\varepsilon}$ is the partition of $I$ generated by all the intervals on which the function $y_{\varepsilon}$ is constant.\\

From inequalities $(3.9)$ and $(3.10)$ follow the inequalities
$$0\leq y\leq E(x_{\varepsilon}|\mathcal{F}_{\varepsilon})+\varepsilon\cdot \textbf{1}.\eqno(3.11)$$
By Theorem 3.1 $X$ is a symmetric ideal and therefore $x_{\varepsilon}\in X$. Because, by assumption,  $E(x_{\varepsilon}|\mathcal{F}_{\varepsilon})\in X,$ we have $y\in X$ and therefore part 1 of Theorem  3.2 has been proved.\\
$\Box$\\

2. Proof of part 2. Let $(X,|| \cdot ||_X)$ is a Banach ideal vector subspace of $L^1(I)$ such that $|| E(\cdot |\mathcal{F}||_X\leq 1$ for any at most countable partition $\mathcal{F}$. From the relation $y\prec x$ and inequalities $(3.9) - (3.11)$  follows that
$$||y||_X\leq ||E(x_{\varepsilon}|\mathcal{F}_{\varepsilon})+\varepsilon\cdot \textbf{1}||_X\leq ||E(x_{\varepsilon}|\mathcal{F}_{\varepsilon})||_X\cdot||x_{\varepsilon}||_X+ ||\varepsilon\cdot\textbf{1}||_X\leq ||x_{\varepsilon}||_X+\varepsilon\cdot||\textbf{1}||_X.$$
Because $\varepsilon$ is arbitrary small the second part of Theorem 3.2 follows.\\
$\Box$\\

\textbf{Theorem 3.8.} Let $(X,\|\cdot \|_X)$ be a Banach ideal such that $\|E(\cdot |\mathcal{F})\|_X\leq 1$ for any at most countable partition $\mathcal{F}$. Then on $X$ exists an equivalent symmetric norm $||\cdot||^1$ such taht  $(X,||\cdot||^1)$ is a strongly majorant ideal.\\

\bigskip

\emph{Proof.} First we will prove that there is a constant $K>0$ such that
$$x,\ y\in X,\ y\prec x\Rightarrow || y||_X\leq K\cdot|| x||_X.$$
Assume to the contrary that such a constant does not exist. Notice that for any $r>0$ we have
 $y\prec x \Rightarrow r\cdot y\prec r\cdot x$ and therefore we can find two sequences $\{x_n\}_{n=1}^\infty$ and $\{y_n\}_{n=1}^\infty$ of elements of $X$ such that
 $$y_n=y^*_n\prec x_n=x^*_n,\ || x||_X\leq 2^{-n},\ || y_n||_X\geq n,\ n=1,2,... .$$
Let $x_0:=\sum_{n=1}^\infty x_n$ where the series converges almost everywhere.  Because the space $(X,|| \cdot ||)$ is norm complete we obtain that $x_0\in X$. It follows by the majorant convergence theorem that $y_0:=\sum_{n=1}^\infty y_n\in L^1$ and $y_0\prec x_0$. By Theorem 3.2 $y_0\in X.$\\

On the other hand $\textbf{0}\leq y_n\leq y_0,\ n\geq 1$ whence $||y_0||_X\geq ||y_n||_X\geq n,\ n\geq 1$, a contradiction.\\

Let us continue with the proof of Theorem 3.8. For any $x\in X$ let
$$||x||^1:=\sup\{||y||_X:\  y\prec x.\}$$
As we have just proved $||x||^1\leq K\cdot ||x||_X$ for some $K>0$. We will show that for $||\cdot||^1$ we have the triangle inequality. Let $x_1,\ x_2\in X,\ \varepsilon>0$. By definition there is $y\in X$, such that $y\prec x_1 +x_2$ and $||x_1 + x_2||^1\leq ||y||_X + \varepsilon$. According to~\cite{BS} there are $y_1,\ y_2\in L^1$ such that $y_i\prec x_i,\ i=1,2$ and $y=y_1+y_2$. Combining it with Theorem 3.2 we obtain that $y_i\in X,\ i=1,2$. Next we have
$$||x_1+x_2||^1\leq ||y||_X + \varepsilon\leq ||y_1||_X+||y_2||_X + \varepsilon \leq ||x_1||^1 +||x_2||^1+\varepsilon,$$
where $\varepsilon$ is arbitrary small. That proves the triangle inequality for $||\cdot ||^1$.\\

It is now trivial to verify that
 $||\cdot ||^1$ is a monotonic symmetric norm on $X$. Because for any $x\in X$ we have $x\prec x$ therefore $||x||_X\leq ||x ||^1$ and both norms are equivalent.\\
$\Box$\\

\textbf{ 3.9.} The fact that any $\sigma$-subalgebra in $\Lambda$ strongly averages the Banach ideal $(X, ||\cdot||_X)$ implies that $X$ is a symmetric vector ideal though the original norm $||\cdot||_X$ might be not symmetric. Indeed, let $(X, ||\cdot||_X)$ be a strongly interpolation space. For any $x\in X$ let
$$||x||^*:=||x||_X+\sup\{\int_I|x|\textbf{P}_{i=1}^nE(1_{(0,\frac{1}{2}]}|\mathcal{A}_i)d\lambda\},$$
where the supremum is taken over all $\sigma$-subalgebras $\mathcal{A}_i$ in $\Lambda,\ i=1,2,...,n;\ n=1,2,...$ and by $\textbf{P}_{i=1}^nE(f|\mathcal{A}_i)$ we denote $E(E(...E(f|\mathcal{A}_n))$. It is not difficult to verify that $||\cdot||^*$ is a monotonic Banach norm on $X$ equivalent to $||\cdot||_X$ and it follows from the properties of averaging operators that $||E(\cdot|\mathcal{A})||^*=1$ for any $\sigma$-subalgebra  $\mathcal{A}$ in $\Lambda.$ It is also clear that $||1_{(0,\frac{1}{2}]}||^*>||1_{(\frac{1}{2},1]}||^*.$\\
$\Box$\\

\bigskip

\newpage

\bigskip

\centerline{\textbf{Chapter 4. On averaging principal symmetric ideals by countable partitions.}}
\bigskip

\centerline{\textbf{The results of this chapter are based on papers~\cite{BM},\ ~\cite{Me14}.}}

\bigskip

\emph{In this Chapter we discuss when, given a function $f$, a countable partition  averages the principal symmetric ideal $\mathcal{N}_f$ or the symmetric space $\mathcal{M}_f^1$}.\\


\bigskip

First of all notice that without loss of generality we can consider only countable interval partitions
(see Remark 0.11.2)\\

Thus, in the sequel we will consider an interval partition $\mathcal{B}=(b_n)$, that belongs to the stochastic vector $\vec{\beta} =(\beta_n)$, and the principal symmetric ideal $\mathcal{N}_f$  generated by the function $f,\ f=f^*\in L^1(I)$.\\

\textbf{Theorem 4.1.} If an interval partition $\mathcal{B}$ averages the principal symmetric ideal $\mathcal{N}_f$ then for some constant $Q>1$ we have the inequality
$$f^{**}(b_n)\leq Q\cdot f(Q^{-1}\cdot b_n),\ n\geq 1.\eqno(4.1)$$

Moreover, if the partition $\mathcal{B}$ is monotonic then condition $(4.1)$ is not only necessary but also sufficient for $\mathcal{B}$ to average $\mathcal{N}_f$
\bigskip

The proof of Theorem 4.1 will be based on some axillary results.

\bigskip

\textbf{ Proposition 4.2.} If interval partition $\mathcal{B}$ is monotonic and for some constant $Q>1$ we have the inequalities
$$f^{**}(b_n)\leq Qf(b_n),\ n\geq 1,\eqno(4.2)$$
then the partition $\mathcal{B}$ averages the principal symmetric ideal $\mathcal{N}_f$.\\

The proof of Proposition 4.2 will be obtained in several steps which for convenience we single out as separate lemmas.\\

\textbf{Lemma 4.3.} $E(f|\mathcal{B})\in \mathcal{N}_f.$\\

\textit{Proof}. For a fixed $n,\  n\geq 1,$ for any $t,\ t\in (b_n,b_{n-1}],$ we have
$$E(f|\mathcal{B})(t)=\frac{1}{b_{n-1}-b_n}\int_{b_n}^{b_{n-1}}fd\lambda=\frac{1}{b_{n-1}-b_n}[b_{n-1}f^{**}(b_{n-1})-b_nf^{**}(b_n)].\eqno (4.3)$$
If $\frac{b_{n-1}}{b_{n-1}-b_n}\leq 2,$ then we can continue $(4.3)$ as
$$...\leq \frac{b_{n-1}}{b_{n-1}-b_n}f^{**}(b_{n-1})\leq 2 f^{**}(b_{n-1})\leq 2Qf(b_{n-1})\leq 2Qf(t).$$
On the other hand, if $\frac{b_{n-1}}{b_{n-1}-b_n}> 2$ then by Lemma 0.1 we can continue (4.3) as
$$...\leq\frac{1}{b_{n-1}-b_n}[b_{n-1}f^{**}(b_n)-b_nf^{**}(b_n)]=f^{**}(b_n)\leq f^{**}(2^{-1}\cdot b_{n-1})\leq 2 f^{**}(b_{n-1})\leq 2Q f(b_{n-1})$$
$$\leq 2Q f(t).$$
Therefore, $E(f|\mathcal{B})(t)\leq 2Q f(t),\ t\in I.$\\
$\Box$\\

\textbf{Remark 4.4.} The statement of Lemma 4.3 remains true without the assumption that the partition $\mathcal{B}$ is monotonic.\\

In the sequel we will always assume that the function $f$ is $ \mathcal{D}_2$-measurable. Relations $(0.4)$ show that this assumption does not result in loss of generality.\\

\textbf{Lemma 4.5.} For any subsequence $\bar{\beta}=\{\beta_{k_i}\}_{i\geq 1}$ of coordinates of the stochastic vector $\vec{\beta}=(\beta_n)$ we can find a natural $N_{\bar{\beta}}$ such that if
 $n\geq N_{\bar{\beta}}$ then we have
$$f^{**}(\sum_{i\geq n}\beta_{k_i})\leq 4\cdot Q\cdot f(\sum_{i\geq n}\beta_{k_i}),\ n\geq N_{\bar{\beta}}. \eqno(4.4)$$
\emph{Proof.} Let $\vec{\beta}^* =(\beta_m^*)$ is the stochastic vector of the monotonic interval partition $\mathcal{B}^*:=(b_m^*).$ For an interval partition $\mathcal{B}_1$ that belongs to the stochastic vector $\tilde{\beta}=(\beta_{k_i}, 1-\sum_{n\neq k_i,\ i\geq 1}\beta_n)$ and is a sample from $\bar{\beta}$, we can find another interval partition $\mathcal{B}:=(b_n)\sim \mathcal{B}^*$ such that the points of $\mathcal{B}_1=(\sum_{i\geq n}\beta_{k_i})$ coincide with the points $b_n$ if $n\geq N_{\bar{\beta}},$ where the number $N_{\bar{\beta}}$ is defined by the sample. But the monotic rearrangement of $\mathcal{B}$ is the interval partition $\mathcal{B}^*$, and therefore we can apply Theorem 2.1. Because $2>\frac{5^{1/2}+1}{2}$, by assuming that in $(2.1)\ \delta=2,$ and applying inequalities $(0.2)$ and  $(4.2)$ we obtain that
$$f^{**}(\sum_{i\geq n}\beta_{k_i})=f^{**}(b_n)\leq f^{**}(2^{-1}\cdot b^*_m)\leq 2\cdot f^{**}( b^*_m)\leq 2\cdot Q\cdot f(b^*_m)\leq$$
$$\leq 4\cdot Q\cdot f(b_n)=4\cdot Q\cdot f(\sum_{i\geq n}\beta_{k_i}),\ n\geq N_{\bar{\beta}}.$$
$\Box$\\

\textbf{Corollary 4.6.} Assume that the interval partition $\mathcal{S},\ \mathcal{S}=\sigma(s_n,\ n\geq 0),$ belongs to the stochastic vector $\vec{\mathfrak{b}}$ that is coarser than $\vec{\beta}=(\beta_n)$. Then for some constant $Q\geq 1$ we have
$$f^{**}(s_n)\leq 4Qf(s_n),\ n\geq 0.\eqno(4.5)$$
\emph{Proof.} The proof follows directly from Lemma 4.5 because every point $s_n$ of the partition $\mathcal{S}$ is equal to $\sum_{k\geq 1}\beta_{n_k}$ for an appropriate sequence $\{\beta_{n_k}\}$ of coordinates of the stochastic vector $\vec{\beta}$.\\
$\Box$\\

\textbf{Lemma 4.7.} Let $\mathcal{F}=\sigma(F_n,\ n\geq 1)$ be a partition of $I$ that belongs to a rearrangement  $\vec{\gamma}$ of the stochastic vector $\vec{\beta}$. Then $E(f|\mathcal{F})\in \mathcal{N}_f.$\\

\emph{Proof.} It is enough to prove that $h:=[E(f|\mathcal{F})]^*\in \mathcal{N}_f.$ Let $\mathcal{S}=(s_n)$ be an interval partition such that elementary function $h=h^*$ is $\mathcal{S}$-measurable . If $\mathcal{S}\in \vec{\alpha}$ then the stochastic vector $\vec{\alpha}$ is coarser than the stochastic vector  $\vec{\gamma}$ and therefore the stochastic vector $\vec{\alpha}$ is coarser than the stochastic vector $\vec{\beta}$. Now we can apply Corollary 4.6 from which follows inequality $(4.5)$. Finally, it follows from Lemma 4.3 and Remark 4.4 that the function $H:=E(f|\mathcal{S})$ is in $\mathcal{N}_f.$\\

On the other hand, it follows from Proposition 0.4 that $h=E(\tilde{f}|\mathcal{S})$ where $\tilde{f}=f\circ\pi$ is a function equimeasurable with $f$ and $\pi$ is an automorphism of the interval $I$ such that $\mathcal{S}=\mathcal{F}\circ\pi$ (see~\cite{Ro} for the proof of existence of such an automorphism). Applying the properties of averaging operators, the inequality from Lemma 0.1 for $f^{**}$, the monotonicity of $f$, and inequality $(4.5)$ we obtain that for any $n,\ n\geq 1$ and for  $t\in(s_n,s_{n-1}]$ we have
$$\frac{h(t)}{H(t)}=\frac{h(s_{n-1})}{H(s_{n-1})}\leq \frac{h^{**}(s_{n-1})}{H(s_{n-1})}=\frac{s^{-1}_{n-1}\int_0^{s_n^{-1}}E(\tilde{f}|\mathcal{S})d\lambda}{H(s_{n-1})}=
\frac{s^{-1}_{n-1}\int_0^{s_n^{-1}}\tilde{f}d\lambda}{H(s_{n-1})}\leq$$
$$\leq\frac{f^{**}(s_n^{-1})}{H(s_{n-1})}\leq\frac{f^{**}(s_{n^-1})}{f(s_{n^-1})}\leq 4Q.$$
Therefore, $h\in \mathcal{N}_f.$\\
$\Box$\\

\textit{Proof of Proposition} 4.2. From Lemma 4.7 (applying again an appropriate automorphism of $I$) we can easily obtain that
$$E(g|\mathcal{B})\in \mathcal{N}_f\eqno(4.6)$$
for any elementary function $g$ such that
$$g^*\leq Cf,\ C>0.\eqno(4.7)$$
Finally, using uniform approximation by elementary functions and positivity of the averaging operator we can prove that inclusion $(4.6)$ remains valid for any function $g$ that satisfies inequality $(4.7)$ with an appropriate positive constant $C$.

$\Box$\\

\textit{Proof of Theorem 4.1} To prove the sufficiency of condition (4.1) in Theorem 4.1 we will define for any $m,\ m\geq 1$ a $\mathcal{D}$-measurable non-increasing function $f_{(m)},\ f_{(m)}(t)=f(2^{-m}t),\ t\in I.$ For any function $f_{(m)}$ condition $(4.2)$ of Proposition 4.2 is satisfied with a constant $Q_m,\  Q_m\leq 2^mQ$ and therefore we can apply to it our previous reasoning. It remains to notice that $\mathcal{N}_f=\mathcal{N}_{f_{(m)}},\ m\geq 1$.\\

The necessity in Theorem 4.1 follows directly from Lemma 0.4 and Proposition 4.8 below.\\

\textbf{Proposition 4.8.} Assume that for the function $f,\ f=f^*\in L^1(I)$, for the interval partition $\mathcal{B},\ \mathcal{B}=(b_n)$ (we do not assume here that this interval partition is monotonic), and for any $m\geq 1$ we have
$$\sup_{n\geq 0}\frac{f^{**}(b_n)}{f(\frac{b_n}{m})}=\infty.\eqno(4.8)$$
Then there is an interval partition $\mathcal{V},\ \mathcal{V}=(v_k)$ that is coarser than $\mathcal{B}$ and such that $E(f|\mathcal{V})\notin\mathcal{N}_f$.\\

\emph{Proof.} It ckearly follows from $(4.8)$ that there is a subsequence $\{u_m:=b_{n_m},\ n_0=0\}$ of the sequence $\{b_n\}$ such that
$$\lim_{m\rightarrow\infty}\frac{f^{**}(u_m)}{f(\frac{u_m}{m})}=\infty.\eqno(4.9)$$
We will define a subsequence $\{v_k\}$ of the sequence $\{u_m\}$ by induction. Let $v_0=b_{n_0}=u_{m_0}=1$ and assume that we have already constructed the point $v_k=u_{m_k},\ k\geq 0$. Because the function $f^{**}(t),\ t\in I$ is continuous and monotonic in some neighborhood of $0$ and by using $(4.9)$ we can chose a point $u_{m_{k+1}}$ so close to zero that $$0<u_{m_{k+1}}<\frac{v_k}{k}\ \textrm{ and } \frac{\frac{\int_{u_{m_{k+1}}}^{v_k}}{v_k-u_{m_{k+1}}}}{f(\frac{v_k}{k})}>\frac{f^{**}(v_k)}{f(\frac{v_k}{k})}-1.$$ We finish the induction by letting $v_{k+1}:=u_{m_{k+1}}.$ We obtain that $$\lim_{j\rightarrow\infty}\frac{E(f|\mathcal{V})(v_k)}{f(\frac{v_k}{k})}=\infty$$
and therefore $E(f|\mathcal{V})\notin\mathcal{N}_f.$\\

Поскольку и.р. $\mathcal{V}$ крупнее, чем и.р. $\mathcal{B}$, то по лемме 0.11.5) справедливо: $E(f|\mathcal{B})\notin\mathcal{N}_f$, т.е. необходимость в теореме 4.1 доказана.\\
$\Box$\\




\textbf{Remark 4.9.} Similar to the proof of Proposition 4.8 we can show that if a countable partition $\sigma\Big((c^{(n)}_k,c^{(n)}_{k-1}]\Big)_{n,k\geq 1}$ averages the principal symmetric ideal $\mathcal{N}_f$ then
$$\sup_{n\geq 1,k\geq 0}\frac{f^{**}(c_k^{(n)})}{f(c_k^{(n)})}<\infty.$$

Let $\mathcal{F}$ be a countable partition. A function $f\in L^1(I)$ is called $\mathcal{F}$-\emph{regular} if the monotonic interval partition $\mathcal{F}^*$ averages the principal symmetric ideal $\mathcal{N}_f.$\\

Let $\mathcal{B}=(b_n)$ be an interval partition and let  $f,\ f=f^*\in L^1(I)$. Then we introduce
$$f_{\mathcal{B}}:=\sum_{n\geq 1}f(b_{n-1})\cdot 1_{(b_n,b_{n-1}]}\ ;\ f^{**}_{\mathcal{B}}:=\sum_{n\geq 1}f^{**}(b_n)\cdot 1_{(b_n,b_{n-1}]}.$$
Now we can state Theorem 4.1 in the following way\\

\textbf{Theorem 4.10.} For an interval partition $\mathcal{B}$ and a function $f,\ f=f^*\in L^1(I)$ the following three conditions are equivalent\\

1). The function $f$ is $\mathcal{B}$-regular; \\

2). $f^{**}_{\mathcal{B}^*}\simeq f^*_{\mathcal{B}^*}$;\\

3). $f^{**}_{\mathcal{B}^*}\in\mathcal{N}_f.$\\

It follows from (0.11.3) that these conditions are also equivalent to the following one\\

4). $f$ is $\mathcal{B}^*_{(2)}$-regular.\\

\bigskip

\textbf{Lemma 4.11.}  Let $f=f^*\in L^1(I)$ and let $\mathcal{F}$ be a countable partition.
The following conditions are equivalent

(a) $f$ is $\mathcal{F}$-regular.\\
(b) $$ E(\mathcal{M}_f|\mathcal{F}^*)\subseteq \mathcal{N}_f.\eqno(4.10)$$

\emph{Proof.} The implication $(b) \Rightarrow (a)$ is obvious.

Assume $(a)$. Then $E(\mathcal{N}_f|\mathcal{F}^*)\subseteq\mathcal{N}_f$ and by Theorem 5.6 \footnote{Here we refer the reader to a result that will appear only in the next Chapter, but the proof of Theorem 5.6 does not depend on the current lemma.}
$$\mathcal{N}_{E(\mathcal{M}_f|\mathcal{F}^*)}=\mathcal{N}_{E(\mathcal{N}_f|\mathcal{F}^*)}\subseteq\mathcal{N}_{\mathcal{N}_f}=
\mathcal{N}_f.$$
$\Box$\\

Let $ f=f^*\in L^1(I)$ and let $\textsl{M}_f$ be the corresponding Marcinkiewicz space. Let also $\textsl{M}_f^1$ be the closure in the norm $\|\cdot \|_{\textsl{M}_f}$ of the symmetric ideal $\mathcal{N}_f$ in the symmetric Banach space $\textsl{M}_f$.\\

\bigskip

\textbf{Theorem 4.12.} For an interval partition $\mathcal{B}$, $\mathcal{B}= \sigma(B_n),\ B_n=(b_n, b_{n-1}],\ n\geq 1,\ b_0=1$ the following conditions are equivalent.\\

$$1)\ \mathcal{B}\ \textrm{averages the symmetric ideal }\ \mathcal{N}_f.$$
$$ 2)\ \mathcal{B}\ \textrm{averages the symmetric Banach space}\ \mathcal{M}^1_f.$$.

\bigskip

\emph{Proof}.  1)$\Rightarrow$ 2). Let $g\in \textsl{M}_f^1$, i.e. there is a sequence $\{g_n\}_{n\geq 1}\subseteq \mathcal{N}_f$ such that $\|g_n-g\|_{\textsl{M}_f}\rightarrow_{n\rightarrow\infty}0.$ Because the double stochastic projection $E(\cdot|\mathcal{B})$ is a continuous operator on the Banach space $(\textsl{M}_f,\ \|\cdot\|_{\textsl{M}_f})$ we have
$$\lim_{n\rightarrow\infty} \|E(g|\mathcal{B})- E(g_n|\mathcal{B})\|_{\textsl{M}_f}=0.$$ In virtue of our conditions $E(g_n|\mathcal{B})\in \mathcal{N}_f,\ n\geq 1,$ and therefore $E(g|\mathcal{B})\in \textsl{M}_f^1$.\\

To prove the implication 2)$\Rightarrow$ 1) we will show that if $\mathcal{B}$ does not average the symmetric ideal $\mathcal{N}_f$ then $\mathcal{B}$ also does not average the symmetric space $M^1_f$. (What follows is a modification of the proof of Theorem 4.1, see also~\cite{BM}).\\

Assume that $\mathcal{B}$ does not average $\mathcal{N}_f$ i.e. despite equality $(4.1)$ we have
$$\limsup_{n\rightarrow\infty}\frac{f^{**}(b_n)}{f(b_n)}=\infty.$$
Without loss of generality we can assume that $b_0=1,\ b_1=\frac{1}{2}.$ Let us denote by $\{s_m\}_{m\geq 0}:=\{b_{n_m}\}$ a subsequence of the sequence  $\{b_n\}$ such that $1=s_0>\frac{1}{2}=s_1>s_2>...>s_m>...;\ s_m\downarrow 0$ and also
$$\frac{f^{**}(\frac{1}{2})}{f(1)}=:a_1<a_m:=\frac{f^{**}(s_m)}{f(\frac{s_m}{m})} \uparrow_{m\uparrow\infty}\infty.$$
Let us fix a $\varepsilon>0$ and let $m_1=1,t_1=2s_{m_1}.$ On the interval $[0,\frac{1}{2}]$ we define the following function
$$\Phi_1(v)=[f(\frac{s_{m_1}}{m_1})(s_{m_1}-v)]^{-1}\int_v^{s_{m_1}}f(u)du.$$
Clearly the function $\Phi_1(v)$ is continuous, decreasing, and $\Phi_1(0)=a_{m_1}.$ There is an index $m_2$ such that
$$\ s_{m_2}\leq \frac{s_{m_1}}{4m_1},\  \Phi_1(2s_{m_2})\geq a_{m_1}-\varepsilon.$$
Let $t_2=2s_{m_2}.$ Then
$$a_{m_1}-\varepsilon\leq \Phi_1(t_2)=[f(\frac{s_{m_1}}{m_1})(s_{m_1}-t_2)]^{-1}\int_{t_2}^{s_{m_1}}f(u)du\leq [f(\frac{s_{m_1}}{m_1})(s_{m_1}-t_2)]^{-1}\int_{t_2}^{t_1}f(u)du.$$
Similarly we can construct indices $\{m_i\}_{i=1}^k:\ 1=t_1=2s_{m_1}>s_{m_1}>t_2=2 s_{m_2}>...>s_{m_k}.$\\
If $\Phi_k:=[f(\frac{s_{m_k}}{m_k})(s_{m_k}-v)]^{-1}\int_v^{s_{m_k}}f(u)du,\ 0\leq v\leq s_{m_k}$ then the function $\Phi_k(v)$ is defined on $[0,s_{m_k}]$, continuous, decreasing, and $\Phi_k(0)=a_{m_k}.$ Let $t_k:=2s_{m_k}$. By reasoning as above we can find an index $m_{k+1}$ such that
$$i)\ s_{m_{k+1}}\leq \frac{s_{m_k}}{4m_k};\ ii)\ s_{m_k}-t_{k+1}<t_k-s_{m_k}\Leftrightarrow s_{m_k}<\frac{t_k+t_{k+1}}{2},$$ где $t_{k+1}:=2s_{m_{k+1}};$
$$iii)\ a_{m_k}- \varepsilon\leq\Phi_k(s_{m_{k+1}})\leq [f(\frac{s_{m_k}}{m_k})(s_{m_k}-s_{m_{k+1}})]^{-1}\int_{s_{m_{k+1}}}^{s_{m_k}}f(u)du$$
$$\leq [f(\frac{s_{m_k}}{m_k})(s_{m_k}-t_{k+1})]^{-1}\int_{t_{k+1}}^{t_k}f(u)du.$$
We are done with the construction of the sequence of indices $\{m_k\}$ .\\

It follows from $i)$ that $m_k^{-1}s_{m_k}>t_{k+1}$. Therefore, in virtue of $iii)$ we have $m_k^{-1}s_{m_k}\in (t_{k+1},t_k),\ k\geq 1.$

Also,
$$(s_{m_k}-t_{k+1})^{-1}\leq 3(t_k-t_{k+1})^{-1},\ k\geq 1.\eqno(4.11)$$
Indeed, $s_{m_k+1}\leq \frac{s_{m_k}}{4m_k}\leq \frac{s_{m_k}}{4},\ k\geq 1,$ откуда $2s_{m_k}-2s_{m_{k+1}}
\leq 3s_{m_k}-6s_{m_{k+1}},$ или $t_k-t_{k+1}\leq 3(s_{m_k}-t_{k+1}).$\\
We will now define the interval partition $\mathcal{T}=\sigma(T_k)$ where $T_k=(t_{k+1},t_k],\ k\geq 1,$ and the averaging operator
$$E(f|\mathcal{T})=\sum_{k\geq 1}\alpha_k 1_{(t_{k+1},t_k]},\ \textrm{where}\ \alpha_k=\frac{1}{t_k-t_{k+1}}\int_{t_{k+1}}^{t_k}fdu,\ k\geq 1.$$
It follows from the construction of indices $s_{m_k},\ k\geq 1$ that
$$c_k:=\frac{\alpha_k}{f(\frac{s_{m_k}}{m_k})}\rightarrow\infty.\eqno(4.12)$$
We need to show that $inf\{\|E(f|\mathcal{T})-y\|_{\mathcal{M}_f}:\ y\in \mathcal{N}_f\}>0.$ In virtue of Lemma 0.16 and the equality $E(f|\mathcal{T})=(E(f|\mathcal{T}))^*$ it is enough to prove that $inf\{\|E(f|\mathcal{T})-y\|_{\mathcal{M}_f}:\ y=y^*\in \mathcal{N}_f\}>0.$\\

Let us fix an $y=y^*\in \mathcal{N}_f$ and let $k(y)$ be the smallest of all indices $k$ such that
$$vraisup_{t\in I}\frac{y(t)}{f(\frac{t}{m_k})}<\infty, $$
where $\{m_k\}$ is the sequence of indices constructed above. Let
$$c:=vraisup_{t\in I}\frac{y(t)}{f(\frac{t}{m_{k(y)}})};$$
$$H(c)=\{t\in I:\ \frac{E(f|\mathcal{T})(t)}{f(\frac{t}{m_{k(y)}})}\geq 2c;$$
$$\tau_k=\lambda T_k,\ \eta_k=\lambda (T_k\cap H(c)),\ k\geq 1.$$
Let us now prove that there is an index $k_1$ such that
$$\frac{\eta_k}{\tau_k}\geq\frac{1}{2}.\eqno(4.13)$$
if $k>k_1$. Indeed, it follows from $(4.11)$ that there is an index $k_0$ that for $k>k_0$
$$c_k=\frac{\alpha_k}{f(\frac{s_{m_k}}{m_k})}
=\frac{E(f|\mathcal{T})(s_{m_k})}{f(\frac{s_{m_k}}{m_k})}\geq 2c.$$
Therefore if $k>k_1:=\max[k_0,k(y)]$ we have
$$\frac{E(f|\mathcal{T})(s_{m_k})}{f(\frac{s_{m_k}}{m_{k(y)}})}\geq \frac{E(f|\mathcal{T})(s_{m_k})}{f(\frac{s_{m_k}}{m_k})}\geq 2c.$$
Thus, $s_{m_k}\in H(c).$ From the monotonicity of $f$ we conclude that $[s_{m_k},t_k]\subset H(c)$ and therefore
$$\frac{\eta_k}{\tau_k}\geq \frac{t_k-s_{m_k}}{t_k-t_{k+1}}\geq \frac{1}{2}$$
 $k>k_1$. That proves inequality $(4.13)$.\\

We continue with the proof of the theorem. For any $t\in H(c)$ we have
$E(f|\mathcal{T})(t)-y(t)\geq E(f|\mathcal{T})(t)- cf(\frac{t}{n_{k(y)}})\geq E(f|\mathcal{T})(t)-\frac{E(f|\mathcal{T})(t)}{2}=\frac{E(f|\mathcal{T})(t)}{2}.$ Therefore,
$$\|E(f|\mathcal{T})-y\|_{M_f}=\sup_{0<t\leq 1}\frac{\int_0^t(E(f|\mathcal{T})-y)^*(u)du}{\int_0^tf(u)du}\geq\frac{1}{2}\sup_{0<t\leq 1}\frac{\int_{(0,t)\cap H(c)}E(f|\mathcal{T})d\lambda}{\int_0^tfd\lambda}.\eqno(4.14)$$
Let $k_1$ be an index such  that $(4.13)$ is fulfilled. Because $\eta_k\geq\frac{\tau_k}{2}$ for $k\geq k_1$ it follows that
$$\int_{(0,t_{k_1})\cap H(c)}E(f|\mathcal{T})d\lambda=\sum_{k=k_1}^\infty\int_{T_k\cap H(c)}E(f|\mathcal{T})d\lambda=\sum_{k\geq k_1}\alpha_k\eta_k\geq$$
$$\geq \frac{1}{2}\sum_{k=k_1}^\infty\tau_k\frac{\int_{T_k}E(f|\mathcal{T})d\lambda}{\tau_k}=\frac{1}{2}\int_0^{t_{k_1}}f(u)du.$$
By continuing inequality  $(4.14)$ we get
$$\|E(f|\mathcal{T})-y\|_{M_f}\geq\frac{1}{2}\sup_{0<t\leq 1}\frac{\int_{(0,t)\cap H(c)}E(f|\mathcal{T})d\lambda}{\int_0^tfd\lambda}\geq\frac{1}{2}\frac{\int_{(0,t_{k_1})\cap H(c)}E(f|\mathcal{T})d\lambda}{\int_0^{t_{k_1}}f(u)du}\geq\frac{1}{4}.$$
Thus, $\inf\{\|E(f|\mathcal{T})-y\|_{M_f}:\ y=y^*\in \mathcal{N}_f\}\geq\frac{1}{4}.$\\
Therefore the interval partition $\mathcal{T}$ does not average the symmetric space $M_f^1.$ Now notice that $\mathcal{T}$ is a multiple of the interval partition $\mathcal{S}:=\sigma(s_{m_k})$ and the later is coarser than $\mathcal{B}.$ Applying Remark 0.11;3), 5) we see that the interval partition $\mathcal{B}$ also does not average the symmetric space $M_f^1.$\\
$\Box$
\newpage
\bigskip

\centerline{\textbf{Chapter 5. On symmetric ideals generated by averaging}}

\centerline{\textbf{ of principal symmetric ideals over countable partitions.}}
\bigskip

\centerline{\textbf{The results of this chapter are based on papers~\cite{Me1},~\cite{Me2},~\cite{Me4},~\cite{Me15},\ ~\cite{Me16}}}

\bigskip

\emph{Theorems 5.10 - 5.12 provide conditions on an interval partition $\mathcal{B}$ and a function $f$  suffucient for the symmetric ideal $\mathcal{N}_{E(\mathcal{N}_f|\mathcal{B})}$ generated by the image of a principal symmetric ideal $\mathcal{N}_f$ under the action of the $\mathcal{B}$-averaging operator, is itself a principal symmetric ideal.}\\

\bigskip

\bigskip

Let $\mathcal{F}_1$ and $\mathcal{F}_2$ be at most countable partitions of $I$,  and let $\mathcal{B}$, $\mathcal{B}_1$, and  $\mathcal{B}_2$ be interval partitions. In the sequel we will need the following simple lemma.\\

\textbf{Lemma 5.1.} Let $X,\ X\subseteq L^1(I),$ be a symmetric ideal. Then\\

1). If $\mathcal{F}_1\subseteq\mathcal{F}_2$ then

$\mathcal{N}_{E(X|\mathcal{F}_1)}\subseteq\mathcal{N}_{E(X|\mathcal{F}_2)}$  and

 $\mathcal{M}_{E(X|\mathcal{F}_1)}\subseteq\mathcal{M}_{E(X|\mathcal{F}_2)}$;\\

2). If $\mathcal{F}_1\sim\mathcal{F}_2$ then

 $\mathcal{N}_{E(X|\mathcal{F}_1)}=\mathcal{N}_{E(X|\mathcal{F}_2)}$ and

 $\mathcal{M}_{E(X|\mathcal{F}_1)}=\mathcal{M}_{E(X|\mathcal{F}_2)}$;\\

3). If $\mathcal{B}_1\simeq\mathcal{B}_2$ then

 $\mathcal{N}_{E(X|\mathcal{B}_1)}=\mathcal{N}_{E(X|\mathcal{B}_2)}$ and

 $\mathcal{M}_{E(X|\mathcal{B}_1)}=\mathcal{M}_{E(X|\mathcal{B}_2)}$;\\

4). $\mathcal{N}_{E(X|\mathcal{B})}=\mathcal{N}_{E(X|\mathcal{B}_{(2)})}$.\\

5). Let $X^* :=\{x^*:\ x\in X\}$. Then for an arbitrary countable partition $\mathcal{F}$ we have the inclusion
$$\Big(\mathcal{N}_{E(X^*|\mathcal{F}^*)}\Big)^*\subset \Big(\mathcal{N}_{E(X|\mathcal{F})}\Big)^*.$$
$\Box$\\

Let us recall the following notation from Chapter 4. If $\mathcal{B}=(b_n)$ is an interval partition and $f=f^*$ is a monotonic function, then
$$f_{\mathcal{B}}:=\sum_{n\geq 1}f(b_{n-1})\cdot 1_{(b_n,b_{n-1}]}.$$
Clearly $f_{\mathcal{B}}=(f_{\mathcal{B}})^*.$ It is also easy to see that
$$(f^*)_{\mathcal{B}}\leq E(f^*|\mathcal{B})=\Big(E(f^*|\mathcal{B})\Big)^*=\Big(E(f^*|\mathcal{B})\Big)_{\mathcal{B}}\leq \Big(E(f^*|\mathcal{B})^{**}\Big)_{\mathcal{B}}\leq \Big(f^{**}\Big)_{\mathcal{B}}.\eqno(5.1)$$
Because it will not cause any ambiguity, we will denote the left and the right part of $(5.1)$ by $f^*_{\mathcal{B}}$ and $f^{**}_{\mathcal{B}}$, respectively.\\

\textbf{Remark 5.2.} We can rewrite the statement of Theorem 2.2 in the form of the following inequality
$$f^{**}_{\mathcal{B}}\leq \frac{4}{3}\cdot f^{**}_{\mathcal{B}^*}, \eqno (5.2)$$
that is true for any $f\in L^1(I)$ and any interval partition $\mathcal{B}$. Also recall that $\frac{4}{3}$ is the smallest possible constant in $(5.2)$.\\
$\Box$\\

\textbf{Theorem 5.3.} Let $f\in L^1(I)$ and  $\mathcal{B}$ be an interval partition. Then
$$f^{**}_{\mathcal{B}}=\sup\{ E(f^*|\mathcal{S}):\  interval\ partition\ \mathcal{S}\subseteq \mathcal{B}\}. \eqno(5.3)$$

\emph{Proof}. Let $u$ be an arbitrary point in $I$ and let $\mathcal{S}$ be an interval partition such that $\mathcal{S}:=(s_n)\subset (b_n):=\mathcal{B}.$ Let $n$ be such an index that $s_{n+1}<u\leq s_n$. In view of $(5.1)$ we have
$$E(f^*|\mathcal{S})(u)=E(f^*|\mathcal{S})(s_n)\leq f^{**}(s_n)=f_{\mathcal{S}}^{**}(s_n)=f_{\mathcal{B}}^{**}(s_n)=f_{\mathcal{B}}^{**}(u),$$
and therefore the left part here is less or equal to the right part. To prove the opposite inequality let us fix $\varepsilon>0$ and a natural $k$ such that $b_{k+1}<u\leq b_k.$ Let us now chose an $m>k$ such that
$$(b_k-b_m)^{-1}\int_{b_m}^{b_k}f^*d\lambda\geq f^{**}(b_k)-\varepsilon.$$
Therefore, if $\mathcal{S}_{\varepsilon}=(s_n)$ is an interval partition coarser than $\mathcal{B}$ and such that $s_i=b_k,\ s_{i+1}=b_m$ for some $i\geq 0$ then
$$E(f^*|\mathcal{S})(u)\geq f^{**}(b_k)-\varepsilon,$$
Because $\varepsilon$ is arbitrary small we are done.\\
$\Box$\\

\textbf{Theorem} 5.4. Let $\mathcal{F}$ be a countable partition and $f\in L^1(I).$ Then.\\

 1).  $$E(f|\mathcal{F})^*\leq 4/3\cdot f^{**}_{\mathcal{F}^*}.\eqno(5.4)$$
2). The constant $4/3$ in $(5.4)$ is the smallest possible one.\\

\emph{Proof.} From Proposition 0.6.1) we have $E(f|\mathcal{F})^*=E(\bar{f}|\mathcal{B})$ where  $\mathcal{B}\sim \mathcal{F},\ \bar{f}\sim f.$ By using the formulas $(0.2),\ (5.1), (5.2),$ we obtain that
$$E(f|\mathcal{F})^*=E(\bar{f}|\mathcal{B})\leq\Big(E(\bar{f}|\mathcal{B})\Big)^{**}_{\mathcal{B}}\leq\Big(E(f^*|\mathcal{B})\Big)_{\mathcal{B}}^{**}
=f^{**}_{\mathcal{B}}\leq 4/3\cdot f^{**}_{\mathcal{B}^*}=4/3\cdot f^{**}_{\mathcal{F}^*},$$
whence the first part of the theorem follows. To prove the second part assume to the contrary that for any countable partition $\mathcal{F}$ and any $f\in L^1(I)$ we have
$$E(f|\mathcal{F})^*\leq b\cdot f^{**}_{\mathcal{F}^*},\eqno(*)$$
where $0<b<4/3.$ In virtue of 0.6.1) and 0.6.2) inequality $(*)$ is valid for any interval partition $\mathcal{S}$ that is coarser than the interval partition $\mathcal{F}^*.$ Combining it with Theorem 5.3 we see that  $f^{**}_{\mathcal{B}}\leq b\cdot f^{**}_{\mathcal{F}^*}$ for any interval partition  $\mathcal{B}$ that is equimeasurable with $\mathcal{F}^*.$ But that contradicts Remark 5.2.\\
$\Box$\\

From Theorem 5.4 and Lemma 0.14 follows\\

\textbf{Corollary 5.5.} For any countable partition $\mathcal{F}$ and any function $f\in L^1(I)$ we have
$$\mathcal{N}_{E(\mathcal{N}_f|\mathcal{F})}\subseteq\mathcal{N}_{f^{**}_{\mathcal{F}^*}}.$$

\textbf{Theorem 5.6.} For any countable partition $\mathcal{F}$ and any function $f\in L^1(I)$ we have
$$\mathcal{N}_{E(\mathcal{N}_f|\mathcal{F})}=\mathcal{N}_{E(\mathcal{M}_f|\mathcal{F})}.$$
\emph{Proof.} The inclusion $\mathcal{N}_{E(\mathcal{N}_f|\mathcal{F})}\subseteq\mathcal{N}_{E(\mathcal{M}_f|\mathcal{F})}$ follows from the inclusion $\mathcal{N}_f\subseteq\mathcal{M}_f.$ Because a majorant ideal is averaged by any $\sigma$-subalgebra the inverse inclusion follows from the relation
$$E(\mathcal{M}_f|\mathcal{F})\subseteq \mathcal{N}_{E(\mathcal{N}_f|\mathcal{F})},\eqno (5.5)$$
which in its turn follows dirctly from Proposition  0.6.1), Lemma 3.7, and Remark 0.11.1).\\
$\Box$\\

\textbf{Theorem 5.7.} For any symmetric ideal $X$ and for any countable partition $\mathcal{F}$ it is true that $\mathcal{F}$ averages the symmetric ideal $\mathcal{N}_{E(X|\mathcal{F})}.$\\

\emph{Proof.} We will first assume that $X=\mathcal{N}_f$ where $ f\in L^1(I)$. It is enough to prove the inclusion
$$E(\mathcal{N}_{E(\mathcal{N}_f|\mathcal{F})}|\mathcal{F})\subseteq\mathcal{N}_{E(\mathcal{N}_f|\mathcal{F})}.$$ Because $\mathcal{N}_f\subseteq \mathcal{M}_f$ we have $E(\mathcal{N}_f|\mathcal{F})\subseteq\mathcal{M}_f$ and our statement follows from $(5.5)$.\\

Consider now the general case. From $(0.9)$ and already considered special case we obtain
$$E(\mathcal{N}_{\mathcal{N}_{E(X|\mathcal{F})}}|\mathcal{F})=E(\sum_{x\in X}\mathcal{N}_{E(\mathcal{N}_x|\mathcal{F})}|\mathcal{F})=\sum_{x\in X}E(\mathcal{N}_{E(\mathcal{N}_x|\mathcal{F})}|\mathcal{F})\subseteq$$
$$\subseteq \sum_{x\in X}\mathcal{N}_{E(\mathcal{N}_x|\mathcal{F})}=\mathcal{N}_{E(\sum_{x\in X}\mathcal{N}_x|\mathcal{F})}=
\mathcal{N}_{E(X|\mathcal{F})}.$$
$\Box$\\

\textbf{Corollary 5.8}, \cite{Me1}. Let $X$ be a symmetric ideal and $\mathcal{F}$  be a verifying partition (see Definition 6.1). Then $\mathcal{M}_X=\mathcal{N}_{E(X|\mathcal{F})}.$ In particular, $\mathcal{M}_f=\mathcal{N}_{E(\mathcal{N}_f|\mathcal{F})}.$\\
$\Box$\\

\textbf{Theorem 5.9.} For any countable partition $\mathcal{F}$ and any function $f\in L^1(I)$ we have
$$\mathcal{M}_{E(\mathcal{M}_f|\mathcal{F})}=\mathcal{M}_{E(\mathcal{N}_f|\mathcal{F})}=\mathcal{M}_{E(f^*|\mathcal{F}^*)}.$$

\emph{Proof.} The first equality follows directly from Theorem 5.6. To prove the second equalty notice that it follows from Proposition 0.6.1), Remark 5.2, and formulas $(0.13)$ and $(0.2)$ that
$$\mathcal{M}_{E(\mathcal{M}_f|\mathcal{F})}=\sum_{\tilde{f}\sim f,\ \mathcal{\tilde{F}}\sim \mathcal{F}}\mathcal{M}_{E(\mathcal{M}_{\tilde{f}}|\mathcal{\tilde{F}})}=$$
$$=\sum_{\tilde{f}\sim f,\ \textrm{int.part.}\ \mathcal{\tilde{B}}\sim \mathcal{F}}\mathcal{M}_{E(\mathcal{M}_{\tilde{f}}|\mathcal{\tilde{B}})}=\sum_{\textrm{int. part.}\ \mathcal{\tilde{B}}\sim \mathcal{F}}\mathcal{M}_{E(\mathcal{M}_{f^*}|\mathcal{\tilde{B}})}\subseteq\mathcal{M}_{E(f^*|\mathcal{F}^*)}. $$
The converse inclusion is obvious.\\
$\Box$\\

The symmetric ideal generted by the set $E(\mathcal{N}_f|\mathcal{F})$ (unlike the majorant ideal generated by the same set) is not always generated by the function $E(f^*|\mathcal{F}^*)$.  In the next theorem we provide a criterion for it to be true (compare with Theorem 1.13)).\\

\textbf{Theorem 5.10} For any countable partition $\mathcal{F}$ and any function $f\in L^1(I)$ the equality $\mathcal{N}_{E(\mathcal{N}_f|\mathcal{F})}=\mathcal{N}_{E(f^*|\mathcal{F}^*)}$ is valid if and only if the function $E(f^*|\mathcal{F}^*)$ is $\mathcal{F}$-regular (recall the definition after the proof of Proposition 4.8 and Theorem 4.10).\\

\emph{Proof.} Assume the equality in the statement of the theorem. From definition and from Theorem 5.7 follows that the function $E(f^*|\mathcal{F}^*)$ is $\mathcal{F}$-regular.\\

Conversely, let the function $E(f^*|\mathcal{F}^*)$ be $\mathcal{F}$-regular. If for any  $\tilde{f}\sim f$ and $\mathcal{\tilde{F}}\sim\mathcal{F}$ we can prove that $E(\tilde{f}|\mathcal{\tilde{F}})\in \mathcal{N}_{E(f^*|\mathcal{F}^*)}$ then, in virtue of Lemma 0.14 we will obtain the inclusion $\mathcal{N}_{E(\mathcal{N}_f|\mathcal{F})}\subseteq\mathcal{N}_{E(f^*|\mathcal{F}^*)}.$ The converse inclusion is trivial and the theorem will be proved.\\
From Proposition 0.6.1) and formulas $(0.2),\ (5.1)$, and $(_5.2)$ we obtain that for appropriately chosen $\bar{f}\sim f,\ \textrm{int.part.}\ \mathcal{B}\sim \mathcal{F}$ we have
$$E(\tilde{f}|\mathcal{\tilde{F}})^*=E(\bar{f}|\mathcal{B})\leq E(\bar{f}|\mathcal{B})^{**}_{\mathcal{B}}\leq E(f^*|\mathcal{B})^{**}_{\mathcal{B}}=f^{**}_{\mathcal{B}}\leq 4/3\cdot f^{**}_{\mathcal{F}^*}=$$
$$=4/3\cdot E(f^*|\mathcal{F}^*)_{\mathcal{F}^*}^{**}\leq 4/3\cdot Q\cdot\rho_Q E(f^*|\mathcal{F}^*),$$
where $Q$ is the constant from $(4.1)$. Thus, $E(\tilde{f}|\mathcal{\tilde{F}})\in \mathcal{N}_{E(f^*|\mathcal{F}^*)}$.\\
$\Box$\\

\bigskip

We will now state and prove the main results of this Chapter. In the first two of the main theorems we provide criteria for the symmetric ideal $\mathcal{N}_{E(\mathcal{N}_f|\mathcal{F})}$ (or, equivalently, for $\mathcal{N}_{E(\mathcal{M}_f|\mathcal{F})}$) to be a principal symmetric ideal. In the last two theorems we investigate the connections and analogies between properties of regularity and $\mathcal{B}$-regularity of functions from $L^1(I).$\\

\textbf{Theorem 5.11.} For any countable partition $\mathcal{F}$ and any function $f\in L^1(I)$ the following conditions are equivalent.\\

1). There exists a function $g\in L^1(I)$ such that $\mathcal{N}_{E(\mathcal{N}_f|\mathcal{F})}=\mathcal{N}_g;$\\

2). $\mathcal{N}_{E(\mathcal{N}_f|\mathcal{F})}=\mathcal{N}_{f^{**}_{\mathcal{F}^*}};$\\

3). $\mathcal{N}_{E(\mathcal{N}_f|\mathcal{F})}=\mathcal{N}_{E(f^*|\mathcal{F}^*)}.$\\

4). $f^{**}_{{\mathcal{F}}^*}$  is an $\mathcal{F}$-regular function;\\

5). $\mathcal{M}_{f^{**}_{{\mathcal{F}}^*}}=\mathcal{M}_{E(f^*|\mathcal{F}^*)};$\\

6). $f^{**}_{\mathcal{F}^*}\in \mathcal{M}_f.$\\

\emph{Proof.} $1)\Rightarrow 2).$ If $\mathcal{N}_{E(\mathcal{N}_f|\mathcal{F})}=\mathcal{N}_g$ then by Theorem 5.7 $g$ is a $\mathcal{F}$-regular function, and therefore by Theorem 4.10 $g^{**}_{\mathcal{F}^*}\simeq g^*_{\mathcal{F}^*}.$ Also notice that by Theorem 5.9 we have the equality $\mathcal{M}_g=\mathcal{M}_{E(f^*|\mathcal{F}^*)}$ whence $g^{**}\simeq E(f^*|\mathcal{F}^*)^{**},$ whence $g^{**}_{\mathcal{F}^*}\simeq f^{**}_{\mathcal{F}^*}.$ Thus, $g^*\geq g^*_{\mathcal{F}^*}\simeq f^{**}_{\mathcal{F}^*},$ and from it follows the inclusion $\mathcal{N}_{f^{**}_{\mathcal{F}^*}}\subseteq \mathcal{N}_g.$ The converse inclusion follows from Theorem 5.4.\\

The implication $2)\Rightarrow 1)$ us trivial.\\

$2)\Rightarrow 3)$. By definition $f^*_{\mathcal{F}^*}\leq E(f^*|\mathcal{F}^*)$, and therefore
 $\mathcal{N}_{E(\mathcal{N}_f|\mathcal{F})}\subseteq \mathcal{N}_{E(f^*|\mathcal{F}^*)}$. The converse inclusion is trivial.\\

The implication $3)\Rightarrow 4)$ follows from Theorem 5.10 and the fact that the functions $E(f^*|\mathcal{F}^*)$ and $f^{**}_{{\mathcal{F}}^*}$ are either both $\mathcal{F}$-regular, or both - not.\\

$4)\Rightarrow 5).$ The functions $f^{**}_{\mathcal{F}^*}$ и $E(f^*|\mathcal{F}^*)^{**}_{\mathcal{F}^*}$ are identical, and therefore by Theorem 4.10 ${(f^{**}_{\mathcal{F}^*})}^{**}_{\mathcal{F}^*}\simeq E(f^*|\mathcal{F}^*)^{**}_{\mathcal{F}^*}$; together with the fact that the functions $E(f^*|\mathcal{F}^*)$ and $f^{**}_{\mathcal{F}^*}$ have finite range relatively to $\mathcal{F}^*$ it implies that $(f^{**}_{\mathcal{F}^*})^{**}\simeq E(f^*|\mathcal{F}^*)^{**},$ i.e.
$\mathcal{M}_{f^{**}_{{\mathcal{F}}^*}}=\mathcal{M}_{E(f^*|\mathcal{F}^*)}.$\\

$5)\Rightarrow 6).$ Modulo the inclusion $\mathcal{M}_{E(f^*|\mathcal{F}^*)}\subseteq \mathcal{M}_f$ this implication is trivial.\\

$6)\Rightarrow 2).$ The inclusion $f^{**}_{{\mathcal{F}}^*}\in \mathcal{M}_f$ can be written in the form $f^{**}_{{\mathcal{F}}^*}\in E(\mathcal{M}_f|\mathcal{F}^*).$ It follows from Theorem 5.6 and Lemma 5.1.2) that
$$\mathcal{N}_{f^{**}_{{\mathcal{F}}^*}}\subseteq \mathcal{N}_{E(\mathcal{M}_f|\mathcal{F}^*)}=\mathcal{N}_{E(\mathcal{N}_f|\mathcal{F}^*)}=\mathcal{N}_{E(\mathcal{N}_f|\mathcal{F})}.$$
The converse inclusion follows from Theorem 5.4.\\
$\Box$\\

\textbf{Theorem 5.12} Let $f=f^*\in L^1,\ \mathcal{B}=\mathcal{B}^*=(b_n)$ be some monotonic interval partition. The following conditions are equivalent.\\

1. The function $f$ is  $\mathcal{B}$-regular;\\

2. The function $E(f|\mathcal{B})$ is $ \mathcal{B}$-regular.\\

3. The function $f^{**}_{\mathcal{B}}$ is $ \mathcal{B}$-regular.\\

4. $\mathcal{N}_{f_{\mathcal{B}}}=\mathcal{N}_{f^{**}_{\mathcal{B}}}$;\\

5. $\mathcal{N}_{E(\mathcal{N}_f|\mathcal{B})}=\mathcal{N}_{E(f|\mathcal{B})}.$\\

\textbf{Remark} It is worth noticing that if we put $f=f^*,\ \mathcal{F}=\mathcal{B}=\mathcal{B}^*$ then we can see that condition 3 of theorem 5.11 is identical to condition 5 of Theorem 5.12, and therefore all the conditions in these theorems are equivalent.\\

\emph{Proof} of Theorem 5.12.  $1.\Rightarrow 2.$ For a $\mathcal{B}$-regular function $f$ in virtue of $(0.2)$ and $(5.1)$ we have\\
 $$E(f|\mathcal{B})^{**}_{\mathcal{B}}=f^{**}_{\mathcal{B}}\simeq f_{\mathcal{B}}\leq E(f|\mathcal{B})=E(f|\mathcal{B})_{\mathcal{B}}\leq E(f|\mathcal{B})^{**}_{\mathcal{B}},$$
 and therefore the function $E(f|\mathcal{B})$ is also $\mathcal{B}$-regular.\\

 $2.\Rightarrow 3.$ If the function $E(f|\mathcal{B})$ is $\mathcal{B}$-regular then by Theorem 5.10 $\mathcal{N}_{E(\mathcal{N}_f|\mathcal{B})}=\mathcal{N}_{E(f|\mathcal{B})}.$ Now we can apply  parts 1 and 4 of Theorem 5.11 and conclude that the function $f^{**}_{\mathcal{B}}$ is $\mathcal{B}$-regular.\\

 $3.\Rightarrow 4.$ The implication follows immediately from Theorem 4.10. \\

 $4.\Rightarrow 1.$ This implication follows directly from the definition of $\mathcal{B}$-regularity (see also Theorem 4.10).\\

 $1.\Rightarrow 5.$ By Theorem 4.10 we have $f^{**}_{\mathcal{B}^*}\simeq f_{\mathcal{B}}$. Next we apply parts 2 - 4 of Theorem 5.11 to obtain that $\mathcal{N}_{E(\mathcal{N}_f|\mathcal{B})}\subseteq \mathcal{N}_{E(f|\mathcal{B})}$. The converse inclusion is trivial.\\

 $5.\Rightarrow 2.$ This implication follows directly from Theorem 5.11.\\

$\Box$\\

One of our main results, Theorem 5.15 below, shows that the functions $g:=f^{**}_{\mathcal{F}^*}$  and $g^{**}_{\mathcal{B}^*}$ are either both $\mathcal{B}$-regular, or both not. It shows again the analogy between the properties of regularity and $\mathcal{B}$-regularity (because the functions $f$ and $f^{**}$ by Corollary  1.5 are both regular or both not). The following Theorem 5.13 underlines some reasons for this analogy..\\

To state this theorem we assume that $\mathcal{T}=(t_n)$ is an arbitrary interval partition and define the following two functions of $k,\ k\geq 1$:\\

$\varpi_{\mathcal{T}}(k):=\sup_{n\geq 0}\frac{t_{n+k}}{t_n};$\\

$q_{\mathcal{T}}(k)=$ \emph{the number of points of the sequence} $\{t_n\}_{n=0}^\infty$ \emph{on the diadic interval}\ $D_k:=(2^{-k},2^{-k+1}].$ \\

\textbf{Theorem 5.13}.  The following conditions are equivalent.\\

i). There is a natural number $m$ such that $\varpi_{\mathcal{T}}(m)<1;$\\

ii). The sequence $\{q_{\mathcal{T}}(k)\}_{k=1}^\infty$ is bounded.\\

\emph{Proof}. The implication $ii)\ \Rightarrow i)$ is trivial. Indeed, let $q_{\mathcal{T}}(k)\leq Q,\ k\geq 1$ for some natural $Q$. Take $m=3Q$ and condition $i)$ is clearly satisfied.\\

Let us prove that $i)\ \Rightarrow ii)$. In every interval $D_k:=(2^{-k},2^{-k+1}],\ k\geq 1$ we will take its midpoint to get the intervals $D_{k|0},D_{k|1},\ k\geq 1.$\\

At this point it is convenient to remark that if the sequence $\{q_{\mathcal{T}}(k)\}_{k=1}^\infty$ is unbounded then, either for $\omega=0$ or for $\omega=1$, the sequence $\{m_k\}$ where $m_k$ is the number of points  of $\mathcal{T}$ on the interval $D_{k|\omega},\ k\geq 1, $ is also unbounded.\\

Next we will half the intervals $D_{k|0},D_{k|1},\ k\geq 1$ and denote the intervals we obtain by $D_{k|0,0},D_{k|0,1},D_{k|1,0},D_{k|1,1},\ k\geq 1$ every time choosing between these intervals one that contains an infinite subset of points from $\mathcal{T}$. By continuing this procedure on step number $i$ for every $k\geq 1$ we obtain the intervals $D_{k|\bar{\omega_i}}$ where $\bar{\omega_i}=(\omega_1,...,\omega_i)$ are the first $i$ digits of the binary representation of $\omega\in I$.\\

By generalizing the remark we made  above in the third paragraph of the current proof we can claim that if the sequence $\{q_{\mathcal{T}}(k)\}_{k=1}^\infty$ is unbounded then we can find a sampling $D_{n_k|\bar{\omega}_{i(n_k)}},\ k\geq 1$ such that the sequence $\{r_k\}$ where $r_k$ is the number of points of $\mathcal{T}$ on the interval $D_{n_k|\omega_{i(n_k)}},\ k\geq 1,$ is also unbounded.\\

Let us denote by $\overleftarrow{t}_{n_k|\omega_{i(n_k)}}$ and $\overrightarrow{t}_{n_k|\omega_{i(n_k)}}$ the extreme to the left and the extreme to the right point of  $\mathcal{T}$ on the interval $D_{n_k|\omega_{i(n_k)}},\ k\geq 1,$ respectively. By $r_{n_k|\omega_{i(n_k)}}$ we denote the number of $\mathcal{T}$-points between these two extreme points.  Our construction guarantees that $\limsup_{k\rightarrow\infty}r_{n_k|\omega_{i(n_k)}}=\infty$ and therefore for any natural $m$ we can find  $r_{n_k|\omega_{i(n_k)}}>m$. On the other hand $\overleftarrow{t}_{n_k|\omega_{i(n_k)}}$ and $\overrightarrow{t}_{n_k|\omega_{i(n_k)}}$ are in the same diadic interval  $D_{n_k|\omega_{i(n_k)}},\ k\geq 1$,

where $\lambda(D_{n_k|\omega_{i(n_k)}})_{k\rightarrow\infty}\rightarrow 0.$ Therefore  assuming that $m<r_{n_k|\omega_{i(n_k)}}$ for the points $t_{n+m},t_n\in D_{n_k|\omega_{i(n_k)}}$ we have contrary to condition $i)$ that
$$\varpi_{\mathcal{T}}(m)=\sup_{n\geq 0}\frac{t_{n+m}}{t_n}\geq
\limsup_{k\rightarrow\infty}\frac{\overleftarrow{t}_{n_k|\omega_{i(n_k)}}}{\overrightarrow{t}_{n_k|\omega_{i(n_k)}}}=1.$$ The obtained contradiction proves the implication $i)\ \Rightarrow ii)$ .\\
$\Box$\\

Let us recall that (see Theorem 1.8) with the function  $f=f^*\in L^1(I)$ and the concave function $\psi(t):=\int_0^tfd\lambda,\ t\in I$ we connect the sequence of natural numbers $q_{\psi(\mathcal{D})}(k)$:

$$q_{\psi(\mathcal{D})}(k)=\ \emph{\textrm{the number\ of points\ of the sequence}}\ \{\psi(2^{-n})\}_{n=0}^\infty$$ $$\emph{\textrm{on\ the  interval}}\ D_k,\ k\geq 1.$$ \\

It follows from Lemma 1.10 (via Corollary 1.11) that $f$ is regular if and only if when the sequence of natural numbers  $q_{\psi(\mathcal{D})}(k)$ is bounded. We will now obtain an analogous statement for $\mathcal{B}$-regularity. In vitrue of 0.11.4) instead of $\mathcal{B}$ we can consider the dyadic interval partition ${\mathcal{B}}^*_{(2)}:=(2^{-m_n}).$\\

We define the following sequence
$$\{q_{\psi(\mathcal{B}^*_{(2)})}(k)\}_{k=1}^\infty = \ \textrm{the number\ of points\ of the sequence}\ \{\psi(2^{-m_n})\}_{n=0}^\infty$$  $$\textrm{on\ the interval}\ D_k,\ k\geq 1.$$

\bigskip

\textbf{ Theorem 5.14. }  The function $f^{**}_{\mathcal{B}^*}$ is $\mathcal{B}^*$-regular if and only if  the sequence of natural numbers $\{q_{\psi(\mathcal{B}^*_{(2)})}(k)\}_{k=1}^\infty$ is bounded.\\

\emph{Proof.} Because $f^{**}_{\mathcal{B}_{(2)}^*}=\Big(E(f^*|\mathcal{B}_{(2)}^*)\Big)_{\mathcal{B}_{(2)}}^{**}$ we can take as function $f$ the function $E(f^*|\mathcal{B}_{(2)}^*)$. It follows from Theorems 5.11, 5.13, and 5.1, as well as from Remark 4.10.4 and from the fact that relatively to  $\mathcal{B}^*_{(2)}$ the functions $f^{**}_{\mathcal{B}_{(2)}^*}$ и $E(f^*|\mathcal{B}_{(2)}^*)$ have finite range, that it is enough to prove that the following two conditions are equivalent.\\
1). There exists a $C>1$ such that
$$\int_0^{s_n}f^{**}_{\mathcal{S}}d\lambda\leq C\cdot \psi(s_n),\ n\geq 0;$$
2). There are a natural $p$ and a real $\gamma,\ 0<\gamma<1$ such that
$$\psi(s_{n+p})\leq \gamma\cdot\psi(s_{n}),\ n\geq 0,\eqno (5.6)$$
where $\mathcal{S}:=(s_n)$ denotes the interval partition $\mathcal{B}_{(2)}.$\\

Let us prove that  2) implies 1). Inequality $(5.6)$ can be written in the form
$$\psi(s_{n+kp+i})\leq \gamma^k\cdot\psi(s_{n+i}),\ n\geq 0,\ k\geq 0,\ i=0,1,...,p-1.$$
Whence, for any $n\geq 0$ we obtain
$$\int_0^{s_n}f^{**}_{\mathcal{S}}d\lambda=\int_0^{s_n}\sum_{k=n}^\infty f^{**}(s_k)1_{(s_{k+1},s_k]}d\lambda
=\sum_{k=0}^\infty\sum_{i=0}^{p-1}f^{**}(s_{n+kp+i})(s_{n+kp+i}-s_{n+kp+i+1})\leq $$
$$\leq \sum_{k=0}^\infty\gamma^k\sum_{i=0}^{p-1}\psi(s_{n+i})\leq \frac{p}{1-\gamma}\cdot\psi(s_{n}).$$
It remains to take $C=\frac{p}{1-\gamma}$.\\

We will prove now that 1) implies 2). It follows from  1) that for any natural $n$
$$\frac{\int_0^{s_{n+1}}f^{**}_{\mathcal{S}}d\lambda}{\int_0^{s_n}f^{**}_{\mathcal{S}}d\lambda}=
1-\frac{\psi(s_n)}{\int_0^{s_n}f^{**}_{\mathcal{S}}d\lambda}
(1-\frac{s_{n+1}}{s_{n}})\leq 1-\frac{1}{2C}:=\gamma\in (0,1).$$
Therefore, for any natural $p$ we have the inequalities
$$\psi(s_{n+p})\leq \int_0^{s_{n+p}}f^{**}_{\mathcal{S}}d\lambda\leq \gamma\cdot\int_0^{s_{n+p-1}}f^{**}_{\mathcal{S}}d\lambda\leq...$$ $$\leq \gamma^p\cdot\int_0^{s_n}f^{**}_{\mathcal{S}}d\lambda\leq C\cdot\gamma^p\cdot\psi(s_{n}).$$
It remains to chose a $p$ such that $C\cdot\gamma^p<1.$\\
$\Box$\\

In addition to being useful for verifying conditions of Theorem 5.11 Theorem  5.14 is of independent interest because it helps to construct principal symmetric ideals that are averaged or not averaged by any given diadic interval partition. We will use it later in proofs of Theorems 6.8 and 6.4.\\

\textbf{Theorem 5.15.} Let $\mathcal{F}$ be an arbitrary countable partition and let $f$ be an arbitrary function in $L^1(I)$. The functions  $g:=f^{**}_{\mathcal{F}^*}$ and $g^{**}_{\mathcal{F}^*}$ are either both $\mathcal{F}$-regular or both not.\\

\emph{Proof.} The implication ($g$ is $\mathcal{F}$-regular) $\Rightarrow$ ($g^{**}_{\mathcal{F}^*}$ is $\mathcal{F}$-regular)
follows directly from Theorem 4.10.2). Let us prove the converse implication.\\

Without loss of generality we can assume that $\mathcal{F}^*=\mathcal{F}_{(2)}^*:=(2^{-m_n})$ (see Remark 0.11.3). Let $\psi_1(t):=\int_0^t g^{**}_{\mathcal{F}^*}d\lambda;\ \psi_2(t):=\int_0^t g_{\mathcal{F}^*}d\lambda=\int_0^t f^{**}_{\mathcal{F}^*}d\lambda,\ t\in I.$ By Theorem 5.13 it is enough to show that if the sequence $\{q_{\psi_1(\mathcal{F}^*_{(2)})}(k)\}_{n=0}^\infty$ is bounded then the sequence $\{q_{\psi_2(\mathcal{F}^*_{(2)})}(k)\}_{n=0}^\infty$ is also bounded.\\

Assume to the contrary that
$$\limsup_{k\rightarrow\infty}q_{\psi_2(\mathcal{F}^*_{(2)})}(k)=\infty,\  \textrm{while}\  \limsup_{k\rightarrow\infty}q_{\psi_1(\mathcal{F}^*_{(2)})}(k)\leq K\eqno(5.7)$$
for some natural $K$. We can assume that $K$ is so large that
$$2K^2+1\leq 2^{K^2-1}.\eqno(5.8)$$
For any $n\geq 0$ we have
$$\int_0^{2^{-m_n}}g^{**}_{\mathcal{F}^*}d\lambda=\int_0^{2^{-m_n}}\sum_{i\geq n}g^{**}(2^{-m_i})\cdot 1_{(2^{-m_{i+1}},2^{-m_i}]}=$$
$$=\sum_{i\geq n}g^{**}(2^{-m_i})(2^{-m_i}-2^{-m_{i+1}})=\sum_{i\geq n}(1-2^{m_i-m_{i+1}})\int_0^{2^{-m_i}}gd\lambda.$$
Because  $m_{i+1}-m_i\geq 1$ for any $i\geq 0$ we have
$$2^{-1}\leq 1-2^{m_i-m_{i+1}}\leq 1,\ i\geq 0,$$
and thus
$$2^{-1}\cdot\sum_{i\geq n} \int_0^{2^{-m_i}}gd\lambda\leq\int_0^{2^{-m_n}}g^{**}_{\mathcal{F}^*}d\lambda\leq\sum_{i\geq n}\int_0^{2^{-m_i}}gd\lambda,\ n\geq 0,\eqno(5.9)$$
In virtue of $(5.7)$ there is an interval $D_j=(2^{-j},2^{-j+1}],\ j\geq 1$ that contains at least $K^2+1$ points of the sequence $\{\int_0^{2^{-m_n}}gd\lambda\}_{n=0}^\infty.$ Let $n$ be a natural number such that $$2^{-j+1}\geq \int_0^{2^{-m_n}}gd\lambda\geq \int_0^{2^{-m_{n+1}}}gd\lambda\geq ...\geq \int_0^{2^{-m_{n+K^2}}}gd\lambda>2^{-j}.\eqno(5.10)$$
It follows from Lemma 0.1, $(5.9)$ and $(5.10)$ that
$$\frac{1}{2^{K^2}}\geq \frac{\int_0^{2^{-m_{n+K^2}}}g^{**}_{\mathcal{F}^*}d\lambda}{\int_0^{2^{-m_n}}g^{**}_{\mathcal{F}^*}d\lambda}\geq \frac{1}{2}\cdot\frac{\sum_{i\geq n+K^2}\int_0^{2^{-m_i}}gd\lambda}{\sum_{i\geq n}\int_0^{2^{-m_i}}gd\lambda}=$$
$$=\frac{1}{2}\cdot\Big(1-\frac{\sum_{i=n}^{n+K^2-1}\int_0^{2^{-m_i}}gd\lambda}{\sum_{i\geq n}\int_0^{2^{-m_i}}gd\lambda}\Big)=\frac{1}{2}\cdot[1-\Big(1+\frac{\sum_{i\geq n+K^2}\int_0^{2^{-m_i}}gd\lambda}{\sum_{i= n}^{n+K^2-1}\int_0^{2^{-m_i}}gd\lambda}\Big)^{-1}].\eqno(5.11)$$
But, according to $(5.9)$ and $(5.10)$ we have
$$\sum_{i= n}^{n+K^2-1}\int_0^{2^{-m_i}}gd\lambda\leq K^2\cdot 2^{-j+1};\ \sum_{i\geq n+K^2}\int_0^{2^{-m_i}}gd\lambda\geq \int_0^{2^{-m_{n+K^2}}}gd\lambda>2^{-j}.$$
Therefore it follows from $(5.11)$ that
$$\frac{1}{2^{K^2}}>\frac{1}{2}[1-\Big(1+\frac{2^{-j}}{K^2\cdot 2^{-j+1}}\Big)^{-1}]=\frac{1}{2}(1-\frac{2K^2}{2K^2+1})=\frac{1}{2}\cdot\frac{1}{2K^2+1}$$
Therefore $2K^2+1>2^{K^2-1}$ in contradiction with $(5.8)$. Theorem 5.15 is proved.\\
$\Box$\\

\textbf{Theorem 5.16.} Let $f=f^*\in L^1(I),\ \mathcal{B}=\mathcal{B}^*$. The equalities
$$\mathcal{N}_{E(\mathcal{N}_{E(\mathcal{N}_f|\mathcal{B})}|\mathcal{B})}=
 \mathcal{N}_{E(\mathcal{N}_f|\mathcal{B})}=\mathcal{N}_{E(f|\mathcal{B})}\eqno(5.12)$$
 take place if and only if the function $f$ is $\mathcal{B}$-regular.\\

\emph{Proof.} Let $f$ be $\mathcal{B}$-regularа. Then the second of equalities  $(5.12)$ is true by Theorem 5.12.5). We need to prove the first of equalities $(5.12)$. By Theorem 5.7 for any function  $f\in L^1(I)$ the interval partition $\mathcal{B}$ averages the symmetric ideal $\mathcal{N}_{E(\mathcal{N}_f|\mathcal{B})}$ and therefore $$E(\mathcal{N}_{E(\mathcal{N}_f|\mathcal{B})}|\mathcal{B})\subseteq\mathcal{N}_{E(\mathcal{N}_f|\mathcal{B})}.$$
Therefore by plugging in the left part of $(5.12)$ $\mathcal{N}_{E(\mathcal{N}_f|\mathcal{B})}$ instead of $E(\mathcal{N}_{E(\mathcal{N}_f|\mathcal{B})}|\mathcal{B})$ we obtain
$$\mathcal{N}_{E(\mathcal{N}_{E(\mathcal{N}_f|\mathcal{B})}|\mathcal{B})} \subseteq \mathcal{N}_{\mathcal{N}_{E(\mathcal{N}_f|\mathcal{B})}}
=\mathcal{N}_{E(\mathcal{N}_f|\mathcal{B})}.$$
 On the other hand, obviously $E(E(\mathcal{N}_f|\mathcal{B})|\mathcal{B})\subseteq E(\mathcal{N}_{E(\mathcal{N}_f|\mathcal{B})}|\mathcal{B})$ and therefore
 $$E(f|\mathcal{B})\in E(\mathcal{N}_f|\mathcal{B})=E(E(\mathcal{N}_f|\mathcal{B})|\mathcal{B})\subseteq E(\mathcal{N}_{E(\mathcal{N}_f|\mathcal{B})}|\mathcal{B})\subseteq
 \mathcal{N}_{E(\mathcal{N}_{E(\mathcal{N}_f|\mathcal{B})}|\mathcal{B})}.$$
Thus, the equality  $\mathcal{N}_{E(\mathcal{N}_{E(\mathcal{N}_f|\mathcal{B})}|\mathcal{B})}=
 \mathcal{N}_{E(\mathcal{N}_f|\mathcal{B})}$ is proved.\\

The converse: if we have the equalities  $(5.12)$ then the second of them by Theorem 5.12.5) implies the $\mathcal{B}$-regularity of $f.$\\
$\Box$\\

\newpage

\bigskip

\centerline{\textbf{Chapter 6. Verifying and universal $\sigma$-subalgebras.}}
\bigskip

\centerline{\ \textbf {The results of this chapter are based on papers~\cite{Me1},\ ~\cite{Me10}, ~\cite{Me17}.}}

\bigskip

\emph{In this Chapter we begin the study of \textbf{verifying} $\sigma$-subalgebras of $\Lambda$, i,e. such subalgebras that a symmetric ideal is majorant if it is invariant under the action of \textbf{just one operator} - the averaging operator corresponding to a verifying subalgebra. We also introduce the class of \textbf{universal} $\sigma$-sublageras; for a subalgebra from this class the corresponding averaging operator leaves invariant \textbf{any symmetric ideal}.}

\textit{We will continue the study of verifying and universal $\sigma$-subalgebras also in Chapters 7 and 8}.

\bigskip

Theorem 3.2 states that if some ideal $X,\ X\subseteq L^1(I),$ is averaged by any $\sigma$-subalgebra of $\lambda$ then $X$ is a majorant symmetric ideal. But, if a single $\sigma$-subalgebra averages some symmetric ideal, it is in general not enough to guarantee that it is a majorant ideal. This simple observation provides a rationale for the following definition.\\

\textbf{Definition 6.1.} A $\sigma$-subalgebra $\mathcal{A}$  is called a  \emph{verifying} $\sigma$-subalgebra if from the fact that $\mathcal{A}$ averages a symmetric ideal$X$ follows that $X$ is a majorant ideal.\\

It is obvious that if a  $\sigma$-subalgebra is equimeasurable with a verifying one then it itself is verifying. It follows from Remark 0.11.3),4) that a countable partition is verifying if and only if its diadic projection is verifying, and also that if a countable partition is finer than a verifying one (or equivalent to a finer countable partition) then it is itself a verifying one.\\

It follows from Remark 0.10 that the notions of verifying and universal $\sigma$-subalgebras can be introduced for any continuous probability space $(\Omega,\Sigma,\mu)$ and that these two properties are invariant under the relation of $e$-equivalence.\\

\textbf{Theorem 6.1.} Let $(\Omega,\Sigma,\mu)$ be a continuous probability space. Let $\mathcal{A}$ be a continuous $\sigma$-subalgebra of $\Sigma$ such that $\mathcal{A}$ is complemented in $\Sigma$ and the complementary $\sigma$-subalgebra $\mathcal{A}^\perp$ is also continuous. Then $\mathcal{A}$ is a verifying $\sigma$-subalgebra.

\emph{Proof. }  To avoid cumbersome notations we will assume that $(\Omega,\Sigma,\mu)=(I,\Lambda,\lambda)\times (I,\Lambda,\lambda)$ is the unit square with the Lebesgue measure and that the $\sigma$-subalgebra $\mathcal{A}$ is generated by all Lebesgue measurable vertical strips, i.e. by the sets of the form $A\times I$ where $A\in \Lambda$. The arguments below can be repeated verbatim in the general case.\\

Notice that the operator of averaging by $\mathcal{A}$ acts as averaging of a function $f(t,s),\ f\in L^1(I^2)$ by the second variable:
$$E(f|\mathcal{A})(u,v)=\int_0^1f(u,v)dv.\eqno(6.1)$$
Assume that $\mathcal{A}$ averages a symmetric ideal $X\subseteq L^1(I^2)$ that is not majorant:
$$E(X|\mathcal{A})\subseteq X.\eqno(6.2)$$
The symmetric ideal $\tilde{X}\subseteq L^1(I)$ is $\stackrel{e}{\sim}$equivalent to $X$ and therefore also is  not a majorant ideal. By Theorem 3.2 we can find an interval partition of $I$, $\mathcal{\tilde{P}}:=\sigma(\tilde{P}_n,\ n\geq 1)$ and an $\tilde{x},\ \tilde{x}\in\tilde{X}$ such that
$$E(\tilde{x}|\mathcal{\tilde{P}})\notin \tilde{X}.\eqno(6.3)$$
Without loss of generality we can assume that the function $\tilde{x}$ is an elementary one and $\tilde{x}_n:=\tilde{x}\cdot 1_{\tilde{P}_n}=\sum_{i\geq 1}r_i^n\cdot 1_{\tilde{A}_i^n}$ where the sets $\tilde{A}_i^n\subseteq I$, $i\neq j$, are pairwise disjoint, $\lambda(\tilde{A}_i^n := \alpha_i^n )>0$, and $\bigcup_{i\geq 1}\tilde{A}_i^n=\tilde{P}_n,\ n\geq 1.$ Let $\alpha_n:=\lambda(\tilde{P}_n),\ n\geq 1$ and let $\sigma(P_i^n,\ i\geq 1)$ be the countable partition of $I$ into pairwise disjoint intervals $\tilde{P}_i^n$ of the length $\frac{\alpha^n_i}{\alpha_n},\ i,n\geq 1$. Such a partition exists by Proposition 0.4.\\

Let $P_n:=\tilde{P}_n\times I;\ P_i^n:=\tilde{P}_n\times \tilde{P}_i^n,\ i,n\geq 1.$ It is easy to see that for any $n\geq 1$ the functions $x_n:=\Sigma_{i\geq 1}r_i^n1_{P_i^n}$ and $\tilde{x}_n$ are equimeasurable, whence the functions $\tilde{x}$ and $x:=\Sigma_{n\geq 1}x_n$ are also equimeasurable, and therefore $x\in X$.
Denote by $\mathcal{P}$ the partition $\sigma(P_n,\ n\geq 1)$ of the unit square. Applying formula $(6.1)$ we come to a contradiction with $(6.2)$, namely, for $x\in X$ we have
$$E(x|\mathcal{A})=E(x|\mathcal{P})\sim E(\tilde{x}|\tilde{\mathcal{P}})\notin \tilde{X}.$$
$\Box$\\

Let us now discuss countable verifying partitions.\\

\textbf{Lemma 6.3.} For any $p,\ 0<p<1$ the interval partition $\mathcal{B}^{(p)}:=(p^{n})$ is verifuying.\\

\emph{Proof.} We need to prove that if $\mathcal{B}^{(p)}$ averages a symmetric ideal  $X,\ X\subseteq L^1(I,\Lambda,\lambda)$, then $X$ is a majorant ideal.\\

Let us fix a $\varepsilon>0$ and assume that $y\prec z\in X$. Without loss of generality we can assume that $y$ is an elementary function. By formula $(0.2)$ we have $y^*\simeq E(\rho_py^*|\mathcal{B}^{(p)})$ and therefore it is enough to prove that  $E(\rho_py^*|\mathcal{B}^{(p)})\in X.$ But $\rho_py^*\prec\rho_pz^*\in X$, whence $\int_0^t\rho_py^*d\lambda<\int_0^tx^*,\ t\in I$ where $x^*:=\rho_pz^*+\varepsilon\cdot \textbf{1}\in X.$ It follows from Lemma 3.7 that there is a function $\tilde{x}$ equimeasurable with $x^*$ (whence $\tilde{x} \in X$) and such that  $E(\rho_py^*|\mathcal{B}^{(p)})\leq E(\tilde{x}|\mathcal{B}^{(p)})\in X.$ \\
$\Box$\\

\textbf{Theorem 6.4.} Let $\mathcal{T}=(t_n)$ be an interval partition belonging to the stochastic vector $\vec{a}=(a_n)$. For $\mathcal{T}$ to be verifying it is sufficient and, when  $\vec{a}$ is monotonic, also necessary that
$$\sup_{n\geq 1}\frac{a_n}{\sum_{k>n}a_k}\Big(=\sup_{n\geq 1}\frac{t_{n-1}-t_n}{t_n}\Big)<\infty,$$
$$\textrm{ or equivalently,}\  \sup_{n\geq 1}\frac{t_{n-1}}{t_n}<\infty. \eqno (6.4)$$

To prove Theorem 6.4 we will need the following lemma.\\

\textbf{Lemma 6.5.} Condition $(6.4)$ is satisfied if and only if there is a real number $p,\ 0<p<1$ such that every interval $B_n^{(p)}=(p^{n},p^{n-1}],\ n\geq 1,$ of the interval partition $\mathcal{B}^p$ contains at least one point of the sequence $\{t_n\}.$\\

\emph{Proof.} Necessity. Assume that for any $p$, $\ p=2^{-\nu^2},\ \nu=1,2,...,$ there us an interval $B_m^{(2^{\nu^2})}$ that does not contain a single point of the sequence $\{t_n\}.$  Let us denote by $t_{n-1}^\nu$ and $t_n^\nu$ points of the sequence $\{t_n\}$ that are the closest to this interval from the right and, respectively, from the left. Then $\limsup_{n\rightarrow \infty}\frac{t_{n-1}}{t_n} \geq \limsup_{\nu\rightarrow \infty}\frac {2^{-(\nu-1)^2}}{2^{-\nu^2}}=\limsup_{\nu\rightarrow \infty}2^{2\nu-1}=\infty$, in contradiction with $(6.4)$.\\

Sufficience. Let  $p,\ 0<p<1$ be such a real number that every interval $B_n^{(p)}=(p^{n},p^{n-1}],\ n\geq 1,$ of the interval partition $\mathcal{B}^p$ contains at least one point of the sequence $\{t_n\}.$ Let $t_k$ be an arbitrary point of the sequence $\{t_n\}$.
We have to consider two cases.

 1) $t_k,\ t_{k-1}\in B_n^{(p)}.$ Then $\frac{t_{k-1}}{t_k}\leq \frac{p^{n-1}}{p^{n}}=p<\infty.$
2). $t_k \in B_{n+1}^{(p)},\ t_{k-1}\in B_n^{(p)}.$ Then $\frac{t_{k-1}}{t_k}\leq \frac{p^{n-1}}{p^{n+1}}=p^{-2}.$

Thus, in both cases condition $(6.4)$ is satisfied. \\
$\Box$\\

\textit{Proof} of Theorem 6.4. Sufficiency. Let $C$ be the left part of $(6.4)$ and let $p=(C+1)^{-1}$. From the sequence $\{t_n\}$ we will remove all the points that are contained in the intervals $B_{2n-1}^{(p)},$ while in the intervals $B_{2n}^{(p)}$ we will leave only one point of this sequence in each such interval, $n\geq 1$. The remaining points we will organize in the sequence $\{t'_n\}$ in decreasing order on $I$ assuming that $t'_0=1.$\\

By Proposition 0.6 it is enough to prove that the interval partition  $\mathcal{T}':=(t'_n)$ is verifying.  By Lemma 6.3 every partition that belongs to the stochastic vector $\vec{p}_2=(p^{2n})$ (multiplied by $K:=p^{-2}-1$) is verifying. Therefore, if $X$ is a symmetric but not a majorant ideal then by Remark 0.11.4) there is $x \in X$ such that
$$x=x\cdot 1_{\bigcup_{n\geq 1}B_{2n-1}^{(p^2)}}\ \textrm{and}\ \sum_{n\geq 1} [\lambda(B_{2n-1}^{(p^2)})]^{-1}\int_{B_{2n-1}^{(p^2)}}x\cdot d\lambda 1_{B_{2n-1}^{(p^2)}}\notin X.$$
From our choice of the sequence $\{t'_n\}$  and the element $x$ we have
$$E(x|\mathcal{T}')=K\sum_{n\geq 0}\frac{1}{t'_n-t'_{n+1}}\int^{t'_n}_{t'_{n+1}}x\cdot d\lambda\cdot1_{(t'_{n+1},t'_n]}\geq$$
$$\geq K\rho_{p^{-2}}\Big(\sum_{n\geq 1} [\lambda(B_{2n-1}^{(p^2)})]^{-1}\int_{B_{2n-1}^{(p^2)}}x\cdot d\lambda 1_{B_{2n-1}^{(p^2)}}\Big)\notin X,$$
and the sufficiency is proved.\\

\emph{Necessity.} Assume that interval partition $\mathcal{T}=(t_n)$ belongs to the monotonic stochastic vector $\vec{a}=(a_n)$ and does not satisfy condition $(6.4)$. Our goal is to construct a function $f=f^*\in L^1(I)$  such that the symmetric ideal $\mathcal{N}_f$ would be not majorant but at the same time would be averaged by $\mathcal{T}$.\\

By Remark 0.11.3) without loss of generality we can instead of the interval partition $\mathcal{T}$ consider its diadic projection, $\mathcal{T}_{(2)}:=(2^{-m_n}),$ where $0=m_0<m_1<...<m_n<...$ and by assumption $\limsup_{n\rightarrow\infty}(m_{n+1}-m_n)=\infty.$\\

Let us consider the sequence of natural numbers $\{q_n\},\ q_n=m_{n+1}-m_n,\ n \geq 1,$ and let $q_{m_n}^\prime=1,\ q_{m_n+1}^\prime=q_{m_n+1}+q_{m_n}-1,\ q_k^\prime=q_k\ \textrm{if}\ k\neq m_n,\ k\neq m_{n}+1,\ n\geq 1.$\\

The subsequence $\{q_{m_n}^\prime\}$ of the unbounded sequence of naturals $\{q_n^\prime\}$ is by construction bounded: $q_{m_n}^\prime=1.$  By using the sequence $\{q_n^\prime\}$ we will construct a non-increasing function $f$ by the rules $\mathcal{U}_1)\ -\ \mathcal{U}_4)$ of Lemma 1.9. Lemma  1.10 shows that $f$ is not regular and therefore the symmetric ideal $\mathcal{N}_f$ is not majorant. At the same time, by Theorem 5.14 $\mathcal{T}_{(2)}$ averages the symmetric ideal $\mathcal{N}_f$.\\
$\Box$\\

\textbf{Remark 6.6.} 1. In the necessity part of Theorem 6.4 we cannot drop the condition that the stochastic vector  $\vec{a}$ is monotonic. Indeed, for any stochastic vector there is a permutation of its coordinates such that the obtained stochastic vector would not satisfy condition $(6.4).$\\

2. From Theorem 6.4 and Proposition 0.6 we can easily obtain the following statement.\\

For an interval partition $\mathcal{F}$ to be verifying it is necessary and sufficient that for any $f\in L^1(I)$
$$E(f^*|\mathcal{F}^*)\simeq f^*.$$
\\

\textbf{Corollary 6.7.} For any interval partition $\mathcal{B}=(b_n)$ there is a finer verifying interval partition.\\

\emph{Proof.}  For a positive  $v$ let $[v]:=\textit{integral}\ \textit{part}\ v + 1,\ 0<\frac{v}{[v]}<1.$  Denote by $\vec{\beta}=(\beta_n)$ the stochastic vector to which $\mathcal{B}$ belongs. We will construct an interval partition $\mathcal{A}=(a_n)$ that belongs to the stochastic vector $\vec{\alpha}=(\alpha_n)$, and such that о $\sum_{k>n}\alpha_k=a_n,\ n\geq 1$ and then we will use the criterion $(6.4)$.\\

Let\\
$$a_0:=1,a_1:=1-\frac{\beta_1}{[\frac{\beta_1}{b_1}]},...,a_{[\frac{\beta_1}{b_1}]}:=1-[\frac{\beta_1}
{b_1}]\cdot \frac{\beta_1}{[\frac{\beta_1}{b_1}]}=1-\beta_1=b_1;\ $$
$$a_{[\frac{\beta_1}{b_1}]+1}:=b_1-\frac{\beta_2}{[\frac{\beta_2}{b_2}]},...,\ a_{[\frac{\beta_1}{b_1}]+[\frac{\beta_2}{b_2}]}:=b_1-[\frac{\beta_2}{b_2}]\cdot \frac{\beta_2}{[\frac{\beta_2}{b_2}]}:=b_1-\beta_2=b_2,...$$
$$a_{[\frac{\beta_1}{b_1}]+[\frac{\beta_2}{b_2}]+...+[\frac{\beta_{m-1}}{b_{m-1}}]+1}:=
b_{m-1}-\frac{\beta_m}{[\frac{\beta_m}{b_m}]},...,
a_{[\frac{\beta_1}{b_1}]+[\frac{\beta_2}{b_2}]+...+[\frac{\beta_{m-1}}{b_{m-1}}]+[\frac{\beta_m}{b_m}]}:=$$
$$b_{m-1}-[\frac{\beta_m}{b_m}]\cdot \frac{\beta_m}{[\frac{\beta_m}{b_m}]}=b_{m-1}-\beta_m=b_m,...$$
Because the interval partition $\mathcal{A}=(a_n)$  has points inside each interval $B_k,\ k\geq 1,$ of the partition $\mathcal{B}$, it is finer than $\mathcal{B}.$ Notice also that $0<\sup_{m\geq 1}\frac{\alpha_m}{a_m}\leq \sup_{n\geq 1}\frac{\frac{\beta_n}{[\frac{\beta_n}{b_n}]}}{b_n}\leq 1.$ Thus condition $(6.4)$ is satisfied and therefore $\mathcal{A}$ is a verifying partition.\\
 $\Box$\\

The rest of this Chapter is devoted to universal $\sigma$-subalgebras.\\

\textbf{Theorem}. (Abramovich - Lozanovsky) Let  $\mathcal{A},\ \mathcal{A}\subseteq\Lambda,$ be a $\sigma$-subalgebra of $\Lambda$.
The following conditions are equivalent.
1. $\mathcal{A}$ averages any \textit{ideal} $X\subseteq L^1(I)$.
2. $\mathcal{A}$ is a finite partition of $I$.\\

We provide the proof of this theorem in Appendix (see also~\cite{Me1}). Currently we are interested in $\sigma$-subalgebras averaging an arbitrary \textit{symmetric} ideal.\\

\textbf{Definition 6.2.}\\

1. A $\sigma$-subalgebra $\mathcal{A}$ of $\Lambda$ is called  \emph{universal} if it averages any symmetric ideal $X\subseteq L^1(I).$\\

2. A $\sigma$-subalgebra $\mathcal{A}$ in $\Lambda$ is called \emph{strongly universal} if it is universal and the averaging operator $E(\cdot |\mathcal{A})$ acts on any Banach symmetric ideal $(X,\|\cdot \|_X)$ as a contractive projection, $\|E(\cdot |\mathcal{A})\|_X\leq 1.$\\

If two $\sigma$-subalgebras are equimeasurable then they are either both universal or both not. Therefore to verify that a countable partition is universal it is enough to check it for an equimeasurable interval partition.\\

It follows from Definitions 6.1 and 6.2 that no $\sigma$-subalgebra can be at the same time universal and verifying.\\

Using Theorems 4.1 and 6.4 it is not difficult to construct an example of a countable partition that is neither universal nor verifying (see also Theorem 6.8 below).\\

It is woth noticing that if a  $\sigma$-subalgebra $\mathcal{A}$ is coarser than some universal one, then it itself is a universal $\sigma$-subalgebra.\\

Trivial examples of universal $\sigma$-subalgebras are provided by finite partitions. Independent complements of finite partitions (as follows from considered in the next Chapter formula $(7.1)$) provide examples of continuous universal $\sigma$-subalgebras.\\

\textbf{Theorem 6.8.} No countable partition $\mathcal{B}$ can be universal.\\

\emph{Proof.} If $\mathcal{B}$ is a verifying interval partition then it cannot be universal because it does not average non-majorant symmetric ideals.  Assume that $\mathcal{B}$ is a non-verifying interval  partition. We will construct a function $f=f^*\in L^1(I)$ such that $\mathcal{N}_f$ is not averaged by $\mathcal{B}$. To construct $f$ we will use reasons similar to the ones in the proof of necessity in Theorem 6.4.\\

As in that proof we can assume that  $\mathcal{B}$ is a diadic interval partition, $\mathcal{B}=(2^{-m_n})$ where (according to $(6.4)$) $ 0=m_1<m_2<...<m_n<...,\ \limsup_{n\rightarrow\infty}m_{n+1}-m_n=\infty.$\\

Let us consider the sequence of natural numbers, $\{q_n\},\ q_n=m_{n+1}-m_n,\ n \geq 1,$ and let
 $ q_{m_n}^\prime=q_{m_n}+q_{m_n+1},\ q_k^\prime=q_k\ \textrm{при}\ k\neq m_n,\ n\geq 1$. Thus, $\limsup_{n\rightarrow\infty}q^\prime_{m_n}=\infty.$

By the sequence $\{q^\prime_n\}$ we construct the non-increasing function $f$ using rules $\mathcal{U}_1)\ -\ \mathcal{U}_4)$ of Lemma 1.9. The subsequence $\{q_{m_n}^\prime\}$ is not bounded and therefore by Theorem 5.14 the interval partition $\mathcal{B}$ does not average the symmetric ideal $\mathcal{N}_f$.\\
$\Box$

\newpage
\bigskip
\centerline{ \textbf {Chapter 7. Independent complement of an interval partition.}}
\bigskip
\centerline{ \textbf {The results of this chapter are based on paper~\cite{Me1}.}}
\bigskip

 \emph{In this Chapter we study the independent complement $\mathcal{B}^\bot$ of the interval partition $\mathcal{B}$, and the action of the averaging operator $E(\cdot|\mathcal{B}^\bot).$ Our results can be applied to arbitrary countable partitions.}
\bigskip

Throughout this Chapter $\mathcal{B}=\sigma(B_n=(b_n,b_{n-1}],\ 1=b_0>b_1>...>b_n\downarrow_{n\uparrow\infty}0$ will denote an interval partition of interval $I$, $\vec{\beta}=(\beta_n),\ \beta_n=b_{n-1}-b_n,\ n\geq 1,$ will  denote the stochastic vector corresponding to $\mathcal{B}$, and $\mathcal{B}^\bot$ will denote the independent complement of $\mathcal{B}$ which always exists, ~\cite{Ro}. The constructions we use for representation $\mathcal{B}^\bot$  are analogous to the ones used in~\cite{Me1}.\\
In the sequel for any subset $C$ of the real line $\mathbb{R}$ and
for any $r,\ a\in \mathbb{R}$ we denote $a+r\cdot C:=\{a+r \cdot c,\ c\in C\}.$\\

For each $n\geq 1$ define the linear function
$\varphi_n:B_n\mapsto I$ and its inverse function
$\varphi_n^{-1}:I\mapsto B_n$ putting
$$\varphi_n(t)=\frac{t-b_n}{\beta_n},\
\varphi_n^{-1}(s)=b_n+\beta_n\cdot s.$$Now define the
$\mathcal{B}$-digit function
$\nu:I\mapsto \mathbb{N}$ and superposition $\varphi:I\mapsto I$
as follows.$$\nu(t):=\sum_{n\geq 1}n1_{B_n}(t);\ \varphi(t):=
\varphi_{\nu(t)}(t).$$ Note that both of the functions $\nu$ and
$\varphi$ are not one-to-one maps of $I$ onto $\mathbb{N}$ and
onto $I,$ respectively, namely $\nu^{-1}(n)=B_n$ for every $n,\
n\geq 1,$ whereas for the piecewise linear function $\varphi$ we
have $$\varphi^{-1}(s)=\bigcup_{n\geq
1}\{\varphi_n^{-1}(s)\}=\{b_n+ \beta_n\cdot s\}_{n\geq
1}\eqno(7.1)$$ for every $s,\ s\in I.$ The former equality
justifies notation $$\varphi^{-1}(A):=\bigcup_{n\geq
1}\{\varphi_n^{-1}(A)\}=\bigcup_{n\geq 1}\{b_n+\beta_n\cdot A\},\
A\in \Lambda.\eqno(7.2)$$ Moreover the identity
$$t=\varphi_{\nu(t)}^{-1}\Big(\varphi(t)\Big),\
t\in I,\eqno(7.3)$$is fulfilled by definition.\\

{\bf Lemma 7.1.} The identity
$$f(t)=f\Big(\varphi_{\nu(t)}^{-1}(\varphi(t))\Big)=f(b_{\nu(t)}+
\beta_{\nu(t)}\cdot\varphi(t)),\ t\in I,\eqno(7.4)$$ is fulfilled for each
$f\in L^0(I,\Lambda,\lambda).$\\

{\bf Lemma 7.2.} $\varphi$ is a measure endomorphism of the
probabilistic space $(I,\Lambda,\lambda)$ onto itself.\\

\textbf{Proof.} By formula (7.2)
\begin{center}
$\lambda\Big(\varphi^{-1}((a,b))\Big)=\lambda\Big(\bigcup_{n\geq
1}(b_n+\beta_n a,b_n+\beta_n b)\Big)=\lambda\Big(a,b)\Big)$ \end{center}
for every $a,b,\ 0< a\leq b\leq 1.$ The conclusion of Lemma 7.2 follows now by
means of standard measure-theoretic arguments.\\
$\Box$

 {\bf Theorem 7.3.}\begin{enumerate}
 \item The preimage $$\mathcal{B}^{\bot}:=\varphi^{-1}(\Lambda)$$ is a
non-atomic $\sigma$-subalgebra of $\Lambda.$ Moreover, measure
space $(I,\mathcal{B}^\bot,\lambda|_{\mathcal{B}^\bot})$ is  isomorphic to
$(I,\Lambda,\lambda).$
\item  A function $f,\ f\in L^0(I,\Lambda,\lambda),$ is
equimeasurable with function $f\circ \varphi\in
L^0(I,\mathcal{B}^\perp,\lambda|_{\mathcal{B}^\bot}).$ Moreover, on $I$ the
conditions are equivalent as follows.\\
$i).\ f\in L^1(I,\mathcal{B}^\perp,\lambda);$\\
$ii).\ f=f\circ\varphi;$\\
$iii).\ f(t)=f(b_n+\beta_n\cdot t)=f(b_m+\beta_m\cdot t),\
m,n\geq 1,\ t\in I.$
\end{enumerate}
\textbf{Proof.} Part (1) and the first
assertion in part (2) follow immediately from Lemma 7.2.  If
$ii)$ is fulfilled then the function $f$ is
$\varphi^{-1}(\Lambda)$-measurable since it is evidently true for
the function $f\circ \varphi.$ Conversely suppose that a
countable valued function $f\in L^0(I,\Lambda,\lambda)$ is
$\varphi^{-1}(\Lambda)$-measurable i.e. $f=\sum_{n\geq
1}r_n1_{C_n}(t)\in L^1(I,\mathcal{B}^{\bot},\lambda)$ where
$C_n\in \mathcal{B}^{\bot},\ n\geq 1.$ Then by part 1 of the
present theorem we have $C_n=\varphi^{-1}(D_n)$ for some $D_n\in
\Lambda,\ n\geq 1,$ hence $\varphi(C_n)=D_n,\ n\geq 1$ and
 $f(\varphi^{-1} (s))=\sum_{n\geq 1}r_n1_{D_n}(s),\ s\in I.$
 That means $ f(t)=r_n$ {\bf iff} $f\circ \varphi^{-1}(s)=r_n,\
n\geq 1,$ or equivalently $f=f\circ  \varphi.$\\
Now in the case of an arbitrary $\varphi^{-1}(\Lambda)$-measurable
function $f\in L^1(I,\Lambda,\lambda)$ it is enough to approximate $f$ uniformly by countable valued
functions.\\
We have proved that $i)\Leftrightarrow ii).$ Assume $ii)$. Then
by formula (7.4) for each $n\geq 1$ we have
$f(b_n+\beta_n\cdot t)=f(b_n+\beta_n\cdot
\varphi(t))=f(b_n+\beta_n\cdot \varphi1_{B_n}(t)=f(t)$ and
one obtains the same with $m\geq 1$ instead of $n.$ Therefore
$ii)$ implies $iii).$ Conversely if $iii)$ is fulfilled then by
formula (7.3) for any $t\in I$ we have
$(f\circ\varphi)(t)=f(\varphi^{-1}_{\nu(t)}(\varphi(t)))=f(t)$,
i.e. $iii)$ implies $ii).$ \\
$\Box$

{\bf Theorem 7.4.} $\sigma$-subalgebras $\mathcal{B}$ and
$\mathcal{B}^{\bot}$ are mutually independent and mutually
complement in $\Lambda.$\\

\textbf{Proof.} Omitting the routine details we give a sketch of proof
that
\begin{enumerate}
\item $\mathcal{B}$ and $\mathcal{B}^\perp$ are independent and
\item  Taken together $\sigma$-subalgebras $\mathcal{B}$ and
$\mathcal{B}^\perp$ generate  $\Lambda.$\end{enumerate}

(1). Each $A\in \mathcal{B}$ is obviously a countable union of
sets of the form $B_n,\ n\geq 1.$ Suppose at first $A=B_{n_0},\
n_0\geq 1,$ and take us an arbitrary $B\in \mathcal{B}^\perp.$ By
definition there exists $C,\ C\in \Lambda,$ such that
$B=\varphi^{-1}(C)=\{b_n+C\cdot \beta_n\}_{n\geq 1}$ so that $
\lambda(B_{n_0}\cap B)=\lambda (B_{n_0}\cap\{b_{n_0}+C\cdot
\beta_{n_0}\})=\beta_{n_0}\cdot\lambda(C)
=\lambda(B_{n_0})\cdot\lambda(C)= \lambda
(B_{n_0})\cdot\lambda(B).$ Now (1) follows in a routine way.\\

(2). Let $\Lambda\ni A\subset B_{n_0}$ for some $n_0\geq 1.$ Put
$B=A-b_{n_0}$ so that $A=b_{n_0}+B$ where by our assumption we have
$\lambda(B)=\lambda(A)\leq \beta_{n_0}.$ Put
$C=\frac{1}{\beta_{n_0}}\cdot B\in \Lambda.$ By definition
$A=\varphi^{-1}(C)\cap B_{n_0}\in \sigma(\mathcal{B}\cup
\mathcal{B}^\perp).$ Now (2) follows easily by passing from the
considered case to the case of an arbitrary union of the sets of
the form $\{B_n,\ n\geq 1\} $.\\
$\Box$

{\bf Theorem 7.5.} The  operator $$\sum_{n\geq
1}\beta_nf(b_n+\beta_n\cdot t),\ t\in I.\eqno (7.5)$$
maps the space $L^1(I,\Lambda,
\lambda)$ onto its subspace $L^1(I,\mathcal{B}^{\perp},\lambda)$
and is an averaged operator $E(\cdot|\mathcal{B}^{\perp})$ with respect to $\sigma$-subalgebra $\mathcal{B}^{\perp}$.\\

\textbf{Proof.}
First let us show that the two identities $f(t)=(f\circ
\varphi)(t)$ and $\sum_{n\geq
1}\beta_nf(b_n+\beta_n\cdot t)=f(t), \ t\in I,$
are equivalent. Suppose the first one is true. Then we
have$$\sum_{n\geq 1}\beta_nf(b_n+\beta_n\cdot t)=\sum_{n\geq
1}\beta_nf\Big(b_n+\beta_n\cdot \varphi (t)\Big)=$$ $$=\sum_{n\geq
1}\beta_nf\Big(\varphi_n^{-1}(\varphi (t)\Big)=\sum_{n\geq
1}\beta_nf(t)=f(t),\ t\in I,\ \eqno
(7.6) $$ where the second equality follows from formula (7.4). \\

Considering formula $(7.6)$ in the reversed order we obtain the converse implication.\\

Theorem 7.3 shows that the sum in formula (7.5) is an idempotent
linear operator acting from the space $L^1(I,\Lambda, \lambda)$
onto its subspace $L^1(I,\mathcal{B}^{\perp},\lambda|_{\mathcal{B}^\bot}).$  Moreover,
obviously the image of $\bf 1$ is $\bf 1.$ Besides using change
of variable $t=\varphi_n^{-1}(s),\ n\geq 1,$ we obtain similarly
to formula (7.6)

$$\int_0^1 \sum_{n\geq 1}\beta_nf (b_n+\beta_n\cdot
t)dt
=\sum_{n\geq 1}\beta_n^{-1}\int_0^1\beta_nf(\varphi_n^{-1}(s))ds=
\int_0^1f(t)dt.$$
Therefore $\sum_{n\geq
1}\beta_nf(b_n+\beta_n\cdot t)$ is an
averaging operator $E(f|\mathcal{B}^{\perp})$ with respect to the $\sigma$-subalgebra
$\mathcal{B}^{\perp}$ of the $\sigma$-algebra $\Lambda$,
see~\cite[Theorem 5.45, p.218]{AA}.\\
$\Box$

\textbf{Remark 7.6.} Let us list some properties of the operator $E(\cdot|\mathcal{B}^\bot)$ that follow directly from its definition.\\

1. It follows from Remark 7.2.2 that for any function $f\in L^1(I,\mathcal{B}^\bot,\lambda|_{\mathcal{B}^\bot})$ we have
$$E(f|\mathcal{B}^\bot)(t)=\sum_{i\geq 1}\beta_i\cdot f(b_i+\beta_i\cdot t)=\sum_{i\geq 1}\beta_i\cdot (f\circ\varphi\circ \varphi_i^{-1})(t)=f(t),\ t\in I.$$

Therefore, $E(L^1(I)|\mathcal{B}^\bot)=L^1(I,\mathcal{B}^\bot,\lambda|_{\mathcal{B}^\bot}).$
(We have to distinguish between $E(\cdot|\mathcal{B}^\bot)$ as an operator from $L^1(I)$ into $L^1(I)$ and $E(\cdot|\mathcal{B}^\bot)$ as an operator from $L^1(I)$ onto $L^1(I,\mathcal{B}^\bot,\lambda|_{\mathcal{B}^\bot})$).\\

2. For any $f\in L^1(I)$ it follows from the equivalence $ii)\Leftrightarrow iii)$ of Theorem 7.3 that
$$E(f\circ\varphi|\mathcal{B})=\int_0^1fd\lambda.$$
Similarly we obtain that
$$E\Big(E(f|\mathcal{B}^\bot)|\mathcal{B}\Big)=E\Big(E(f|\mathcal{B})|\mathcal{B}^\bot\Big)=\int_0^1fd\lambda.$$

3. As $E(\cdot|\mathcal{B}^\bot)$ is a double stochastic projection then for any $f\in L^1(I)$
$$E(f|\mathcal{B}^\bot)^{**}\leq f^{**}.$$
Moreover, we have $$\Big(E(f^*|{\mathcal{B}}^\bot)\Big)^{**}=E(f^{**}|{\mathcal{B}}^\bot).\eqno(7.7)$$
Indeed,
$$\Big(E(f^*|{\mathcal{B}}^\bot)\Big)^{**}(t)=\Big(\sum_{n\geq 1}\beta_nf^*(b_n+\beta_nt)\Big)^{**}(t)=\frac{1}{s}\int_0^s\sum_{n\geq 1}\beta_nf^*(b_n+\beta_nt)ds=$$
$$=\sum_{n\geq 1}\beta_n\frac{1}{s}\int_0^sf^*(b_n+\beta_nt)ds=\sum_{n\geq 1}\beta_nf^{**}(b_n+\beta_nt)   =E(f^{**}|{\mathcal{B}}^\bot)(t).$$
$\Box$\\

Averaging operators are defined only for integrable functions and therefore the right part of $(7.7)$ is undefined if the function $f^{**}$ is not integrable. In this case we understand the expression $E(f^{**}|{\mathcal{B}}^\bot)$ as the sum $\sum_{n\geq 1}\beta_nf^{**}(b_n+\beta_nt).$\\
$\Box$\\

In the sequel for a $\mathcal{W}\subseteq L^1(I)$ we denote by $E(\mathcal{W}|\mathcal{B})$ the set $\{E(w|\mathcal{B}):\ w\in \mathcal{W}\}.$\\

\textbf{Lemma 7.7.} For the operator $E(\cdot|\mathcal{B}^\bot)$ and for any $f\in L^1(I)$ and $0<a\leq 1$ we have the equality
$$E(\rho_af|\mathcal{B}^\bot)(t)=\rho_aE(f|\mathcal{B}^\bot)(t),\ t,a\in I.\eqno(7.8)$$
\emph{Proof.} By the definition of compressing-dilating operators we have
$$E(\rho_af|\mathcal{B}^\bot )(t) =\sum_{n\geq 1}\beta_n\cdot \rho_af(b_n+\beta_n\cdot t)
=\sum_{ n\geq 1}\beta_n\cdot f(b_n+\beta_n\cdot a\cdot t),\ t\in I.$$
$$\rho_aE(f|\mathcal{B}^\bot)(t)=E(f|\mathcal{B}^\bot)(a\cdot t)=\sum_{ n\geq 1}\beta_n\cdot f(b_n+\beta_n\cdot a\cdot t),\ t\in I.$$
$\Box$\\

\textbf{Theorem 7.8.} Let $f\in L^1(I),\ \mathcal{B}=\mathcal{B}^*.$ Then  $\mathcal{N}_{E(\mathcal{N}_f|\mathcal{B}^\bot)}$ is a principal symmetric ideal in  $L^1(I)$:
$$\mathcal{N}_{E(\mathcal{N}_f|\mathcal{B}^\bot)}=\mathcal{N}_{E(f^*|\mathcal{B}^\bot)}.\eqno(7.9)$$
\\
\emph{Proof.} All we need to do is to prove the inclusion $E(\mathcal{N}_f|\mathcal{B}^\bot)\subseteq\mathcal{N}_{E({f^*}|\mathcal{B}^\bot)}$,
which is equivalent to the inclusion of the left part of $(7.9)$ into its right part, because the converse inclusion in  $(7.9)$ is trivial.

Assume that $g\in E(\mathcal{N}_f|\mathcal{B}^\bot)$, i.e.
 $g=E(\tilde{f}|\mathcal{B}^\bot),\ \tilde{f}\in \mathcal{N}_f.$ It follows from $(7.6)$, from the Hardy inequality (see e.g. ~\cite[(9.1)]{CR}, and from the fact the sequence $(\beta_n)$ and the function $\tilde{f}^*$ are not monotonically increasing that
$$g(t)=\sum_{ n\geq 1}\beta_n\cdot \tilde{f}(b_n+\beta_n\cdot t)\leq$$
$$\leq\sum_{ n\geq 1}\beta_n\cdot \tilde{f}^*(b_n+\beta_n\cdot t),\ t\in I.$$
We can find $C, C>1$, and $a, 0<a<1$, such that $\tilde{f}^*\leq C\cdot \rho_af^*$ and therefore in virtue of $(7.5)$ we have
$$g^*(t)\leq C\cdot\sum_{n\geq 1}\beta_n\cdot \rho_af^*(b_n+\beta_n\cdot t)=C\cdot E(\rho_a f^*|\mathcal{B}^\bot)=C\cdot\rho_aE(f^*|\mathcal{B}^\bot)(t),\ t\in I,$$
i.e. $\Big(E(\tilde{f}|\mathcal{B}^\bot)\Big)^*\in \mathcal{N}_{E(f^*|\mathcal{B}^\bot)}.$ That proves the inclusion $\mathcal{N}_{E(\mathcal{N}_f|\mathcal{B}^\bot)}\subseteq\mathcal{N}_{E({f^*}|\mathcal{B}^\bot)}$ . \\
$\Box$\\

\textbf{Remark 7.9.} It follows from Remark 7.6.1 that $(\mathcal{B}^\bot)\mathcal{N}_{E({f^*}|\mathcal{B}^\bot)}=L^1(I,\Lambda,\lambda|_{\mathcal{B}^\bot})\cap \mathcal{N}_{E({f^*}|\mathcal{B}^\bot)}$ and a similar equality is valid for the left part of $({7.6}).$ Therefore it follows from Theorem 7.5 that $$(\mathcal{B}^\bot)\mathcal{N}_{E(\mathcal{N}_f|\mathcal{B}^\bot)}=(\mathcal{B}^\bot)\mathcal{N}_{E({f^*}|\mathcal{B}^\bot)}.$$
$\Box$\\

We will now develop some techniques that we will use in the sequel; in particular, in the proof of Theorem 7.11 below. Let $\mathcal{B}$ be an interval partition.
We construct by induction the  \textbf{$\mathcal{B}$-adic representation of points of the interval $I$.} The original interval partition, $\mathcal{B}$ we will denote $\mathcal{B}^{(0)}=\sigma(B_n^{(0)}:=(b_n^{(0)},b^{(0)}_{n-1}])$ and will call it the partition of rank 0. Its stochastic vector will be denoted  by $\beta^{(0)}:=[\beta^{(0)}_n];\ \beta^{(0)}_n:=b^{(0)}_{n-1}-b_n^{(0)}>0,\ n\geq 1;\ \sum_{n\geq 1}\beta^{(0)}_n=1$.\\
Let us fix some point $u\in (0,1]$ and denote by $B^{(0)}_{i_0}(u)$ the interval of the interval partition  $\mathcal{B}^{(0)}$ that contains $u$. We will call $i_0$ \emph{the zero digit of the $\mathcal{B}$-adic representation of point $u.$}\\
Assume that for some $k \geq 1$ the $(k-1)^{th}$ digit of $\mathcal{B}$-adic representation of point $u$ is defined, i.e. the interval partition of rank $(k-1)$, $\mathcal{B}^{(k-1)}$ is constructed. We construct the interval partition of rank $k$, $\mathcal{B}^{(k)}=\varphi(\mathcal{B}^{(k-1)}):=\sigma(B^{(k)}_{i_0,...,i_k}:=(b_{i_0,...i_{k-1},i_k},b_{i_0,...,i_{k-1},i_k-1}])$, by partitioning from the right to the left every interval of the previous rank $B_{i_0,...,i_{k-1}}^{( k-1)}=(b^{(k-1)}_{i_0,...,i_{k-1}},\ b^{(k-1)}_{i_0,...,i_{k-1}-1}]$ using the sequence of points $b^{(k)}_{i_0,...,i_{k-1},i_k}:=b^{(k)}_{i_0,...,i_{k-1},i_{k-1}-1}-\prod_{\nu=0}^k\beta_{i_\nu}$ if  $i_k>1;\ b^{(k)}_{i_0,...,i_{k-1},0}:=b^{(k-1)}_{i_0,...,i_{k-1}-1}.$

Notice that the stochastic vector of the countable partition $\mathcal{B}^{(k)}$ is $[\prod_{\nu=0}^k\beta_{i_\nu}]\in [\beta^{(0)}_n]^k$.\\
If the point $u$ is in $B_{i_0,...,i_k}^{(k)}$ then, denoting this interval by $B_{i_0,...,i_k}^{(k)}(u),$ we will call the number $i_k$ the $k^{th}$ \emph{digit} of the $\mathcal{B}$-adic representation of point $u.$
Thus we have the one-to-one correspondence between the point  $u\in I$ and the infinite sequence $u\leftrightarrow\{i_0(u),...,i_k(u),...\}$ of digits of its  $\mathcal{B}$-adic representation.\\

\textbf{Lemma 7.10.} 1). The following equalities are valid
 $$\mathcal{B}^{(k)}=\mathcal{B}\circ \varphi^{(k)},\ k\geq 0, \eqno(7.10)$$
where $\varphi^{(k)}:=\underbrace{\varphi\circ\varphi\circ\cdots\circ\varphi}_{\mbox{$k$}}.$\\

2). For any symmetric ideals $X$ and $Y$ in $L^1(I)$ and any $k\geq 1$ we have the equivalence
$$E(X\circ\varphi^{(k-1)}|\mathcal{B}^{(k-1)})\subseteq Y\Leftrightarrow E(X|\mathcal{B}^{(k-1)})\subseteq Y,$$
implying another equivalence
$$E(X|\mathcal{B})\subseteq Y\Leftrightarrow E(X|\mathcal{B}^{(k-1)})\subseteq Y,\eqno(7.11)$$
3). The $\sigma$-subalgebra $\mathbb{B}:=\sigma(\bigcup_{k\geq 0}\mathcal{B}^{(k)})$ of $\sigma$-algebra $\Lambda$, generated by all the countable partitions $\mathcal{B}^{(k)},\ k\geq 0,$ coincides with $\Lambda.$\\

\emph{Proof.} 1). Follows directly from the definitions of $\mathcal{B}^{(k)},\ k\geq 0,$ and $\varphi$.\\

 2) Because $\varphi^{(k-1)}$ is an endomorphism the symmetric ideals $X$ and $X\circ\varphi^{(k-1)}$ are $\stackrel{e}{\sim}$equivalent for $k\geq 1$. Therefore the equivalences in 2) follow from the fact that $Y$ is a symmetric ideal and Remark 0.11.1).\\

 3). Every point $u\in I$ is the intersection of decreasing by inclusion sequence of containing it intervals of rank $k$ with diameters $\prod_{\nu=0}^k\beta_{i_\nu}$, i.e., $u=\bigcap_{k=0}^\infty B_{i_0,...,i_k}^{(k)}(u)$. Let us denote by $\mathbb{R}_I$ the set of all rational points of the interval $I$ and let us fix an open interval $(c,d)\subseteq I.$ For every rational point from this interval: $r_n:\ r_n\in \mathbb{R}_I:=\{r_j\}_{j\geq 1},$ let $K(r_n)$ be the smallest rank of the intervals from the sequence $\{B_{i_0,...,i_k}^{(k)}(r_n)\}_{k\geq 1}$ that are subsets of  $(c,d)$. let $u$ be an arbitrary point from $(c,d)$ and let $u=\lim_{n(u)\rightarrow\infty} r_{n(u)}$. Clearly
$$(c,d)=\bigcup_{r_{n(u)}\in\mathbb{R}_I}B^{(K(r_{n(u)}))}_{i_0,...,i_{K(r_{n(u)})}}\in \mathbb{B}.$$
It remains to apply the standard scheme of transaltion of the Lebesgue measure $\lambda$ from open intervals to  $I.$\\
$\Box$\\

Let us fix a natural $k$. Every function $z=z(t)\in L^1(I,\mathcal{B}^{(k)},\lambda),\ k>0,$ can be written in the form $\sum_{i_0,i_1,...,i_k\geq 1}r_{i_0,i_1,...,i_k}1_{B_{i_0,i_1,...,i_k}^{(k)}}(t)$. Returning to the original notation $\mathcal{B}=\sigma(B_j=(b_j,b_{j-1}],\ b_0=1,\ \beta_j=b_{j-1}-b_j,\ j\geq 1),$ for any $t\in I$ by definition $(7.5)$ we have
$$E(z|\mathcal{B}^\bot)(t)=\sum_{j\geq 1}\beta_j\cdot z(b_j+\beta_j\cdot t)=\sum_{j\geq 1}\beta_j\cdot\sum_{i_0,i_1,...,i_k\geq 1}r_{i_0,i_1,...,i_k}\cdot1_{B_{i_0,i_1,...,i_k}^{(k)}(t)}(b_j+\beta_j\cdot t).$$
Notice that if $j\neq i_0$ then $1_{B_{i_0,i_1,...,i_k}^{(k)}(t)}(b_j+\beta_j\cdot t)=0.$ Therefore we can continue the previous equality as
$$=\sum_{i_0,i_1,...,i_k\geq 1}\beta_{i_0}\cdot r_{i_0, i_1,...,i_k}\cdot 1_{B_{i_0,i_1,...,i_k}^{(k)}(t)}(b_{i_0}+\beta_{i_0}\cdot t).\eqno(7.12) $$
Taking into consideration that $1_{B_{i_0,i_1,...,i_k}^{(k)}}(b_{i_0}+\beta_{i_0}\cdot t)=1\Leftrightarrow (b_{i_0}+\beta_{i_0}\cdot t)\in B_{i_0,i_1,...,i_k}^{(k)}$ and that $r_{i_0,i_1,...,i_k}=\frac{1}{\prod_{\nu=0}^k \beta_{i_{\nu}}}\cdot\int_{B^{(k)}_{i_0,i_1,...,i_k}(t)}zd\lambda$ we see that the sum in $(7.12)$ is equal to
$$=\sum_{i_0,i_1,...,i_k\geq 1}\frac{\beta_{i_0}}{\prod_{\nu=0}^k \beta_{i_{\nu}}}\cdot\int_{B^{(k)}_{i_0,i_1,...,i_k}(t)}zd\lambda\cdot 1_{B^{(k)}_{i_0,i_1,...,i_k}}(t)=$$
$$=\sum_{i_1,i_2,...,i_k\geq 1}\frac{1}{\prod_{\nu=1}^k \beta_{i_{\nu}}}\sum_{i_0\geq 1}\int_{B^{(k)}_{i_0,i_1,...,i_k}(t)}zd\lambda\cdot 1_{B^{(k)}_{i_0,i_1,...,i_k}(t)}(t)=$$
$$=\sum_{i_1,i_2,...,i_k\geq 1}\frac{1}{\prod_{\nu=1}^k \beta_{i_{\nu}}}\int_{_{\bigcup_{i_0\geq 1}}B^{(k)}_{i_0,i_1,...,i_k}(t)}zd\lambda\cdot 1_{B^{(k)}_{i_0,i_1,...,i_k}(t)}(t)=$$
$$=\sum_{i_1,...,i_k\geq 1}\frac{1}{\lambda(H_{i_1,...,i_k})}\int_{H_{i_1,...,i_k}}zd\lambda\cdot 1_{H_{i_1,...,i_k}}(t),$$
where the sets $H_{i_1,...,i_k}:=\bigcup_{i_0\geq 1}B_{i_0,i_1,...,i_k}$ as well as the sets $B^{(k)}_{i_0,i_1,...,i_k}$ are pairwise disjoint, and also $\lambda(H_{i_1,...,i_k})=\sum_{i_0\geq 1}\lambda(B_{i_0,i_1,...,i_k})=\sum_{i_0\geq 1}\prod_{\nu=0}^k \beta_{i_{\nu}}=\prod_{\nu=1}^k \beta_{i_{\nu}}.$ Next, because $\sum_{i_1,...,i_k\geq 1}\prod_{\nu=1}^k \beta_{i_{\nu}}=\prod_{\nu=1}^k\sum_{i_\nu\geq 1}\beta_{i_\nu}=1,$ we conclude that
$$E\Big(E(z|\mathcal{B}^{(k)})|\mathcal{B}^\bot\Big)(t)= E(z|\mathcal{H}_{i_1,...,i_k})(t),\ t\in I,\eqno(7.13)$$
where the countable partitions $\mathcal{H}_{i_1,...,i_k}=\sigma(H_{i_1,...,i_k},\ i_1,...,i_k\geq 1)$ и $\mathcal{B}^{(k-1)},$ belong to the same stochastic vector $\beta^k=[\prod_{\nu=1}^k \beta_{i_{\nu}}]$ and are equimeasurable, $\mathcal{H}_{i_1,...,i_k}\sim \mathcal{B}^{(k-1)}$. By  0.6.3) for any $x \in L^1(I)$ and for an appropriate $y \in L^1(I),\ y\sim x,$ we have
$$E\Big(E(x|\mathcal{B}^{(k)})|\mathcal{B}^\bot\Big)= E(x|\mathcal{H}_{i_1,...,i_k})\sim E(y|\mathcal{B}^{(k-1)}).\eqno(7.14)$$
Thus, for any two symmetric ideals $X, Y\subseteq L^1(I)$
$$E\Big(E(X|\mathcal{B}^{(k)})|\mathcal{B}^\bot\Big)\subseteq Y\Leftrightarrow  E(X|\mathcal{H}_{i_1,...,i_k})\subseteq Y\Leftrightarrow E(X|\mathcal{B}^{(k-1)})\subseteq Y.\eqno(7.15)$$
It follows from Paul L\'{e}vy's theorem (see~\cite[Theorem VII.4.3, p. 510]{Shi}) and from Lemma 7.10.3) that for any $x\in L^1(I)$ we have
$$\textrm{a.e.}\ \lim_{k\rightarrow\infty}E(x|\mathcal{B}^{(k)})=E(x|\mathbb{B})=E(x|\Lambda)=x,\eqno(7.16)$$
from whence follows the equality
$$\textrm{a.e.}\ \lim_{k\rightarrow\infty}E(X|\mathcal{B}^{(k)})=E(X|\mathbb{B}) =E(X|\Lambda)=X\eqno(7.17)$$
is available.\\

\textbf{\emph{Assume that for any $x\in X$ there is $y=y(x)\in L^1(I)$ such that for the sequence of countable partitions $\{\mathcal{B}^{(k)})\}_{k\geq 1}$ we have}}
$$E(x|\mathcal{B}^{(k)})\leq y,\ k\geq 1,\eqno(7.18)$$
Therefore for any pair of symmetric ideals $X$ and $Y$ we can apply the theorem of taking limits under the conditional  expectation sign~\cite{Shi},\ II.7.2 and ~\cite{KPS},\ II, 12. That justifies going to the limit a.e. when $k\uparrow\infty$ in the following chain.
$$E(X|\mathcal{B})\subseteq Y\stackrel{0.5.2)}\Leftrightarrow E(X\circ\varphi^{(k-1)}|\mathcal{B}^{(k-1)})\subseteq Y\stackrel{(7.3.2)}\Leftrightarrow E(X|\mathcal{B}^{(k-1)})\subseteq Y=$$
$$\stackrel{(7.15)}=E \Big(E(X|\mathcal{B}^{k})|\mathcal{B}^\bot\Big)\subseteq Y$$
$$\stackrel{[9],[38]} \Rightarrow
 E\Big(\textrm{a.e.}\lim_{k\rightarrow\infty}E(X|\mathcal{B}^{k})|\mathcal{B}^\bot\Big)\subseteq Y\stackrel{(7.17),(7.18)}\Rightarrow $$
$$E\Big(X|\mathcal{B}^\bot\Big)\subseteq Y.$$
We have proved the following\\

\textbf{Theorem 7.11.} Let $\mathcal{B}$ be an arbitrary  interval partition, $Y$ be an arbitrary symmetric ideal in $L^1(I),$ and $X$ be a symmetric ideal from $L^1(I)$ satisfying  condition $(7.18)$. Then
$$E\Big(X|\mathcal{B}\Big)\subseteq Y\Rightarrow E\Big(X|\mathcal{B}^\bot\Big)\subseteq Y.$$
$\Box$\\

\textbf{Corollary 7.12.} If an interval partition $\mathcal{B}$ averages a symmetric ideal $X$ satisfying $(7.18)$ then its independent complement $\mathcal{B}^\bot$ also averages $X$.\\
$\Box$\\

\textbf{Lemma 7.13.} Assume that the interval partition $\mathcal{B}=(b_n)$ averages a principal symmetric ideal $\mathcal{N}_f$. Then for any function $0<g\in\mathcal{N}_f$ the sequence $\{E(g|\mathcal{B}^{(k)})\}_{k\geq 1}$ is order bounded in $L^1(I),$ and therefore the ideal $\mathcal{N}_f$ satisfies condition $(7.18).$\\

\emph{Proof.} 1. Let us recall that for a function $f=f^*$ we use the notation $f_{\mathcal{B}}:=\sum_{n\geq 1}f(b_n)\cdot 1_{B_n}.$ It follows from Theorem 4.1  that if an interval partition $\mathcal{B}$ averages the principal symmetric ideal $\mathcal{N}_f$ then $f_{\mathcal{B}}\simeq f^{**}_{\mathcal{B}}.$ From it and from the fact that  $f^*_{\mathcal{B}}\leq f^*$ follows that $f^{**}_{\mathcal{B}}\in L^1(I)$.\\
The construction of the sequence $\mathcal{B}^{(k)}$ guarantees that for any $k\geq 1$ on the interval $\mathcal{B}^{(k)}_{i_0,...,i_k},\ i_0,...,i_k\geq 1,$ we have the inequality $E(f|\mathcal{B}^{(k)})(t)\leq f^{**}(i_0),\ t\in \mathcal{B}^{(k)}_{i_0,...,i_k}.$ From this inequality follows that $E(f|\mathcal{B}^{(k)})\leq f^{**}_{\mathcal{B}},\ k\geq 1.$\\

2. Consider now an arbitrary function $0<g\in\mathcal{N}_f,\ g^*(t)\leq [q(g)]^{-1}\cdot f^*[t\cdot q(g)],\ t\in I,$ where $q(g)<1.$ Let $\pi$ be an endomorphism of $I$ such that $g=g^*\circ \pi.$ Because $f^{**}([q(g)]^{-1}\cdot t)\leq q(g)\cdot f^{**}(t),\ t\in I,$ part 1 of the proof ensures that $E(g|\mathcal{B}^{(k)})=E(g^*\circ \pi|\mathcal{B}^{(k)})\leq q(g)\cdot f^{**}_{\mathcal{B}}\circ \pi\in L^1(I),\ k\geq 1.$\\
$\Box$\\

It follows from this Lemma and \textbf{7.11}\\

\textbf{Corollary 7.14.}  If the interval partition $\mathcal{B}$ averages the principal symmetric ideal $\mathcal{N}_f$ then its independent complement $\mathcal{B}^\bot$ also averages $\mathcal{N}_f$. If we also assume that $\mathcal{B}$ is monotonic, $\mathcal{B}=\mathcal{B}^*$, then by $(7.6)$
$$\mathcal{N}_{E(\mathcal{N}_f|\mathcal{B}^\bot)}=\mathcal{N}_{E(f^*|\mathcal{B}^\bot)}\subseteq \mathcal{N}_f.\eqno(7.19),$$
$\Box$\\

\textbf{Remark 7.15.} By Paul L\'{e}vy's theorem the condition $(7.18)$ is equivalent to $(o)$-convergence $E(x|\mathcal{B}^{(k)})\rightarrow x$ for every $x\in X$, see ~\cite{KA}.\\

\textbf{Problem.} Does the conclusion of Theorem 7.11 remain true without the assumption that for every $x\in X$ $E(x|\mathcal{B}^{(k)})\rightarrow x$, where the arrow means $(o)$-convergence?\\

\textbf{Theorem 7.16.} Assume that symmetric ideal $X$ is contained in the space  $L\log^+L(I)$ and that $Y$  is an arbitrary symmetric ideal in $L^1(I).$ If for an interval partition $\mathcal{B}$ we have the inclusion  $E(X|\mathcal{B})\subseteq Y$ then for its independent complement $\mathcal{B}^\bot$ also $E(X|\mathcal{B}^\bot)\subseteq Y$.\\

\emph{Proof.} It is well known (see e.g.~\cite{BS}, IV.6) that $L\log^+L(I)=\{x\in L^1(I):x^{**}\in L^1(I)\}.$ Consequently, for any $x\in X$ we have the condition $(7.18)$ with $y=x^{**}$.\\
$\Box$\\

\textbf{Lemma 7.17.} The projections  $E(\cdot|\mathcal{B}^\bot)$ and $E(\cdot|\mathcal{B}^{(k)})$ commute, i.e. for any $x\in L^1(I)$
$$E(E(x|\mathcal{B}^\bot)|\mathcal{B}^{(k)})=E(E(x|\mathcal{B}^{(k)})|\mathcal{B}^\bot),\ k\geq 1,$$
and therefore, in virtue of $(7.11)$
$$E(E(x|\mathcal{B}^\bot)|\mathcal{B}^{(k)})=E(E(x|\mathcal{B}^{(k)})|\mathcal{B}^\bot)\sim E(x|\mathcal{B}^{(k-1)}),\ k\geq 1.\eqno (7.20)$$
\emph{Proof.} The first equality is obvious; $(7.20)$ follows directly from it and $(7.14)$.\\
$Box$\\

\textbf{Theorem 7.18.} Let $X$ and $Y$ be arbitrary symmetric ideals in $L^1(I)$ and let $\mathcal{B}$  be an arbitrary interval partition. If $E(X|\mathcal{B}^\bot)\subseteq Y$ then $E(X|\mathcal{B})\subseteq Y$.\\

\textit{Proof}.
$$E(X|\mathcal{B}^\bot)\subseteq Y\stackrel{7.10.3)}\Rightarrow E(E(X|\mathcal{B}^\bot)|\lim_{k\uparrow\infty}\mathcal{B}^{(k)})\subseteq Y\stackrel{[38],\ (7.17)}\Rightarrow \lim_{k\uparrow\infty}E(E(X|\mathcal{B}^\bot)|\mathcal{B}^{(k)})\subseteq Y\Rightarrow$$ $$\stackrel{(7.20)}\Rightarrow \lim_{k\uparrow\infty}E(E(X|\mathcal{B}^{(k)}))|\mathcal{B}^\bot)\subseteq Y\stackrel{(7.15)}\Rightarrow \lim_{k\uparrow\infty}E(X|\mathcal{B}^{(k-1)})\subseteq Y$$
$$\stackrel{(7.11)}\Leftrightarrow \lim_{k\uparrow\infty}E(X\circ\varphi^{(k-1)}|\mathcal{B}^{(k-1)})\subseteq Y\stackrel{7.10.1)}\Leftrightarrow\lim_{k\uparrow\infty}E(X\circ\varphi^{(k-1)}|\mathcal{B}\circ\varphi^{(k-1)})\subseteq Y$$
$$\stackrel{0.5.2)}\Leftrightarrow \lim_{k\uparrow\infty}E(X|\mathcal{B})\subseteq Y\Leftrightarrow E(X|\mathcal{B})\subseteq Y.$$
$\Box$\\

\textbf{Corollary 7.19.} If the independent complement $\mathcal{B}^\bot$ of the interval partition $\mathcal{B}$ averages a symmetric ideal $X$ then $\mathcal{B}$ also averages $X$.\\
$\Box$\\

We can now complement Corollary 7.14 by the following two corollaries\\

\textbf{Corollary 7.20.} The interval partition $\mathcal{B}$ averages the principal symmetric ideal $\mathcal{N}_f$ if and only if its independent complement $\mathcal{B}^\bot$ averages $\mathcal{N}_f.$\\
$\Box$\\

\textbf{Corollary 7.21.} Let $f=f^*\in L^1(I)$ and let $\mathcal{B}=\mathcal{B}^*$ be a monotonic interval partition. Then the following conditions are equivalent.
$$ i).\ f \textrm{is a} \; \mathcal{B}-\textrm{regular function}.$$
$$ ii).\ E(f|\mathcal{B})\in \mathcal{N}_f.$$
$$ iii).\ E(f|\mathcal{B}^\bot)\in \mathcal{N}_f.$$\\
\emph{Proof.} Assume $i)$. Then by Definition 4.1. $\mathcal{B}$ averages $\mathcal{N}_f$, i.e. $E(\mathcal{N}_f |\mathcal{B})\subseteq \mathcal{N}_f$, whence follows that firstly  $E(f|\mathcal{B})\in \mathcal{N}_f,$ and secondly, by Corollary 7.20, that $E(\mathcal{N}_f |\mathcal{B}^\bot)\subseteq \mathcal{N}_f.$ It follows from Theorem  7.8 that $\mathcal{N}_{E(\mathcal{N}_f |\mathcal{B}^\bot)}=\mathcal{N}_{E(f|\mathcal{B}^\bot)},$ and therefore we have  $iii):\ E(f|\mathcal{B}^\bot)\in \mathcal{N}_f.$\\

Conversely, we have $iii)\stackrel{7.8}\Rightarrow E(\mathcal{N}_f |\mathcal{B}^\bot)\subseteq \mathcal{N}_f\stackrel{7.20}\Rightarrow E(\mathcal{N}_f |\mathcal{B})\subseteq \mathcal{N}_f. $ That by definition means that $f$ is $\mathcal{B}$-regular.\\
$\Box$\\

\textbf{Theorem 7.22.} Let $f=f^*\in L^1(I)$ and $\mathcal{B}=\mathcal{B}^*$. Then
 $$\mathcal{N}_{E(\mathcal{N}_{E(\mathcal{N}_f|\mathcal{B}^\bot)}|\mathcal{B}^\bot)}=
 \mathcal{N}_{E(\mathcal{N}_f|\mathcal{B}^\bot)}.\eqno (7.21)$$\\
 \emph{Proof.} Let us prove the inclusion
 $$\mathcal{N}_{E(\mathcal{N}_f|\mathcal{B}^\bot)}\subseteq \mathcal{N}_{E(\mathcal{N}_{E(\mathcal{N}_f|\mathcal{B}^\bot)}|\mathcal{B}^\bot)}.$$
We have $E(f|\mathcal{B}^\bot)\in E(\mathcal{N}_f|\mathcal{B}^\bot)$ and therefore, being $\mathcal{B}^\bot$-measurable, the function $E(f|\mathcal{B}^\bot)$ satisfies the relation  $$E(f|\mathcal{B}^\bot)\in\mathcal{N}_{E(\mathcal{N}_{E(\mathcal{N}_f|\mathcal{B}^\bot)}|\mathcal{B}^\bot)}.\eqno (7.22)$$
It follows now from Theorem 7.8 that
$$\mathcal{N}_{E(\mathcal{N}_f|\mathcal{B}^\bot)}=\mathcal{N}_{E(f|\mathcal{B}^\bot)}\subseteq \mathcal{N}_{E(\mathcal{N}_{E(\mathcal{N}_f|\mathcal{B}^\bot)}|\mathcal{B}^\bot)}.$$
The converse inclusion follows directly from ~\cite{Me1}, Th. 6.7.(3)\\
$\Box$\\

\textbf{Theorem 7.23.} For an arbitrary interval partition $\mathcal{B}= \sigma(B_n),\ B_n=(b_n, b_{n-1}],\ n\geq 1,\ b_0=1,$ the following conditions are equivalent.\\

1) $\mathcal{B}$ averages $\mathcal{M}^1_f$;\\

2) $\mathcal{B}^\bot$ averages $\mathcal{M}^1_f;$\\

3) $\mathcal{B}$ averages $\mathcal{N}_f$;\\

4) $\mathcal{B}^\bot$ averages $\mathcal{N}_f.$\\

\emph{Proof}. It was proved in Theorem 4.12 that if $\mathcal{B}$ averages $\mathcal{N}_f$ then $\mathcal{B}$ averages $\mathcal{M}^1_f$. By substituting $E(\cdot|\mathcal{B}^\bot)$ for $E(\cdot|\mathcal{B})$ and applying the same reasoning we can prove the implication 4)$\Rightarrow$2).

Next, by Corollary 7.19 (and via Theorem 4.12) we have 1)$\Rightarrow$3)$\Rightarrow$4), and therefore 1)$\Rightarrow$3)$\Rightarrow$4)$\Rightarrow$2). Whence, by Corollary 7.19 we have the chain of implications 1)$\Rightarrow$3)$\Rightarrow$4)$\Rightarrow$2)$\Rightarrow$1).\\
$\Box$\\
\newpage

\bigskip

\centerline{\textbf {Chapter 8. Verifying and universal complemented $\sigma$-subalgebras.}}

\bigskip

\centerline{ \textbf {The results of this chapter are based on papers~\cite{Me10},\ ~\cite{Me17}.}}

\bigskip

\emph{In this Chapter we provide a complete classification of verifying, universal, and strongly universal complemented $\sigma$-subalgebras in $\Lambda$.}\\


\textbf{Theorem 8.1.} The independent complement to a verifuing interval partition is a verifying $\sigma$-subalgebra.\\

\emph{Proof.} Let  $\mathcal{B}$ be a verifying interval partition. By definition for any symmetric ideal $X,\ X\subseteq L^1(I),$ that is not majorant we have $E(X|\mathcal{B})\nsubseteq X.$ Then, by Corollary 7.19 $E(X|\mathcal{B}^\bot)\nsubseteq X,$ and therefore  $\mathcal{B}^\bot$ also is a verifying $\sigma$-subalgebra.\\
$\Box$\\

\textbf{Theorem 8.2}. If the interval partition $\mathcal{B}=\mathcal{B}^*$ is not verifying then its independent complement $\mathcal{B}^\bot$ also is not a verifying $\sigma$-subalgebra.\\

\emph{Proof.} When we proved the necessity in Theorem 6.4, assuming that a monotonic interval partition $\mathcal{B}=(b_n)$ does not satisfy condition $(6.4)$ we constructed a function $f=f^*\in L^1(I)$ such that the symmetric ideal $\mathcal{N}_f$ is not majorant but is averaged by $\mathcal{B}$. Therefore, by Corollary 7.20 the symmetric ideal $\mathcal{N}_f$ is averaged also by $\mathcal{B}^\bot$ whence $\mathcal{B}^\bot$ is not a verifying $\sigma$-subalgebra.\\
$\Box$\\

Let us notice that for any  $\sigma$-subalgebra of mixed type (neither continuous, nor discrete) its independent complement is always continuous and, in particular, cannot be of mixed type. From Theorems 6.1, 8.1, and 8.2 follows\\

\textbf{Corollary 8.3.} In the class of complemented $\sigma$-subalgebras the following ones \textbf{and only they} are verifying.\\

1). Continuous $\sigma$-subalgebras, which have in $\Lambda$ independent continuous complement;\\

2). $\sigma$-subalgebras of mixed type, which have in $\Lambda$ independent continuous complement;\\

3). Verifying countable partitions characterised in Theorem 6.4 and their independent complements.\\
$\Box$\\

\textbf{Lemma 8.4.} Let $\mathcal{B}=(b_n)$ be a monotonic and  non-verifying countable partition. Then there is a function $f=f^*\in L^1(I)$ that is  $\mathcal{B}^{(1)}$-measurable (see \S7) but not $\mathcal{B}$-measurable, and such that $E(f|\mathcal{B})\notin \mathcal{N}_f.$\\

\emph{Proof.} By $(6.4)$ for a monotonic interval partition (i.e.  $\beta_n\downarrow_{n\uparrow \infty}0$) that is not verifying we have
 $$\sup_{n\geq 1}\frac{b_{n-1}}{b_n}=\infty,$$
We will inductively construct the subsequence $\{b_{n_k}\}$ as follows.  Let $b_{n_0}=1$. If the point $b_{n_{k-1}},\ k\geq 1,$ has been already constructed then we chose a point  $b_{n_k}$ satisfying the following conditions
$$\begin{cases}b_{n_k-1}\leq \frac{1}{2}\cdot b_{n_{k-1}};\\
\frac{b_{n_k-1}}{b_{n_k}}>b_{k+1}^{-2};\\
b_{n_k-1}\geq k\cdot s_k;\\
\end{cases}k\geq 1,\eqno(8.1)$$
where $s_0=1;\ s_k=b_k^{(n_k-1)}=b_{n_k-1}-\beta_{n_k-1}\cdot \beta_k,\ k\geq 1$ (and therefore $s_k\downarrow_{k\uparrow \infty}0$) and the countable partition $\mathcal{B}^{(1)}=\sigma(B_k^{(n)}=(b_k^{(n)},b_{k-1}^{(n)}])$ is constructed as above.\\

Define the sequence $\{u_k\}$ by letting $u_1=1,\ u_{k+1}=(2^kb_{n_k-1})^{-1},\ k>1.$ Because $(8.1)$ $b_{n_k}<s_k<b_{n_k-1},\ k>1,$ we have
$$\sum_{k=1}^\infty u_k(s_{k-1}-s_k)\leq 1-s_1+\sum_{k=1}^\infty u_{k+1}s_k\leq 1-s_1+\sum_{k=1}^\infty u_{k+1}b_{n_{k-1}-1}\leq 2-s_1<\infty.\eqno(8.2)$$
Notice that in virtue of the first inequality in $(8.1)$ we can write
$$u_{k+1}> (2^{k-1}b_{n_k-1})^{-1}\geq (2^{k-1}b_{n_k-1})^{-1}b_{n_{k-1}-1}=u_k,\ k\geq 1.\eqno(8.3)$$
If $k\geq 1$ let
$$\begin{cases}f(t)=u_2,\ \mbox{if}\  s_2<t\leq s_0;\\
f(t)=u_k,\ \mbox{if}\  s_k<t\leq s_{k-1}.\end{cases}$$
Note that $f=f^*\in L^1(I)$  is a $\mathcal{B}^{(1)}$-measurable but not $\mathcal{B}$-measurable function.\\

 It follows from our construction that for the sequence $\{s_k\}$ we have  $s_k<b_{n_k-1}<\frac{1}{2}b_{n_{k-1}}<s_{k-1},\ k\geq 1.$ By applying again the definition of $\{s_k\}$  and formulas $(8.1)$ we obtain
$$\frac{E(f|\mathcal{B})(b_{n_k-1})}{f(\frac{b_{n_k-1}}{k})}=\frac{E(f|\mathcal{B})(b_{n_k-1})}{f(b_{n_k-1})}=
\frac{1}{u_k}\cdot\frac{u_{k+1}(s_k-b_{n_k})+u_k(b_{n_k-1}-s_k)}{b_{n_k-1}-b_{n_k}}\geq $$
$$\geq \frac{u_{k+1}}{u_k}\cdot\frac{s_k-b_{n_k}}{b_{n_k-1}-b_{n_k}}=\frac{1}{2}\cdot\frac{b_{n_{k-1}-1}}{b_{n_k-1}}\cdot b_k\geq\frac{b_{n_{k-1}-1}}{b_{n_{k-1}}}\cdot b_k\geq\frac{b_k}{b_k^2}\rightarrow_{k\uparrow\infty}\infty.$$
Whence, by Lemma 0.13 we have $E(f|\mathcal{B})\notin \mathcal{N}_f.$\\
$\Box$\\

\textbf{Corollary 8.5.} For an arbitrary interval partition $\mathcal{B}$ there is a symmetric ideal not averaged by the independent complement $\mathcal{B}^\bot$.\\

\emph{Proof}. If $\mathcal{B}$ is a verifying interval partition then by Theorem 8.1 we can take as our symmetric ideal any non-majorant symmetric ideal. If on the other hand the interval partition $\mathcal{B}$ is not a verifying one then by Lemma 8.4 the symmetric ideal $\mathcal{N}_f$ is not averaged by $\mathcal{B}^\bot$ (where $f$ is the $\mathcal{B}^{(1)}$-measurable function constructed in the proof of the said lemma).\\
$\Box$\\

\textbf{Corollary 8.6.} In the class of complemented $\sigma$-subalgebras the finite partitions and their independent complements are the only universal $\sigma$-subalgebras.\\

\emph{Proof.} It is trivial that finite partitions and their independent complements are universal. By Theorem 6.8 a countable partition cannot be universal. If a $\sigma$-subalgebra is of mixed type then by Corollary 8.3 it is verifying and therefore cannot be universal. In the case when the $\sigma$-subalgebra is continuous and its complement is a countable partition we can apply Corollary 8.5.\\
$\Box$\\

Recall that a universal $\sigma$-subalgebra $\mathcal{A}$ is called  \emph{strongly universal} if for any symmetric space $(X,\|\cdot\|_X)$ the projection  $E(\cdot|\mathcal{A})$ is a \textit{contraction} on $X$, i.e. $\|E(\cdot|\mathcal{A})\|_{X\rightarrow X}\leq 1.$\\

\textbf{Remark 8.7.} Finite partitions are not only universal but \textit{strongly universal} $\sigma$-subalgebras in $\Lambda$,~\cite{MS}. We will prove below that the same is true for $\sigma$-subalgebras that are independent complements of finite partitions.\\
$\Box$\\

In the sequel we will need one of equivalent (in the sense of measure theory) representations of a continuous $\sigma$-subalgebra that is the independent complement to an at most countable partition of the interval $I$.\\

Let us fix a countable partition  $\mathcal{A}$  of the interval $I$ that belongs to the stochastic vector $\vec{\alpha}=(\alpha_n).$ We will denote the set of all atoms in $\mathcal{A}$ by $\textsf{A}:=\{a_n\}$; the set of all subsets of $\textsf{A}$ will be denoted by $\mathfrak{A}$; finally,  we define the measure $\alpha$ on $\mathfrak{A}$ by letting $\alpha(a_n)=\alpha_n,\ n\geq 1.$ The direct product of measure spaces $(\Omega,\Sigma,\mu):=(\textsf{A},\mathfrak{A},\alpha)\times(I,\Lambda,\lambda)=(\textsf{A}\times I,\mathfrak{A}\times \Lambda,\alpha\times\lambda)$ will be called  \emph{the joint realisation} of the partition $\mathcal{A}$ and its independent complement $\mathcal{A}^\bot.$ We will identify every point $a_n\in \textsf{A}$ with the set  $A_{(n)}:=a_n\times I\in \Sigma,\ n\geq 1.$ Thus, we identify at most countable partition $\mathcal{A}$ (more precisely, the $\sigma$-алгебру $\mathfrak{A}$) with the $\sigma$-subalgebra in $\Sigma$ generated by all the sets of the form $A\times I,\ A\in \mathcal{A}$; at the same time we identify the independent complement to $\mathcal{A}$ - the $\sigma$-subalgebra $\mathcal{A}^\bot$ - with the $\sigma$-subalgebra of all sets from $\Sigma$ of the form $\textsf{A}\times e,\ e\in \Lambda.$ Under this joint reralisation the averaging operator $E(\cdot|\mathcal{A}^\bot):L^1(\Omega,\Sigma,\mu)\rightarrow L^1(I,\Lambda,\lambda)$ acts as the integral of the first variable, namely:
 $$E(f|\mathcal{A}^\bot)(a,t)=\sum_{n\geq 1}\alpha_n\cdot f(a_n,t);\ f\in L^1(\Omega,\Sigma,\mu),\ (a,t)\in \Omega.\eqno(8.4).$$
 In this formula the function $E(f|\mathcal{A}^\bot)$ of two variables actually depends only on the first of them.\\

 Notice that described above joint realisation of a partition and its independent complement is unique modulo isomorphisms of measure spaces.\\

\textbf{Theorem 8.8.} For any finite partition $\mathcal{B}$ its independent complementе $\mathcal{B}^\bot$ is a strongly universal $\sigma$-subalgebra.\\

\textit{Proof.} Let $m$  be the cardinality of the partition $\mathcal{B}$. Denote by $\vec{\beta}=(\beta_k)$ the stochastic vector to which the partition $\mathcal{B}$ belongs. Consider the case when coordinates of this stochastic vector are rational numbers.
$$\beta_k=\frac{n_k}{n},\ n,n_k>0,\ k=1,...,m,\ \sum_{k=1}^mn_k=n.$$

We will use the joint realisation of the partition and its independent complement. We introduce the $n$-dimensional stochastic vector $\vec{\alpha}=(\frac{1}{n})$ and consider a belonging to it finite partiiton $\mathcal{A}$. we will consider two joint realisations: one for $\mathcal{B}$ and $\mathcal{B}^\bot$, the second for $\mathcal{A}$ and $\mathcal{A}^\bot$.
$$(\Omega,\Sigma,\mu):=(\textsf{B},\mathfrak{B},\beta)\times(I,\Lambda,\lambda);\ (\bar{\Omega},\bar{\Sigma},\bar{\mu}):=(\textsf{A},\mathfrak{A},\alpha)\times(I,\Lambda,\lambda).$$
Above $\textsf{B}=\{b_k\}_{k=1}^m$ means the set of $m$ atoms $\{b_1,b_2,...,b_m\}$; $\mathfrak{B}$ is the $\sigma$-algebra of all its subsets, and the measure $\beta$ is defined on $\mathfrak{B}$ by the equalities $\beta(b_k)=\beta_k,\ k=1,2,...,m.$ Similarly we define the triple $(\textsf{A},\mathfrak{A},\alpha)$ where $\textsf{A}=\{a_k\}_{k=1}^n$ consists of $n$ atoms, $\mathfrak{A}$ is the $\sigma$-algebra of all its subsets, and measure $\alpha$ is defined on $\mathfrak{A}$ by the equalities $\alpha(a_k)=\alpha_k=\frac{1}{n},\ k=1,2,...,n.$\\

Let $k_3:=n_1+n_2,\ k_4:=k_3+n_3,...,k_m=\sum_{i=1}^{m-1}n_i$ and let us aprtition $\textsf{A}$ into pairwise disjoint sets $\textsf{A}_1=\{a_1,...,a_{n_1}\},\ \textsf{A}_2=\{a_{n_1+1},...,a_{k_3}\},\textsf{A}_m=\{a_{k_3+1},...,a_{k_4}\}...,\ \{a_{k_m+1},...,a_n\}.$ Then we can see that  $\textsf{A}_1$ contains  $n_1$ points of the set $\textsf{A}$,  $\textsf{A}_2$ contains $n_2$ such points, ..., and $\textsf{A}_m$ contains $n_m$ points of the set $\textsf{A}.$\\

We define on $\textsf{A}$ the bijection $\varphi:\ \varphi(a_n)=a_1,...,\varphi(a_k)=a_{k+1};\ k=1,2,...,n-1.$ It is immediate that $\varphi$ is an automorphism of the measure space $(\textsf{A},\mathfrak{A},\alpha)$. This automorphism generates the map $\bar{\varphi}:\bar{\Omega}\rightarrow\bar{\Omega}$ defined by the formula $\bar{\varphi}(a,t):=(\varphi(a),t),\ (a,t)\in \bar{\Omega}$. this map in its turn is an automorphism of the measure space $(\bar{\Omega},\bar{\Sigma},\bar{\mu}).$
The linear isometry $\bar{\Phi}$ of the space  $L^1(\bar{\Omega},\bar{\Sigma},\bar{\mu})$ is generated by the automorphism $\bar{\varphi}$ according to the formula
$$\bar{\Phi}\bar{f}(a,t)=\bar{f}(\bar{\varphi}(a,t),\ \bar{f}\in L^1(\bar{\Omega},\bar{\Sigma},\bar{\mu}),\ (a,t)\in \bar{\Omega}.$$
Therefore the functions $\bar{f}$ and $\bar{\Phi}\bar{f}$ are equimeasurable.\\

Let $(X,\|\cdot\|_X)$ be a symmetric space that is a subset of $L^1(\Omega,\Sigma,\mu)$ and let $(\bar{X},\|\cdot\|_{\bar{X}})\  \stackrel{e}{\sim}$ be an equivalent to it symmetric space that is a subset of $L^1(\bar{\Omega},\bar{\Sigma},\bar{\mu})$ (see 0.10). Let us fix an arbitrary element $x\in X,\ \|x\|_X\leq 1$ and define $\bar{x}\in L^1(\bar{\Omega},\bar{\Sigma},\bar{\mu})$ by the formula
$$\bar{x}(a,t)=x(b_i,t),\ \textrm{если}\ a\in\textsf{A}_i,\ i=1,...,m;\ t\in I.$$
Clearly $\bar{x}\sim x$ whence $\bar{x}\in \bar{X},\ \|\bar{x}\|_{\bar{X}}\leq 1.$ Consider
$$\bar{y}:=\frac{1}{n}\sum_{i=1}^{n-1}\bar{\Phi}^i\bar{x}.$$
It follows from the previous text that $\bar{y}\in \bar{X}$ and
$$\bar{y}(a,t)=\frac{1}{n}\sum_{i=1}^{n-1}\bar{x}(a,t)=\frac{1}{n}\sum_{i=1}^m\sum_{a\in \mathcal{A}_i}\bar{x}(a,t)=\sum_{i=1}^m\frac{n_i}{n}x(b,t),\ (a,t)\in \bar{\Omega}.$$
Thus, $\bar{y}\sim E(x|\mathcal{B}^\bot)$ and therefore $\|\bar{y}\|_{\bar{X}}=\|E(x|\mathcal{B}^\bot)\|_X.$\\

On the other hand
$$\|\bar{y}\|_{\bar{X}}\leq \frac{1}{n}\sum_{i=1}^{n-1}\|\bar{\Phi}^i\bar{x}\|_{\bar{X}}=
\frac{1}{n}\sum_{i=1}^{n-1}\|\bar{x}\|_{\bar{X}}=\|\bar{x}\|_{\bar{X}}=\|x\|_X\leq 1.$$
Thus, $\|E(x|\mathcal{B}^\bot)\|_X\leq 1$ and we are done in the case of rational coordinates of the stochastic vector  $\vec{\beta}.$ The general case now can be obtained by using the routine approximation arguments.\\
$\Box$\\

\textbf{Corollary 8.9.} In the class of complemented $\sigma$-subalgebras only finite partitions and their complements are strongly universal.\\
$\Box$\\

\newpage

\centerline{\textbf{APPENDIX}}
\bigskip
Many years ago G. Ya. Lozanovsky and Yu. A. Abramovich stated and solved the problem about averaging of an arbitrary (not necessarily symmetric) ideal $X,\ X\subseteq L_1(\Omega,\Sigma,\mu)$ where $(\Omega,\Sigma,\mu)$ is a probability space with a continuous measure $\mu$. More precisely their questions were as follows.

\textbf{1. What properties of a $\sigma$-subalgebra $\mathcal{A}$ guarantee that it averages an arbitrary ideal?}

\textbf{2. What Banach and/or order properties of an ideal $X$ guarantee that it is averaged by any $\sigma$-subalgebra of $\Sigma$?}\\

The answers to these questions obtained, respectively, by Lozanovsky and Abramovich were never published by them and first appeared in the author's paper~\cite{Me1}. We reproduce these answers in the following theorem that, in particular, asserts that even a very "good" Banach ideal (a $KB$ space) is not averaged by some $\sigma$-subalgebras. Moreover, if a $\sigma$-subalgebra is not generated by a finite set of atoms (thus, in case of ($I, \Lambda, \lambda$), if it is not a finite partition)  then there is an order ideal in $L_1(\Omega,\Sigma,\mu)$ that is not averaged by it.\\

\textbf{Theorem A} (Abramovich)  If the $\sigma$-subalgebra $\mathcal{A}$ contains a countable family of atoms, i.e. $\mathcal{A} \supseteq \mathcal{F} = \sigma(F_n), n \in \mathds{N}$, then there exists an ideal $X$ in $L_1(\Omega,\Sigma,\mu)$
such that $\mathcal {F}$ does not average $X.$
Moreover, $X$ endowed with an appropriate lattice norm is a $KB$-space~(see e.g.~\cite{KA}).\\

\textit{Proof}. Denote $\mu_n=\mu(F_n),\ n\geq 1.$ Since
$\sum_{n=1}^{\infty}\mu_n=1$ there exists a subsequence
$0<\mu_{n_k}\leq 1/k^3,\ k\geq 1.$ Because $\mu$ is a continuous measure it is possible to
divide each of the sets $F_{n_k}$ into two disjoint subsets
$F_{n_k}^1$ and  $F_{n_k}^2$ such that
$\mu(F_{n_k}^1)=\mu(F_{n_k}^2)=\mu_{n_k}/2,\ k\geq 1.$ \footnote{This is a well known fact that follows easily e.g. from the Saks' lemma,see~\cite[Lemma IV.9.3.7, p.808]{DS}. See also~\cite{Ro}.} Put
$F^1=\bigcup_{k=1}^{\infty}F_{n_k}^1,\ F^2=\Omega\setminus F^1,$
$$X=\{x\in L_1(\Omega,\Sigma,\mu){\bf |}\ ||x||_X:=\int_{
F^1}|x|d\mu+\int_{ F^2}x^2d\mu<\infty.\}$$Evidently $X$ is an
ideal in $L^1(\Omega,\Sigma,\mu).$ Let us define on $\Omega$ a
measurable function $x_0$ as follows
$$x_0(\omega)=\left\{\begin{array}{lll}\frac{1}{\sqrt{\mu_{n_k}}},\
{\rm if}\ \omega\in F_{n_k}^1,\ k\geq 1;\\0,\qquad {\rm if}\ \omega\in F^2.\\
\end{array}\right. $$ So we have $\int_{F^2}x_0^2d\mu=0,\ \int_{F^1}x_0d\mu=
\sum_{k=1}^{\infty}\frac{\mu_{n_k}}{2\sqrt{\mu_{n_k}}}<
\sum_{k=1}^{\infty}\frac{1}{k\sqrt{k}}<\infty$ therefore $x_0\in
X.$

On the other hand, $E(x_0|\mathcal
{F})=\sum_{k=1}^{\infty}r_n1_{F_n}$, where
$$r_n=\left\{\begin{array}{lll}\frac{\mu_{n_k}}{2\mu_{n_k}\sqrt{\mu_{n_k}}}=\frac{1}{2\sqrt{\mu_{n_k}}},\
k\geq 1,\ {\rm if}\ n\in\{n_k\}_{k=1}^{\infty};\\0,\qquad {\rm if}\ n\notin \{n_k\}_{k=1}^{\infty}. \\
\end{array}\right. $$ Therefore $\int_{F^2}[E(x_0|\mathcal {F})]^2d\mu=1/2\sum_{k=1}^{\infty} r_{n_k}^2\mu_{n_k}=
1/2\sum_{k=1}^{\infty}\frac{\mu_{n_k}}{4\mu_{n_k}}=\infty$ so
$E(x_0|\mathcal {F})\notin X,$ thus c.p. $\mathcal {F}$ do not
average ideal $X.$\\
$\Box$

\textbf{Remark.} Theorem A follows, of course, from Proposition 3.1 and thus currently is of historic interest. Moreover, we see from Proposition 3.1 that as soon as an ideal $X$ is not symmetric, it can be otherwise arbitrary ``nice", e.g. reflexive, or even isometrically and order isomorphic to a some $L^2$-space, but there is a countable partition that does not average it.\\

\textbf{Theorem B}. (Lozanovsky)  In order for a $\sigma$-subalgebra $\mathcal
{A}$ of the $\sigma$-algebra $\Sigma$ to average every ideal $X$ in
$L_1(\Omega,\Sigma,\mu)$ it is necessary and sufficient that
$\mathcal {A}$ contained at most finite number  of $\mathcal
{A}$-atoms $\{e_k\}_{k=1}^n$ and therefore $\mathcal {A}$ must be of the form
$\mathcal {A}=\sigma(e_1,\cdots,e_n,\Sigma(e))$ where
$e=\Omega\setminus\bigcup_{k=1}^ne_k$ while $\Sigma(e)$ denotes
the trace of $\sigma$-algebra $\Sigma$ on the set $e$ i.e.
$$\Sigma(e)=\{S\cap e:S\in \Sigma\}.$$
\textit{Proof.} Sufficiency follows immediately from the definition of an ideal
and the basic properties of conditional expectation
operator. \\

Let $\mathcal {A}$ be any $\sigma$-subalgebra of $\Sigma$ which
averages every ideal $X\subseteq L_1(\Omega,\Sigma,\mu).$
According to Theorem A $\mathcal {A}$ can contain only finite number of
atoms $e_1,\ \cdots,\ e_n;$ denote
$e=\Omega\setminus\bigcup_{k=1}^ne_k\in \mathcal {A}.$ It is
clear that the trace $\mathcal {A}(e)$ is a continuous
$\sigma$-subalgebra of the trace $\sigma$-algebra $\Sigma(e).$\\

We have to prove the equality $\mathcal {A}(e)=\Sigma(e)$ or,
equivalently, the inclusion $$\Sigma(e)\subseteq \mathcal {A}.
\eqno (*)$$The proof arises
from the sequence of auxiliary Lemmas 1 - 4 as follows.\\

{\bf LEMMA 1.} Let $F\in \mathcal {A},\ F\subseteq e$ and
$\Sigma(F)\not\subseteq \mathcal {A}.$ Then there are two mutually
disjoint subsets $F_1,\ F_2\in \mathcal {A}$ such that $F_1\cup
F_2=F$ and $\Sigma(F_1)\not\subseteq \mathcal {A},\
\Sigma(F_2)\not\subseteq \mathcal {A}.$\\

\textbf{Proof.} Assume the contrary. Using ~\cite{DS}, v.I,IV,9, let us take any two $F_1^1,\
F_2^1\in \mathcal {A}$ such that $\mu(F_1^1\cap F_2^1)=0,\
F_1^1\cup F_2^1=F,\ \mu(F_1^1)=\mu(F_2^1)=\frac{1}{2}\mu(F).$ For
definiteness we can assume that $\Sigma(F_1^1)\subseteq \mathcal
{A},\ \Sigma(F_2^1)\not\subseteq \mathcal {A}.$ Continuing in the
same way we can construct two sequences of subsets $\{F_i^j,\
j=1,2,\cdots,\ ;i=1,2,\}$ such that the following conditions are satisfied
$$\mu(F_1^j\cap F_2^j)=0,\ F_1^{j+1}\cup
F_2^{j+1}=F_2^j,\ \Sigma(F_1^j)\subset \mathcal {A},\
\Sigma(F_2^j)\not\subseteq \mathcal {A},\
\mu(F_i^j)=\frac{1}{2^j}\mu(F).$$It follows from the last
equalities that $\mu(\bigcap_{j=1}^{\infty}F_2^j)=0$ therefore
$\mu$-a.e. equality $E=\bigcup_{j=1}^{\infty}(E\cap F_1^j)$ holds
for every $E\subseteq F$ with $\mu(E)>0.$ On the other hand by our
construction $E\cap F_1^j\in \mathcal {A},\ j=1,2,\ \cdots,$
therefore $E\in \mathcal {A}.$\\
Thus we prove that $\Sigma(F)\subseteq \mathcal {A},$ - a
contradiction.\\
$\Box$

For each $E\in \Sigma$ denote $$\hat{E}=\bigcap_{F\in \mathcal
{A},\ F\supseteq E}F;\ \check{E}=\bigcup_{F\in \mathcal {A},\
F\subseteq
E}F.$$Immediately from these definitions follows\\

{\bf LEMMA 2.}  Both of the sets $\hat{E}\setminus E$ and
$E\setminus \check{E}$ do not contain any $\mathcal
{A}$-measurable subset of a positive measure.\\

{\bf LEMMA 3.} {\it Let $F\in \mathcal {A}$ be given such that
$F\subseteq e,\ \Sigma(F)\not\subseteq \mathcal {A}.$ Then there
exists a set $E\in \Sigma(F)$ such that $E\notin \mathcal {A},\
\mu(\check{E})=0$ and $\mu(E)\geq \frac{1}{2}\mu(\hat{E}).$}\\
\textbf{Proof.} Fix an arbitrary $G\in \Sigma(F)$ such that $G\notin {\mathcal
A}$ and put $E=G\setminus \check{G},\ E_1=\hat{E}\setminus E;$
obviously  $E\in \Sigma(F).$ By Lemma 3.3 we have $E\notin
\mathcal {A}$ and moreover $\mu(\check{E})=0.$ Analogously
$E_1\in \Sigma(F),\ E_1\notin \mathcal {A}$ and
$\mu(\check{E}_1)=0.$ Furthermore at least one of the two
inequalities $\mu(E)\geq \frac{1}{2}\mu(\hat{E}),\ \mu(E_1)\geq
\frac{1}{2}\mu(\hat{E})$ has to be true. If it is the first one
then $E$ is the desired set. In the second case inclusion
$\hat{E_1}\subset \hat{E}$ implies inequality $\mu(E_1)\geq
\frac{1}{2}\mu(\hat{E_1})$ so $E_1 $ is the desired set.
$\Box$

Note that immediately from definition of averaging operator
follows that $ supp\ E(x|\mathcal {A})\subseteq F$ provided that
$ supp\
x\subseteq F\in \mathcal {A}.$ Our next auxiliary assertion is\\

{\bf LEMMA 4.} Let $B\in \Sigma(F)$ satisfy the following conditions.
$$\mu(B)>0,\ \mu(\check{B})=0,\ \mu(B)\geq
\frac{1}{2}\mu(\hat{B}).$$ Put $x=1_B.$ Then
$vraisup_{\hat{B}\setminus B}E(x|\mathcal {A})\geq 1/2.$\\
\textbf{Proof.}
Denote  $S=supp\ E(x|\mathcal {A}).$ Clear that $B\subset S
\subseteq \hat{B},\ S\in \mathcal {A}.$ Suppose that
$$p:=vraisup_{\hat{B}\setminus B}E(x|\mathcal {A})<1/2$$ and put
$G=\{\omega\in \hat{B}:\ E(x|\mathcal {A})>p \}.$ Evidently
$G\subseteq B,\ G\in \mathcal {A},$ so by Lemma 3.3 and in view of
$\mu(\check{B})=0$ we have $\mu(G)=0.$ That means that inequality
$E(x|\mathcal {A})(\omega)\leq p$ is true a.e. on $\hat{B}.$
Therefore $$\int_{\hat{B}}E(x|\mathcal {A})d\mu\leq
p\int_{\hat{B}}d\mu=p\mu(\hat{B})<1/2\cdot 2\mu(B)\leq \mu(B).$$
On the other hand by definition of conditional expectation we have
$$\int_{\hat{B}}E(x|\mathcal {A})d\mu=\int_SE(x|\mathcal
{A})d\mu=\int_Bxd\mu=\mu(B),$$ - a contradiction.
$\Box$

Let us now proceed with the proof of (*). If it is false then
according to Lemma 3.2 there is a sequence of the mutually
disjoint sets $F_n\in\Sigma(e)$ such that $F_n\in \mathcal {A},\
\Sigma(F_n)\not\subseteq \mathcal {A},\ n\geq 1.$ Choosing an appropriate subsequence we can assume the inequalities $\mu(F_n)\leq n^{-3},\
n\geq 1,$. In view of Lemma 3.4 there are sets
$B_n\in\Sigma(F_n)$ such that $$B_n\notin \mathcal {A},\
\mu(\check{B_n})=0,\ \mu(B_n)\geq \frac{\mu(\hat{B_n})}{2},\
n\geq 1.$$ Put $x_n=\frac{1}{\mu(B_n)}1_{B_n},\ n\geq 1.$
Functions $x_n$ are evidently mutually disjoint thus functions
$E(x_n|\mathcal {A})$ are mutually disjoint, too, since $supp\
E(x_n|\mathcal {A})\subseteq F_n,\ n\geq 1.$ Put $$x_0={\bf
1}+\sum_{n=1}^{\infty}\frac{x_n}{n^2}.$$Then $x_0\in
L_1(\Omega,\Sigma,\mu)$ since all the functions $x_n$ have unit
norms in $L_1(\Omega,\Sigma,\mu).$ Denote $$X_0=\{x\in
L_1(\Omega,\Sigma,\mu){\bf |}\ \exists\ r=r(x)\ {\rm such\ that}\
|x|\leq r\cdot x_0 \}.$$ It is clear that $X_0$ is an ideal in
$L_1(\Omega,\Sigma,\mu).$ According our supposition
$E(X_0|\mathcal {A})\subseteq X_0$ therefore $E(x_0|\mathcal
{A})\leq r_0\cdot x_0$ for some $r_0>0.$ Keeping in mind
properties of conditional expectation operator and the mutual
disjointness of elements $x_n$ we come from the last inequality
to the system of inequalities
$$1_{\hat{B_n}}+\frac{E(x_n|\mathcal {A})}{n^2}\cdot
1_{\hat{B_n}}\leq r_0\Big(1_{\hat{B_n}}+\frac
{x_n}{n^2}\Big)\cdot 1_{\hat{B_n}},\ n\geq 1.$$From here via
projecting onto $\hat{B}_n\setminus B_n$ we obtain according to
Lemma 3.5 $$\Big(
1+\frac{1}{2n^2\mu(B_n)}\Big)1_{\hat{B_n}\setminus
B_n}\leq\Big(1+\frac{E(x_n| \mathcal
{A})}{n^2}\Big)1_{\hat{B_n}\setminus B_n}\leq r_0\cdot
1_{\hat{B_n}\setminus B_n}$$ where $\mu(\hat{B_n}\setminus
B_n)>0,\ n\geq 1.$ Since $\frac{1}{\mu(B_n)}\geq
\frac{1}{\mu(F_n)}\geq n^3,\ n\geq 1,$ we come to a contradiction:
$1+\frac{n}{2}\leq r_0,\ n\geq 1.$\\
$\Box$

\newpage
\centerline{\textbf{INDEX}}

$1_A$ \ldots 8

$\mathbf{1}$ \ldots 8

$(\mathcal{A})\mathcal{M}_X$ \ldots 17

$(\mathcal{A})\mathcal{N}_X$ \ldots 17

$f^\star$ \ldots 8

$f^{\star \star}$ \ldots 8

$\mathcal{F}$-regular function \ldots 55

$I$ \ldots 8

$(I, \Lambda, \lambda)$ \ldots 8

$L^0(I, \Lambda, \lambda) = L^0(I), L^1(I), L^\infty(I)$ \ldots 8

$\mathcal{M}_f$ \ldots 17

$\mathcal{N}_f$ \ldots 17

$\mathcal{M}_Z$ \ldots 17

$\mathcal{N}_Z$ \ldots 17

$\mathcal{M}^1_f$ \ldots 18

$y \prec x$ \ldots 8

$\lambda$ \ldots 8

$\Lambda$ \ldots 5 \ldots 8

$\sigma$-subalgebra \ldots 9

$\sigma$-subalgebra $\mathcal{A}$ averages ideal $X$ \ldots 15

$\sigma$-subalgebra $\mathcal{A}$ strongly averages symmetric space $X$ \ldots 15

 $\sigma(\Gamma)=\sigma$-subalgebra generated by the family $\{G_{\gamma}\}_{\gamma\in \Gamma}$

  of subsets of $\Lambda$ \ldots 15

 $\stackrel{e}{\sim}$
 equivalent sets \ldots 16

  $\stackrel{e}{\sim}$equivalent $\sigma$-subalgebras \ldots 16

atom in a $\sigma$-subalgebra \ldots 9

automorphism \ldots 8

averaging operator \ldots 12

coarser of two $\sigma$-subalgebras \ldots 10

coarser of two stochastic vectors \ldots 10

conditional mathematical expectation \ldots 12

compression-dilation operator, $\rho_s$, \ldots 8

continuous $\sigma$-subalgebra \ldots 9

countable partition \ldots 10

countable partition belonging to a stochastic vector \ldots 10

diadic point \ldots 11

diadic interval partition \ldots 11

discrete $\sigma$-subalgebra \ldots 9

elementary function \ldots 9

endomorphism \ldots 8

equimeasurable functions \ldots 8

equimeasurable $\sigma$-subalgebras \ldots 10

equivalent functions \ldots 8

equivalent interval partitions \ldots 10

finer of two $\sigma$-subalgebras \ldots 10

finer of two stochastic vectors \ldots 10

finite partition \ldots 77

golden ratio \ldots 35

ideal $X$ is averaged by $\sigma$-subalgebra $\mathcal{A}$ \ldots 15

independent complement of a $\sigma$-subalgebra \ldots 10

interpolation vector space between $L^1(I)$ and $L^\infty(I)$ \ldots 13

interval partition \ldots 10

interval partition belonging to a stochastic vector $\phi$ \ldots 10

majorant ideal \ldots 14

Marcinkiewicz spaces $M_f$ and $M_f^1$ \ldots 18

mixed type $\sigma$-subalgebra \ldots 10

monotonic interval partition \ldots 11

monotonic rearrangement of interval partition \ldots 35

monotonic stochastic vector $\vec{\beta}^*=[\beta^*_n]$ \ldots 11

multiple of interval partition \ldots 11

non-majorant ideal \ldots 14

operator of conditional mathematical expectation \ldots 12

orbit, $\Omega_f$, of function $f$ \ldots 13

order ideal in $L^1(I)$ \ldots 14

regular function \ldots 20

right diadic projection, $\mathcal{F}^\star_{(2)}$ of a countable partition $\mathcal{F}$ \ldots 11

Ryff's orbit, $\mathfrak{D}_x$ \ldots 13

sample from a countable partition \ldots 11

stochastic vector \ldots 10

strongly interpolation Banach space between $L^1(I)$ and $L^\infty(I)$ \ldots 13

strongly majorant ideal \ldots 14

strongly majorant symmetric space \ldots 14

strongly universal $\sigma$-subalgebra \ldots 77

sum of a family of ideals \ldots 15

symmetric ideal \ldots 14

symmetric space \ldots 14

universal $\sigma$-subalgebra \ldots 77

verifying $\sigma$-subalgebra \ldots 72

weakly regular function \ldots 20

\end{document}